%% file: Manuscript.tex
\documentclass[11pt,letterpaper,twoside,reqno,nosumlimits]{amsart}

\allowdisplaybreaks

\input{macro-3}

\graphicspath{{Figures}}

\usepackage[numbers]{natbib}

\newcommand{\om}{\omega}
\newcommand{\Om}{\Omega}

\newcommand{\la}{\langle}
\newcommand{\ra}{\rangle}

\newcommand{\cX}{\mathcal X}
\newcommand{\pv}{\operatorname{P.V.}}

\mathtoolsset{showonlyrefs}

\begin{document}

\title[Asymptotically self-similar blowups of the gCLM model and the HL model]
{Asymptotically Self-similar finite-time blowups of \\the generalized Constantin--Lax--Majda model and\\ the 1D Hou--Luo model}
\author[D. Huang, X. Li]{ De Huang$^1$, Xiangyuan Li$^2$}
\thanks{$^1$School of Mathematical Sciences, Peking University. E-mail: dhuang@math.pku.edu.cn}
\thanks{$^2$School of Mathematical Sciences, Peking University. E-mail: lixiangyuan23@stu.pku.edu.cn}

\begin{abstract}
   We establish asymptotically self-similar finite-time blowups of the $a$-parameterized family of the generalized Constantin--Lax--Majda model from $C_c^{\infty}$ initial data for all $a\leq 1$. Our approach is based on a finite-codimensional stability analysis of the dynamic rescaling equation around the exact self-similar blowup profiles constructed in \cite{huang2024self}. Following the same framework, we also give a purely analytic proof of the asymptotically self-similar finite-time blowup of the 1D Hou--Luo model from $C_c^{\infty}$ initial data based on the exact self-similar blowup profile constructed in \cite{huang2025exact}.
\end{abstract}

\maketitle
\section{Introduction}
\subsection{Asymptotically self-similar blowups of the generalized Constantin--Lax--Majda model}
We consider the $1$D generalized Constantin--Lax--Majda (gCLM) equation
\begin{equation}\label{eqt:gCLM-model}
\om_t + au\om_x = u_x \om, \quad u_x = \mtx{H}(\om),
\end{equation}
for $x\in \R$, where $\mtx{H}(\cdot)$ denotes the Hilbert transform on the real line. This equation is a $1$D model for the vorticity formulation of the $3$D incompressible Euler equations, proposed to study the competitive relation between advection and vortex stretching. In particular, $\om$ models the vorticity, and the nonlinear terms $u\om_x$ and $u_x \om$ model the advection term and the vortex stretching term, respectively. The $3$D Biot-Savart law that recovers the velocity from the vorticity is modeled by $u_x= \mtx{H}(\om)$, which has the same scaling as the original Biot-Savart law.

Whether the 3D incompressible Euler equations can develop a singularity in finite time is one of the most fundamental problems in mathematical fluid dynamics. In recent years, substantial progress has been made on this problem. Let us briefly review some of the key developments below. In the axisymmetric no-swirl case, Elgindi \cite{elgindi2021finite,elgindi2022stability} and C{\'o}rdoba, Mart\'inez-Zoroa, and Zheng \cite{cordoba2025finite} constructed finite-time blowups (of vorticity)  from $C^{1,\alpha}$ initial velocity for sufficiently small $\alpha>0$. Elgindi's construction is self-similar (the definition of ``self-similar'' will be explained later), whereas the construction of C{\'o}rdoba, Mart\'inez-Zoroa, and Zheng builds on a cascading of vortices and therefore is not self-similar. Later, independent works by Shkoller \cite{shkoller2026threshold}, Chen \cite{chen2026eulerII}, and Shao, Wei, Zhang, and Zhang \cite{shao2026selfsimilar} obtained three different scenarios of finite-time blowup (with the latter two being typically self-similar) in the axisymmetric no-swirl setting with $C^{1,\alpha}\cap L^2$ initial velocity for $\alpha$ arbitrarily close to the regularity threshold $1/3$. Indeed, finite-time blowup has been excluded for $\alpha\geq 1/3$ (e.g. see \cite{shao2026global}). In the case of smooth initial data, Chen and Hou \cite{chen2022stable,chen2025singularity} first proved with computer-assistance that the 3D axisymmetric Euler equations can develop finite-time, nearly self-similar blowup. Their scenario involves a solid boundary, and the blowup occurs on the boundary. Astonishingly, OpenAI \cite{openai2026euler} recently claimed the existence of finite-time blowup from $C_c^{\infty}$ initial velocity for the 3D incompressible Euler equations on the whole $\R^3$. Their construction builds on a multiscale-amplification mechanism, an idea that was previously developed in the constructions of blowup by C{\'o}rdoba, Mart\'inez-Zoroa, and Zheng \cite{cordoba2025finite} and by C{\'o}rdoba and Mart{\'i}nez-Zoroa \cite{cordoba2023forced}. Nevertheless, the physical relevance of such a construction remains unclear. Therefore, it is still worthwhile to investigate potential finite-time blowup solutions to the 3D incompressible Euler equations that have clearer structural features such as self-similarity.

As a simplified model of the 3D incompressible Euler equations, the gCLM equation \eqref{eqt:gCLM-model} has also been extensively studied. This family includes several classical models corresponding to different values of the parameter $a$.
The original version of \eqref{eqt:gCLM-model} with $a=0$ was proposed by Constantin, Lax and Majda \cite{constantin1985simple} to demonstrate that a nonlocal vortex stretching term can lead to finite-time blowup in the absence of advection. Later, De Gregorio \cite{de1990one} included an advection term $u\omega_x$ back in the equation (known as the De Gregorio model, corresponding to $a=1$) and conjectured the occurrence of a finite-time singularity nonetheless. As a generalization, Okamoto, Sakajo and Wunsch \cite{okamoto2008generalization} introduced the real parameter $a$ to modify the effect of advection in the competition against vortex stretching. Hence, equation \eqref{eqt:gCLM-model} is also referred to as the Okamoto--Sakajo--Wunsch (OSW) model.

There have been numerous results on the finite-time blowup of the gCLM model. We refer the readers to \cite{huang2024self} for a comprehensive survey. Among these blowup results, we are particularly interested in self-similar ones, in view of the scaling property of equation \eqref{eqt:gCLM-model}. More specifically, we say \eqref{eqt:gCLM-model} exhibits an exact self-similar blowup if it admits a solution $\om$ that takes the form
\begin{equation}\label{eqt:gCLM-self-similar-solution}
\om(x,t) = (T-t)^{\lambda}\bar{\Omega}\left(\frac{x}{(T-t)^{\gamma}}\right),
\end{equation}
where $\bar{\Omega}$ is the self-similar profile and $\lambda,\gamma$ are the scaling exponents. Substituting this ansatz into \eqref{eqt:gCLM-model} gives $\lambda=-1$. The value of $\gamma$ determines the spatial feature of the blowup: A positive $\gamma$ corresponds to a spatially focusing blowup at $x = 0$, while a negative $\gamma$ corresponds to an expanding blowup. In contrast to an exact self-similar solution, an asymptotically self-similar finite-time blowup refers to a solution that exhibits clear self-similarity only as $t$ approaches the finite blowup time $T$:   
\begin{equation}\label{eqt:gCLM-asymptotic-selfsimilar}
    \om(x,t)=(T-t)^{\lambda}\Omega\left(\frac{x}{(T-t)^{\gamma}},t\right),
\end{equation}
where $\lambda=-1$, and the time-dependent profile $\Omega(\cdot,t)$  will converge to a rescaled version of the self-similar profile $\bar{\Omega}$ as $t\to T$. 

There are two main motivations for studying asymptotically self-similar blowups. First, a physically relevant blowup mechanism should be robust under small perturbations of the initial data. Such perturbations generally destroy exact self-similarity, so one seeks persistence of the same asymptotic blowup behavior after rescaling. Second, in higher-dimensional fluid equations, exact self-similar solutions often have slowly decaying tails, which leads to infinite kinetic energy. To obtain finite-energy solutions, one therefore seeks to truncate these tails while preserving the local blowup dynamics. This truncation also destroys exact self-similarity and naturally leads to the study of asymptotically self-similar blowups \cite{elgindi2022stability,chen2026eulerII,shao2026selfsimilar}.

A well-established framework for studying asymptotically self-similar blowup, called the dynamic rescaling method, has been widely used in the literature. Roughly speaking, this method uses a suitable time-dependent change of variables to transform the original equation \eqref{eqt:gCLM-model} into an equivalent equation (also known as the dynamic rescaling equation) for the rescaled profile $\Om$. In particular, an exact self-similar profile $\bar{\Om}$ corresponds to a steady state of the dynamic rescaling equation, and vice versa. Therefore, to establish asymptotically self-similar blowup for the original equation, one needs to perform stability analysis of the rescaled variable $\Om$ around some steady state $\bar{\Om}$.

For the gCLM model \eqref{eqt:gCLM-model}, self-similar finite-time blowup solutions have been established for a wide range of the parameter $a$ in the literature \cite{castro2009self,elgindi2020effects,lushnikov2021collapse,chen2020singularity,chen2021finite,huang2023self,huang2024self}. In particular, Huang, Qin, Wang, and Wei \cite{huang2024self} established a general existence result for exact self-similar blowup profiles covering the entire parameter range $a\leq1$. Their unified fixed-point construction yields highly regular profiles in $L^\infty(\R)\cap\dot H^1(\R)$. However, their result does not yield asymptotically self-similar blowups of \eqref{eqt:gCLM-model} for all $a\leq 1$. In fact, existing stability results for gCLM typically rely on explicit steady states (\cite{elgindi2021stable,chen2021finite} for $a$ close to $0$ and \cite{chen2020singularity} for $a$ close to $1/2$) or computer-assisted estimates (\cite{chen2021finite} for $a=1$). For general $a\leq1$, the profiles constructed in \cite{huang2024self} do not have explicit formulas, making a direct stability analysis for all $a\leq 1$ difficult.

The main goal of this paper is to fill the gap between the existence of exact self-similar blowup profiles and asymptotically self-similar blowups of \eqref{eqt:gCLM-model} for all $a\leq 1$. More precisely, we prove the following theorem.
\begin{theorem}\label{thm:main-blowup}
For every $a\leq1$, there exists some nonzero odd initial datum $\om_0\in C_c^\infty(\R)$ such that the corresponding solution of \eqref{eqt:gCLM-model} develops an asymptotically self-similar blowup in the sense of \eqref{eqt:gCLM-asymptotic-selfsimilar} at some finite time $T>0$.
\end{theorem}
Our proof of Theorem \ref{thm:main-blowup} is based on a finite-codimensional stability analysis of the dynamic rescaling equation around the exact self-similar blowup profile constructed in \cite{huang2024self}. More specifically, we prove that for any odd initial data sufficiently close to the steady state in a suitable weighted norm, one can adjust them in finitely many fixed directions so that the corresponding rescaled solution converges to the exact self-similar profile. The corrected initial data form a stable set, which plays the role of a finite-codimensional stable manifold.

The key to our finite-codimensional stability analysis is the out-going property of the flow field in the dynamic rescaling formulation. For the gCLM model, this means that the rescaled flow field associated with the profile $\bar\Omega$ points towards spatial infinity (if $\bar\Omega$ has support equal to $\R$) or towards the boundary of its support (if $\bar\Omega$ is compactly supported). 
This crucial out-going property allows us to construct suitable weighted Sobolev spaces in which we obtain damping estimates for the local part of the linearized operator of the dynamic rescaling equation around $\bar\Omega$. Meanwhile, we show that the remaining nonlocal terms of the linearized operator exhibit some compactness property. These two facts imply finite-codimensional stability of the linearized operator by the semigroup theory. With the linear stability established, we then use a shooting argument to
obtain nonlinear finite-codimensional stability of the dynamic rescaling equation around $\bar{\Omega}$, thereby proving Theorem~\ref{thm:main-blowup}.

We remark that the out-going property is not only a useful tool in establishing finite-codimensional stability of other Euler-type equations \cite{chen2026eulerII,shao2026selfsimilar}, but also an important ingredient in a full stability analysis \cite{elgindi2022stability,elgindi2021stable,chen2020singularity,chen2021boundary,chen2021finite,chen2021slightly,chen2021regularity,chen2022asymptotically,chen2022stable,hou2024blowup}. Here, full stability means stability under all sufficiently small perturbations within an open set in the prescribed class, rather than only for perturbations lying on a finite-codimensional stable (local) manifold. To establish full stability, these works combine the damping effect arising from the out-going property with refined estimates of the nonlocal terms, exploiting their special structure and cancellations, sometimes with computer assistance. In contrast, a finite-codimensional stability analysis essentially requires only the local damping provided by the out-going property, leaving most non-local estimates to a soft compactness argument. This makes the approach particularly useful when sharp quantitative information about the profile is unavailable.

\subsection{Asymptotically self-similar blowup of the 1D Hou--Luo model}
Following the same strategy, we also consider asymptotically self-similar blowup of the 1D Hou--Luo (HL) model:
\begin{equation}\label{eqt:1Dhouluo}
    \begin{aligned}
    &\om_t+u\om_x=\theta_x,\\
    &\theta_t+u\theta_x=0,\\
    & u_x=\mtx{H}(\om),
    \end{aligned}
\end{equation}
for $x\in \R$, where $\mtx{H}$ denotes the Hilbert transform on the real line. This model was first proposed by Luo and Hou \cite{luo2014potentially,luo2014toward} to acquire understanding of the numerically observed self-similar singularity formation of the 3D axisymmetric Euler equations on the solid boundary of an infinitely long cylinder. It also models the boundary induced singularity formation of the 2D Boussinesq equations \cite{luo2014toward,chen2022asymptotically} in the half-space $(x_1,x_2)\in\R\times\R_+$, 
\begin{equation}\label{eqt:Boussinesq}
\begin{aligned}
&\om_t + u_1\om_{x_1} + u_2\om_{x_2} = \theta_{x_1},\\
&\theta_t + u_1\theta_{x_1} + u_2\theta_{x_2} = 0,\\
&(u_1,u_2) = \nabla^{\perp}(-\Delta)^{-1}\om,
\end{aligned}
\end{equation}
for the 2D Boussinesq equations behave similarly to the 3D axisymmetric Euler equations away from the symmetry axis; see e.g. \cite{majda2002vorticity}. In fact, the boundary finite-time blowup of 3D axisymmetric Euler equations can be approximated by the boundary finite-time blowup of the 2D Boussinesq equations up to an asymptotically small perturbation \cite{luo2014toward,chen2022stable}.

As a simplified model of the boundary behavior of the Hou--Luo scenario, the 1D HL model has attracted significant attention regarding its singularity formation. 
Shortly after the original work of Luo and Hou \cite{luo2014toward}, Choi et al. \cite{choi2017finite} used a functional argument to prove the finite-time blowup of the HL model \eqref{eqt:1Dhouluo} and another 1D model known as the CKY model \cite{choi2015finite}. However, their approach was not able to capture the self-similar nature of the blowup. Years later, Chen, Hou, and Huang \cite{chen2022asymptotically} developed a novel analysis framework based on rigorous computer-assisted proofs to establish asymptotically self-similar blowup of \eqref{eqt:1Dhouluo} from smooth initial data. Instead of working with $(\omega,\theta)$, they established finite-time blowup for the variables $(\omega,\theta_x)$. Specifically, they first constructed an approximate self-similar profile $(\tilde{\Omega},\tilde{\Theta}_X)$ using numerical computation, and then performed a stability analysis of the dynamic rescaling equation of \eqref{eqt:1Dhouluo} around $(\tilde{\Omega},\tilde{\Theta}_X)$ and obtained the existence of an exact steady state solution $(\bar{\Omega},\bar{\Theta}_X)$ near $(\tilde{\Omega},\tilde{\Theta}_X)$ using a limit argument. This steady state corresponds to an exact self-similar blowup of \eqref{eqt:1Dhouluo} in the form of 
\begin{equation}\label{eqt:HL_exact_blowup}
    \om(x,t)=(T-t)^{\lambda}\bar{\Omega}\left(\frac{x}{(T-t)^{\gamma}}\right), \quad \theta_x(x,t)=(T-t)^{\mu}\bar{\Theta}_X\left(\frac{x}{(T-t)^{\gamma}}\right),
\end{equation}
with $\lambda=-1$, $\mu=-2$, $\gamma>0$. Moreover, this steady state was shown to be asymptotically stable. That is, for any initial value which is sufficiently close to $(\bar{\Omega},\bar{\Theta}_X)$ in a suitable norm, the corresponding solution of \eqref{eqt:1Dhouluo} would develop an asymptotically self-similar blowup that takes the form 
\begin{equation}\label{eqt:asymptotic_selfsimilar_HL}
    \om(x,t)=(T-t)^{\lambda}\Omega\left(\frac{x}{(T-t)^{\gamma}},t\right), \quad \theta_x(x,t)=(T-t)^{\mu}\Theta_X\left(\frac{x}{(T-t)^{\gamma}},t\right),
\end{equation}
where $\lambda=-1$, $\mu=-2$, $T$ is the blowup time, and the time-dependent profiles $\Omega(\cdot,t)$, $\Theta_X(\cdot,t)$ will converge to rescaled versions of $\bar{\Omega}$, $\bar{\Theta}_X$ as $t\to T$. This computer-assisted approach was finally generalized to the 2D Boussinesq and 3D axisymmetric Euler equations by Chen and Hou \cite{chen2022stable,chen2025singularity}, yielding a rigorous proof of the boundary self-similar singularity formation in the Hou--Luo scenario.

Although the asymptotically self-similar blowup of the HL model has been rigorously established in \cite{chen2022asymptotically}, the proof relies heavily on detailed computer-assisted estimates for a numerically constructed approximate profile, which can be quite uneasy for most readers to digest and reproduce. A natural question is therefore whether such a self-similar blowup can be proved without computer assistance. This question was partially answered by Huang, Qin, Wang, and Wei \cite{huang2025exact}, who proved the existence of a self-similar profile of \eqref{eqt:1Dhouluo} via a purely analytic fixed-point approach, following their work on the gCLM model \cite{huang2024self}. The resulting profile in their construction is smooth and enjoys monotonicity, convexity, and far-field decay estimates. Moreover, their numerical results show that the self-similar profile computed by their fixed-point iteration coincides with that constructed in \cite{chen2022asymptotically}, providing strong numerical evidence that the two constructions yield the same profiles up to rescaling. Although the analytic fixed-point construction in \cite{huang2025exact} provides a more qualitative characterization of the self-similar profile, it does not yield an easier path towards a purely analytic proof of full stability, due to the lack of accurate quantitative estimates (as those obtained with computer assistance in \cite{chen2022asymptotically}). Nevertheless, as in the gCLM case, the associated rescaled flow field constructed in \cite{huang2025exact} has the crucial out-going property, which allows us to establish finite-codimensional stability of the analytically constructed profile by a purely analytic argument. 

Therefore, another main result of this paper is a purely analytic proof of the  asymptotically self-similar blowup for the HL model \eqref{eqt:1Dhouluo} based on the self-similar profile constructed in \cite{huang2025exact} and the out-going property established thereby. More precisely, we prove the following theorem.
\begin{theorem}\label{thm:main-hl-blowup}
There exists some initial datum $(\omega_0,\theta_0)$ with  nonzero and odd $\omega_0,\theta_{0,x}\in C_c^\infty(\R)$ such that the corresponding solution of \eqref{eqt:1Dhouluo} develops an asymptotically self-similar blowup in the sense of \eqref{eqt:asymptotic_selfsimilar_HL} at some finite time $T>0$.
\end{theorem}
Compared with the gCLM case, the dynamic rescaling formulation of \eqref{eqt:1Dhouluo} leads to a more complicated coupled system, with additional cross terms in the energy estimates of the corresponding linearized equation around the analytically constructed steady state. Nevertheless, the out-going property allows us to choose suitable weights to control the principal contributions of these cross terms and obtain damping for the linearized operator after removing a compact part. We thus use arguments similar to those used in the gCLM case to establish finite-codimensional stability of the dynamic rescaling equation of \eqref{eqt:1Dhouluo} and obtain a purely analytic proof of Theorem~\ref{thm:main-hl-blowup}.

We remark that the out-going property also plays a crucial role in the stability analysis of the dynamic rescaling formulation of \eqref{eqt:1Dhouluo} in \cite{chen2022asymptotically}. In their computer-assisted framework, this property is combined with detailed quantitative estimates of the profile to design weights that exploit delicate cancellations in the energy estimates, yielding full asymptotic stability of the self-similar profile. By contrast, our argument relies mainly on the out-going property and leaves some compact terms outside the damping estimate, yielding only finite-codimensional stability. This comparison illustrates how rigorous computer-assisted estimates can yield stronger stability conclusions than the purely analytic estimates used here.

\subsection{Organization of the paper}
The rest of this paper is organized as follows. Sections~\ref{sec:rescaling}--\ref{sec:compact} concern the gCLM model. In Section~\ref{sec:rescaling}, we introduce its dynamic rescaling formulation and recall the existence and properties of the exact self-similar profiles constructed in \cite{huang2024self}. In Sections~\ref{sec:noncompact} and~\ref{sec:compact}, we establish finite-codimensional stability in the non-expanding and expanding cases, respectively, thereby proving Theorem~\ref{thm:main-blowup}. For the HL model, we present the dynamic rescaling formulation and profile properties in Section~\ref{sec:hl-rescaling} and prove Theorem~\ref{thm:main-hl-blowup} in Section~\ref{sec:hl-stability} using a finite-codimensional stability analysis. Appendix~\ref{app:hilbert} collects the properties of the Hilbert transform and the compactness estimate used in the analysis, and Appendix~\ref{app:profile-estimates} contains the proofs of the profile properties stated in Sections~\ref{sec:rescaling} and~\ref{sec:hl-rescaling}.

\section{Dynamic rescaling formulation and exact self-similar blowups of the gCLM model}
\label{sec:rescaling}

In this section, we introduce the dynamic rescaling formulation of the gCLM equation \eqref{eqt:gCLM-model}, which transforms the challenge of investigating finite-time blowup into the study of the convergence and stability of some steady state. We then review the existence and properties of the exact self-similar blowup profiles constructed in \cite{huang2024self}, which will be useful in the subsequent stability analysis.

Before we proceed, let us first clarify how to recover the velocity $u$ from the vorticity $\om$ in \eqref{eqt:gCLM-model}. Note that $u_x$ is determined by the Hilbert transform of $\om$:
\[ u_x(x,t)=\mtx{H}(\om(\cdot,t))(x)=\frac{1}{\pi}\pv\int_\R\frac{\om(y,t)}{x-y}\idiff y,\]
where $\pv$ denotes the Cauchy principal value. We need to recover $u$ from the relation $u_x=\mtx{H}(\om)$. Suppose that $u$ has a zero on $\R$. By \cite[Lemma~A.1]{huang2026singular}, the solutions of \eqref{eqt:gCLM-model} with the same initial data but different choices of zeros of $u$ differ only by a time-dependent spatial translation. Therefore, we may impose $u(0,t)=0$ without loss of generality and recover $u$ by
 \begin{equation}\label{eqt:gCLM-compute-u}
 u(x,t)=\int_0^x\mtx{H}(\om(\cdot,t))(y)\idiff y
 =\frac{1}{\pi}\int_\R\om(y,t)\ln\left|\frac{x-y}{y}\right|\idiff y.
\end{equation}
In the rest of the paper, we will always assume that $\om$ is odd symmetric. In this case, we can further simplify \eqref{eqt:gCLM-compute-u} to
\[ u(x,t)=\frac{1}{\pi}\int_0^\infty\om(y,t)\ln\left|\frac{x-y}{x+y}\right|\idiff y.\]

\subsection{Dynamic rescaling formulation}
Next, we introduce the dynamic rescaling formulation of the gCLM equation \eqref{eqt:gCLM-model}. Consider the following change of variables: 
\begin{equation}
    \omega(x,t)=C_{\om}(\tau)^{-1}\Omega\left(C_l(\tau)x,\tau(t)\right),\quad u(x,t)=C_{\om}(\tau)^{-1}C_{l}(\tau)^{-1}U\left(C_l(\tau)x,\tau(t)\right)
\end{equation}
with
\begin{equation}\label{eqt:gCLM-rescaling-factors}
\begin{aligned}
    C_{\om}(\tau)&=\exp\left(\int_0^\tau c_{\om}(s)\idiff s\right) ,\quad  C_{l}(\tau)&=\exp\left(\int_0^\tau c_l(s)\idiff s\right), \quad  t(\tau)=\int_0^{\tau}C_{\om}(s)\idiff s.
\end{aligned}
\end{equation}
Then \eqref{eqt:gCLM-model} is equivalent to the following dynamic rescaling equations: 
\begin{equation}\label{eqt:gCLM-dynamic-rescaling}
    \begin{aligned}
    &\Omega_\tau + (aU+c_lX)\Omega_X=(U_X+c_{\om})\Omega,\\
    &U_X =\mtx{H}(\Omega),\,\ U(0,\tau)=0,
    \end{aligned}
\end{equation}
where $X=C_l x$. The steady state solutions to \eqref{eqt:gCLM-dynamic-rescaling} satisfy the self-similar profile equation
\begin{equation}\label{eqt:gCLM-steady-state}
    \begin{aligned}
    &(aU+c_lX){\Omega}_X=(U_X+c_{\om})\Omega,\\
    &U_X =\mtx{H}({\Omega}),\,\ U(0)=0.
    \end{aligned}
\end{equation}
One can easily verify that if $(\bar{\Omega},\bar c_l,\bar c_\omega)$ is a solution to \eqref{eqt:gCLM-steady-state}, then it corresponds to exact self-similar finite-time blowups of the gCLM model in the form of \eqref{eqt:gCLM-self-similar-solution} with $\gamma=-\bar{c}_l/\bar{c}_\omega$, $\lambda=-1$.
Note that if $(\bar{\Omega}(X),\bar{c}_l,\bar{c}_{\om})$ is a steady state solution of \eqref{eqt:gCLM-dynamic-rescaling}, then
\begin{equation}
(\alpha \bar{\Omega}(\beta X),\alpha \bar{c}_l,\alpha \bar{c}_{\om})
\end{equation} is also a steady state solution for any $\alpha\in \R$, $\beta>0$, which indicates that this equation admits two degrees of freedom. Hence, in order to uniquely identify the profile up to rescaling, we need to impose two normalization conditions to fix these two degrees of freedom.

Moreover, if under some normalization conditions the solution tuple $(\Omega,c_l,c_{\omega})$ converges to a steady state $(\bar{\Omega},\bar{c_l},\bar{c}_{\om})$ of \eqref{eqt:gCLM-dynamic-rescaling} with $\bar{c}_{\om}<0$ as $\tau\to\infty$, then the solution would blow up in an asymptotically self-similar manner \eqref{eqt:gCLM-asymptotic-selfsimilar}. We refer the readers to \cite[Proposition~2.1]{huang2026singular} for a rigorous proof of this statement. Therefore, the study of asymptotically self-similar blowups of \eqref{eqt:gCLM-model} is reduced to the stability analysis of the dynamic rescaling equation \eqref{eqt:gCLM-dynamic-rescaling} around some steady state.
\subsection{Existence of the steady state solutions}
The existence of steady state solutions, of particular importance for the stability analysis, is itself a problem to be addressed. Due to the nonlinear and nonlocal nature of \eqref{eqt:gCLM-steady-state}, establishing the existence of such solutions for general values of $a$ is highly nontrivial. We refer the readers to \cite[Section~2.2]{huang2026singular} for a survey of previous results on the existence of self-similar profiles. In particular, a general existence result for all $a\leq1$ was established by Huang, Qin, Wang, and Wei \cite{huang2024self} through a unified fixed-point construction, yielding interiorly smooth solutions to \eqref{eqt:gCLM-steady-state}. We briefly restate this existence result below; see \cite[Theorems~2.1 and~4.10]{huang2024self} for details.

\begin{theorem}\label{thm:existence_of_self_similar_profile_gclm}
For each $a\leq 1$, \eqref{eqt:gCLM-steady-state} admits a solution $(\bar{\Omega}(X),\bar{c}_l,\bar{c}_{\om})$ with $\bar{c}_{\om} = -1$. $\bar{\Omega}$ is an odd-symmetric function that satisfies $\bar{\Omega}\in L^\infty(\mathbb{R})\cap \dot{H}^1(\mathbb{R})$. The profile satisfies that $\bar{\Omega}_X(0)<0$ and that $-\bar{\Omega}/X$ is decreasing in $X$ and convex in $X^2$ on $[0,+\infty)$. Depending on the sign of $\bar{c}_l$, one of the following happens:
\begin{enumerate}
\item\label{item:profile-type-1} $\bar c_l> 0$: $\bar{\Omega}$ is strictly negative on $(0,+\infty)$ and smooth on $\mathbb{R}$, and 
\[\lim_{X\to+\infty}X^{1/\bar c_l}\bar\Omega(X)=C_\Omega\quad\text{for some }C_\Omega<0;\]
\item\label{item:profile-type-2} $\bar{c}_l=0$: $\bar{\Omega}$ is strictly negative on $(0,+\infty)$ and smooth on $\mathbb{R}$, and $|\bar{\Omega}(X)|\lesssim \econst^{-\delta X^2}$ for some $\delta>0$;
 \item\label{item:profile-type-3} $\bar{c}_l<0$: $\bar{\Omega}$ is compactly supported on $[-L,L]$ for some $L>0$, strictly negative on $(0,L)$, and smooth in the interior of $(-L,L)$. Moreover, there exist $C_L>0$ and $p\geq1$ such that
\begin{equation}\label{eqt:gCLM-profile-edge-asymptotic}
 \lim_{X\to L^{-}}\frac{-\bar\Om(X)}{(L-X)^p}=C_L, \quad p>1\ \text{if }\, a<1,\quad p=1\ \text{if }\,a=1.
\end{equation}
\end{enumerate}
Moreover, there exist constants $1/2<\underline{a}\leq \overline{a}<1$ such that a solution of \hyperref[item:profile-type-1]{type~(\ref*{item:profile-type-1})} exists for every $a<\underline{a}$, and a solution of \hyperref[item:profile-type-3]{type~(\ref*{item:profile-type-3})} exists for every $\overline{a}<a\leq1$.
\end{theorem}

\begin{remark}
In the above theorem, $\bar c_\om$ is normalized to $-1$, and thus $\gamma=-\bar{c}_l/ \bar c_{\om}=\bar{c}_l$ in \eqref{eqt:gCLM-self-similar-solution}. Therefore, the sign of $\bar c_l$ determines the spatial feature of the blowup: the blowup is focusing if $\bar c_l>0$, expanding if $\bar c_l<0$, and neither expanding nor focusing if $\bar c_l=0$. In what follows, we will consider the non-expanding case $\bar c_l\geq0$ and the expanding case $\bar c_l<0$ separately in Section~\ref{sec:noncompact} and Section~\ref{sec:compact}.
\end{remark}

The following theorem is crucial for the linear stability analysis in Section~\ref{sec:noncompact} and Section~\ref{sec:compact}. Although this result is not stated as a separate theorem in \cite{huang2024self}, it follows directly from the fixed-point representation in \cite[Section~2.2 and Proposition~3.1]{huang2024self}, the positivity and monotonicity properties in \cite[Corollary~3.5]{huang2024self}, and the steady state equation \eqref{eqt:gCLM-steady-state}.

\begin{theorem}\label{thm:out-going}
Let \((\bar\Om,\bar{U},\bar c_l,\bar c_{\om})\) be the self-similar profile given by Theorem~\ref{thm:existence_of_self_similar_profile_gclm} with \(\bar c_{\om}=-1\). Depending on the sign of $\bar c_l$, the following properties hold:
\begin{enumerate}
\item \(\bar c_l\geq0\): The rescaled flow field \(a\bar U(X)+\bar c_lX\) is out-going, i.e.
\[a\bar U(X)+\bar c_lX>0,\quad X>0.\]
At the origin, it holds that
\[ a\bar U_X(0)+\bar c_l=\bar U_X(0)-1>0. \]
\item \(\bar c_l<0\): Let $L$ be the support boundary of \(\bar\Om\) given in Theorem~\ref{thm:existence_of_self_similar_profile_gclm}. Then the rescaled flow field \(a\bar U(X)+\bar c_lX\) is out-going in the sense that
\[ a\bar U(X)+\bar c_lX
 \begin{cases}
  >0,&0<X<L,\\
  =0,&X=L,\\
  <0,&X>L.
 \end{cases}
\]
At the origin,
\[ a\bar U_X(0)+\bar c_l=\bar U_X(0)-1>0.\]
Moreover, at the support boundary $X=L$,
\begin{equation}\label{eqt:gCLM-npc-BL}
 \bar U_X(L)-1\leq a\bar U_X(L)+\bar c_l<0.
\end{equation}
\end{enumerate}
\end{theorem}

\subsection{More information about the steady state solutions}
Theorem \ref{thm:existence_of_self_similar_profile_gclm} and Theorem \ref{thm:out-going} stated above are already established in \cite{huang2024self}. 
To perform the finite-codimensional stability analysis of the dynamic rescaling equation \eqref{eqt:gCLM-dynamic-rescaling}, we need additional properties of the self-similar profiles beyond those given explicitly in \cite{huang2024self}. We summarize these properties in the following lemma, whose proof follows from the steady-state equation \eqref{eqt:gCLM-steady-state} and Theorem~\ref{thm:existence_of_self_similar_profile_gclm}. For the reader's convenience, we provide a detailed proof in Appendix~\ref{app:profile-estimates}.

\begin{lemma}\label{lem:profile-derivative-estimates}
Let \((\bar\Om,\bar{U},\bar c_l,\bar c_{\om})\) be the self-similar profile given by Theorem~\ref{thm:existence_of_self_similar_profile_gclm} with \(\bar c_{\om}=-1\). Then as $X\to +\infty$, it holds that 
\[ \bar{U}_X(X) \to 0, \quad \frac{\bar{U}}{X} \to 0. \]
Moreover, depending on the sign of $\bar{c}_l$, the following estimates hold:
\begin{enumerate}
\item \(\bar c_l>0\):
\begin{equation}\label{eqt:gCLM-profile-two-derivatives}
 \left|\bar\Om(X)\right|+X\left|\bar\Om_X(X)\right|+X^2\left|\bar\Om_{XX}(X)\right| \leq CX^{-1/\bar c_l},\quad X\geq1,
\end{equation}
and 
\[\left\|X\bar U_{XX}(X)\right\|_{L^\infty(\R)}<\infty.\]
\item \(\bar c_l=0\): Let \(\delta >0\) be given in Theorem~\ref{thm:existence_of_self_similar_profile_gclm} (2). There exists a polynomial \(P(X)\) such that
\[ \left|\bar\Om(X)\right|+X\left|\bar\Om_X(X)\right|+X^2\left|\bar\Om_{XX}(X)\right|\leq P(X)\econst^{-\delta X^2},\quad X\geq1,\]
and 
\begin{equation}
 \left\|X\bar U_{XX}(X)\right\|_{L^\infty(\R)}<\infty, \quad a\bar U(X)\asymp X^{-1}\quad \text {as }X\to+\infty.
\end{equation}
Here and below, $f\asymp g$ means that there exist constants $0<c\leq C<\infty$ such that $cg\leq f\leq Cg$.
\item \(\bar c_l<0\): with \(L\) as in
Theorem~\ref{thm:existence_of_self_similar_profile_gclm},
\begin{equation}\label{eqt:gCLM-npc-edge-H2}
 \left\|\bar\Om_X\right\|_{L^\infty(\R)}<\infty,\quad \int_{L/2}^L(L-X)\left|\bar\Om_{XX}(X)\right|^2\idiff X<\infty,
\end{equation}
where \(\bar\Om_{XX}\) denotes the classical interior derivative. Moreover, 
\begin{equation}
 X\bar U_{XX}(X)=O\left(X^{-2}\right)\quad \text {as }X\to+\infty.
\end{equation}
\end{enumerate}
\end{lemma}

\section{Finite-codimensional stability analysis of gCLM: non-expanding case}
\label{sec:noncompact}
In this section, we establish finite-codimensional stability of the dynamic rescaling equation \eqref{eqt:gCLM-dynamic-rescaling} for each $a$ such that \eqref{eqt:gCLM-dynamic-rescaling} admits a steady-state solution \((\bar\Om,\bar U,\bar c_l,-1)\) of \hyperref[item:profile-type-1]{type~(\ref*{item:profile-type-1})} or \hyperref[item:profile-type-2]{type~(\ref*{item:profile-type-2})} given by Theorem \ref{thm:existence_of_self_similar_profile_gclm}. In both cases, we have $\bar c_l\geq 0$ and the self-similar profile $\bar{\Om}$ is smooth on $\R$.
To simplify the notation, we denote the space-time variable $(X,\tau)$ by $(x,t)$ throughout the discussion of the dynamic rescaling equations. We thus rewrite \eqref{eqt:gCLM-dynamic-rescaling} as
\begin{equation}\label{eqt:gCLM-rescaled}
 \Om_t+(aU+c_lx)\Om_x=(U_x+c_\om)\Om,
 \quad U_x=\mtx{H}(\Om),
 \quad U(0,t)=0.
\end{equation}

Denote the rescaled flow field by 
\[\bar{V}(x):=a\bar U(x)+\bar c_lx.\]
Theorem~\ref{thm:out-going} and Lemma~\ref{lem:profile-derivative-estimates} yield
\begin{equation}\label{eqt:gCLM-A-comparable}
 \bar{V}(x)>0,\quad
 \bar{V}(x)\asymp
 \begin{cases}
  x, & \bar c_l>0,\\
  x\la x\ra^{-2}, & \bar c_l=0,
 \end{cases}
 \quad x>0.
\end{equation}
Intuitively, this out-going background flow transports perturbations around the steady state towards spatial infinity, allowing us to construct suitable singular weighted Sobolev spaces and establish finite-codimensional stability.

Let $f$ be the perturbation of $\Omega$ around the steady state, i.e. $f:=\Om - \bar\Om$. We impose the following normalization conditions in \eqref{eqt:gCLM-rescaled}:
\begin{equation}\label{eqt:gCLM-normalization}
 c_l=\bar c_l-a\mtx{H}(f)(0),\quad
 c_\om=-1-\mtx{H}(f)(0).
\end{equation}
With \eqref{eqt:gCLM-normalization}, direct calculation yields
\begin{equation}\label{eqt:gCLM-fixed-slope}
 \Om_{x,t}(0,t) =\big((1-a)U_x(0,t)+c_\om-c_l\big)\Om_x(0,t) =\big((1-a)\bar{U}_x(0)-1-\bar{c}_l\big)\Om_x(0,t)=0.
\end{equation}
Define
\begin{equation}\label{eqt:gCLM-ABQV}
\mtx{H}_0(f):=\mtx{H}(f)-\mtx{H}(f)(0),\quad
 \mtx{V}(f)(x):=\int_0^x\mtx{H}_0(f)(y)\idiff y.
\end{equation}
We rewrite the equation for the perturbation $f$ as 
\begin{equation}\label{eqt:gCLM-perturbation}
 f_t=\mtx{L}(f)+\mtx{N}(f),
\end{equation}
where
\begin{align}
 \mtx{L}(f)&:=-\bar{V}f_x+(\bar{U}_x-1)f+\bar\Om \mtx{H}_0(f)-a\bar\Om_x\mtx{V}(f),\label{eqt:gCLM-perturbation_linear}\\
 \mtx{N}(f)&:=-a\mtx{V}(f)f_x+\mtx{H}_0(f)f.\label{eqt:gCLM-perturbation_nonlinear}
\end{align}
We will study the perturbation equation \eqref{eqt:gCLM-perturbation} in the rest of this section.
\subsection{Choice of singular weights}
Suppose that the initial perturbation satisfies $f_x(0,0)=0$. By \eqref{eqt:gCLM-fixed-slope} we know that 
\begin{equation}
      f_x(0,t)=0\label{eqt:gCLM-slope0}
\end{equation}for all $t>0$.
To capture the out-going effect established in Theorem \ref{thm:out-going}, we introduce singular weighted $L^2$ and $H^1$ spaces which are adapted to \eqref{eqt:gCLM-slope0} as follows. Choose 
\begin{equation}\label{eqt:gCLM-alpha-range-intro}
 \alpha\in \left(-1,\min\left\{1,2/\bar{c}_l-1\right\}\right),
\end{equation}
where we set $2/\bar{c}_l=+\infty$ when $\bar{c}_l=0$. Define 
\[\phi_0(x):= \frac{\la x\ra ^{\alpha+4}}{x^4},\]
where $\la x\ra:= (x^2+1)^{1/2}$. We choose the singular $L^2$ weight as $\phi_0$ times a corrector:
 \begin{equation}
 \phi(x):=\econst^{\kappa(x)}\phi_0(x),
\end{equation}
and define the singular $H^1$ weight as
\[\psi(x):=x^2\phi(x).\]
Here $\kappa(x)$ is a bounded smooth function to be determined later. One can check that $\phi$ and $\psi$ satisfy the asymptotics 
\begin{align}
 \phi(x)&\asymp |x|^{-4},\quad \text{as }x\to0 , \quad \phi(x)\asymp x^\alpha,\quad \,\,\,\,\text{as } x\to+\infty,\\
\psi(x)&\asymp |x|^{-2},\quad \text{as }x\to0 , \quad \psi(x)\asymp x^{\alpha+2},\quad \text{as }x\to+\infty.
\end{align}
We will choose $\kappa$ to decrease sufficiently fast in the $O(1)$ region to capture the out-going property of $\bar{V}$. Together with the above asymptotics of $\phi$, this will give the damping estimates for the local part of the linearized operator $\mtx{L}$. The choice of $\kappa$ is determined in the next lemma.  
\begin{lemma}\label{lem:weight-parameters}
There are a bounded smooth even function $\kappa$, constant near infinity, and a constant $\eta_0>0$ such that, with
\begin{equation}\label{eqt:gCLM-weight}
      \phi=\econst^{\kappa}\phi_0, \quad \psi=x^2\phi,
\end{equation}
the following inequalities hold 
\begin{align}
 d_0(x)&:=(\bar{U}_x-1)+\frac{\bar{V}_x}{2}+\frac{\bar{V}}{2}\frac{\phi_x}{\phi}\leq-\eta_0,\label{eqt:gCLM-d0}\\
 d_1(x)&:=(\bar{U}_x-1)-\frac{\bar{V}_x}{2}+\frac{\bar{V}}{2}\frac{\psi_x}{\psi}\leq-\eta_0.\label{eqt:gCLM-d1}
\end{align}
\end{lemma}
\begin{proof}
Substituting \eqref{eqt:gCLM-weight} into the formulas for $d_0$ and $d_1$ yields
\begin{align*}
 d_0(x)&=(\bar{U}_x-1)+\frac{\bar{V}_x}{2}+\frac {\bar{V}}{2x}\frac{x(\phi_0)_x}{\phi_0}+\frac{\bar{V}}{2}\kappa_x=:d_0^{(0)}(x)+\frac{\bar{V}}{2}\kappa_x,\\
 d_1(x)&=(\bar{U}_x-1)-\frac{\bar{V}_x}{2}+ \frac{\bar{V}}{x}+ \frac {\bar{V}}{2x}\frac{x(\phi_0)_x}{\phi_0}+\frac{\bar{V}}{2}\kappa_x=:d_1^{(0)}(x)+\frac{\bar{V}}{2}\kappa_x.
\end{align*}
By Theorem~\ref{thm:out-going} and Lemma~\ref{lem:profile-derivative-estimates}, we have
\[d_0^{(0)}(0)=d_1^{(0)}(0)=-\frac{\bar{V}_x(0)}2<0.\]
\[d_0^{(0)}(+\infty)=d_1^{(0)}(+\infty)= -1+\frac{(1+\alpha)\bar c_l}{2}<0.\]
Therefore, there exist $\delta>0$ and $\eta_0$ such that 
\[\max\left\{d_0^{(0)},d_1^{(0)}\right\}\leq-\eta_0\] 
holds on $\mathcal U:= (0,\delta)\cup(1/\delta,+\infty)$.
We then take $\mathcal U'=(0,\delta/2)\cup(2/\delta,+\infty)$ and choose a smooth cut-off function 
$\zeta\in C_c^\infty(0,+\infty)$ such that
\[0\leq\zeta\leq1,\quad \zeta=0\ \hbox{on }\mathcal U',\quad \zeta=1\ \hbox{on }(0,+\infty)\setminus\mathcal U.\]
The compact set $\mathcal{K}=(0,+\infty)\setminus\mathcal U$ contains no zero of $\bar{V}$. Hence, we have 
\[ \bar{V}_*:=\min_{x\in\mathcal{K}}\bar{V}>0,\quad
 M_*:=\max_{x\in \mathcal K,\,j=0,1}\left\{d_j^{(0)}\right\}<\infty. \]
Choose a sufficiently large $M>0$ such that $M_*-M\bar{V}_*/2\leq-\eta_0$. We define
\begin{equation}
 \kappa(x)=-M\int_0^x \zeta(s)\idiff s, \quad \text {for } x\geq 0
\end{equation}
and extend the definition to $(-\infty,0)$ evenly. Then one can see that $\kappa$ is a bounded smooth function such that \eqref{eqt:gCLM-d0} and \eqref{eqt:gCLM-d1} hold. The proof is thus completed.
\end{proof}

In the rest of this section we choose the weight functions $\phi$ and $\psi$ according to Lemma \ref{lem:weight-parameters}. With these weight functions, we will establish nonlinear stability in a weighted $H^1$-norm of the perturbation. Let us define the weighted $L^2$ space
\[ \cX_0:=\left\{f\in L^2_{\rm loc}(\R):f\text{ is odd},\ \|f\|_{\cX_0}<\infty\right\},\]
with
\[ \left\|f\right\|_{\cX_0}^2:=\int_\R f(x)^2\phi(x)\idiff x, \]
and the weighted $H^1$ space
\[ \cX_1:=\left\{f\in\cX_0\cap H^1_{\rm loc}(\R):\|f\|_{\cX_1}<\infty\right\},\]
with
\[ \left\|f\right\|_{\dot{\cX}_1}^2:=\int_\R f_x(x)^2\psi(x)\idiff x,\quad
 \left\|f\right\|_{\cX_1}^2:=\left\|f\right\|_{\dot{\cX}_1}^2+\left\|f\right\|_{\cX_0}^2.\]
 Before proving the main results of this section, let us first establish some estimates which will be used in the later stability analysis.

\begin{lemma}\label{lem:coefficient-bounds}
There is \(C=C(\alpha,\phi)\) such that
\begin{enumerate}
  \item     
\begin{equation} \label{eqt:gCLM-low-coeff}
 |\mtx{H}(f)(0)|+\left\|\frac{\mtx{V}(f)}{x}\right\|_{L^\infty(\R)}+\left\|\frac{\mtx{V}(f)}{x^2}\right\|_{L^\infty(\R)} \leq C\left\|f\right\|_{\cX_0}.
\end{equation}
\item For $p=1,2$,
\begin{equation}\label{eqt:gCLM-low-coeff1}
\left\|\frac{\mtx{V}(f)}{x^{p+1}}\right\|_{L^2(\R)}\leq C\left\|f\right\|_{\cX_0}.
\end{equation}
\item
 \begin{equation}\label{eqt:gCLM-high-coeff}
      \left\|f\right\|_{L^\infty(\R)}+\left\|\mtx{H}(f)\right\|_{L^\infty(\R)} \leq C\left\|f\right\|_{\cX_1}.
 \end{equation}
\end{enumerate}
\end{lemma}

\begin{proof}
We assume that $f$ is a smooth function with compact support. The general case can be obtained by approximation. We bound $|\mtx{H}(f)(0)|$ as follows:
\begin{equation}\label{eqt:gCLM-H0-proof}
 |\mtx{H}(f)(0)|=\frac{1}{\pi}\left|\int_{\R}\frac{f(x)}{x} \idiff x\right|\leq\frac1\pi\left(\int_\R\frac{1}{x^2\phi(x)}\idiff x\right)^{1/2}
 \left\|f\right\|_{\cX_0}\lesssim\left\|f\right\|_{\cX_0}.
\end{equation}
For \(|x|\leq1\), Lemma~\ref{lem:core-identities} gives
\[  \mtx{V}(f)(x)=\int_0^xy^2\mtx{H}\left(f/y^2\right)(y)\idiff y,\]
therefore 
 \[ \left|\frac{\mtx{V}(f)(x)}{x^2}\right|\leq \frac{\left|\int_0^x y^4 \idiff y\right|^{1/2}}{x^2}\left\|\mtx{H}\left(\frac{f}{y^2}\right)\right\|_{L^2(\R)} \lesssim |x|^{1/2}\left\|\frac{f}{y^2}\right\|_{L^2(\R)}\lesssim \left\|f\right\|_{\cX_0}. \]
For \(|x|\geq1\), write 
\[ \frac{\mtx{V}(f)(x)}x=-\mtx{H}(f)(0)+\frac1x\int_0^x\mtx{H}(f)(y)\idiff y.\]
Note that $\mtx{H}(f)(0)$ is controlled in \eqref{eqt:gCLM-H0-proof}. Hence we only need to estimate the second term. We can bound the integral as follows:
\begin{align*} 
      \left|\frac1x\int_0^x\mtx{H}(f)(y)\idiff y\right|&\leq \frac{1}{|x|}\left|\int_0^x \la y\ra^{-\alpha} \idiff y \right|^{1/2}\left\|\mtx{H}(f)\la y\ra^{\alpha/2}\right\|_{L^2(\R)}\\
&\lesssim |x|^{-(1+\alpha)/2}\left\|f\la y\ra^{\alpha/2}\right\|_{L^2(\R)}\lesssim \left\|f\right\|_{\cX_0}.
\end{align*}
We have used the $A_2$ estimate \eqref{eqt:common-A2-weight} in the second inequality and \eqref{eqt:gCLM-alpha-range-intro} in the last inequality. This proves \eqref{eqt:gCLM-low-coeff}.

Using Lemma~\ref{lem:core-identities} and Hardy's inequality, we have
\begin{align*}
      \left\|\frac{\mtx{V}(f)}{x^{p+1}}\right\|_{L^2(\R)}&= \left\|\frac{\int_0^x(\mtx{H}(f)(y)-\mtx{H}(f)(0))\idiff y}{x^{p+1}}\right\|_{L^2(\R)}\\ &\lesssim \left\|\frac{\mtx{H}(f)-\mtx{H}(f)(0)}{x^{p}}\right\|_{L^2(\R)}=\left\| \frac{f}{x^p} \right\|_{L^2(\R)} \leq C\left\|f\right\|_{\cX_0},
\end{align*}
which implies \eqref{eqt:gCLM-low-coeff1}.

Note that $f(0)=0$. Therefore for any $x\in \R$,
\begin{align*}
      f^2(x)&=\left|\int_0^x 2f_x(y)f(y)\idiff y \right|\lesssim \left\|f_x \la x\ra^{-\alpha/2}\right\|_{L^2(\R)} \left\|f \la x\ra^{\alpha/2}\right\|_{{L^2(\R)}}\\ 
&\lesssim  \left\|f_x\la x\ra^{\alpha/2+1}\right\|_{L^2(\R)} \left\|f\la x\ra^{\alpha/2}\right\|_{{L^2(\R)}}  \lesssim \left\|f\right\|_{\cX_1}^2.
\end{align*}
We can bound $\|\mtx{H}(f)\|_{L^{\infty}(\R)}$ using a similar argument. Since $f$ has compact support, we have $\mtx{H}f(\infty)=0$. Therefore, 
\begin{align*}
 \mtx{H}(f)(x)^2 & =2\left|\int_{-\infty}^x\mtx{H}(f)(y)\mtx{H}(f)_x(y)\idiff y\right| \lesssim \left\|\mtx{H}(f)\la x\ra^{\alpha/2}\right\|_{L^2(\R)} \left\|\mtx{H}(f_x)\la x\ra^{-\alpha/2}\right\|_{L^2(\R)}\\
 &\lesssim \left\|f\la x\ra^{\alpha/2}\right\|_{L^2(\R)}\left\|f_x\la x\ra^{\alpha/2+1}\right\|_{L^2(\R)}\lesssim \left\|f\right\|_{\cX_1}^2,
\end{align*}
where we have used \eqref{eqt:common-A2-weight} in the second inequality. This proves \eqref{eqt:gCLM-high-coeff}.
\end{proof}

\subsection{Linear stability estimates via semigroups}\label{sec:linear}
In this subsection we study the linearized equation 
\begin{equation}\label{eqt:gCLM-linearized-equation}
      f_t=\mtx{L}(f)
\end{equation}
via semigroups, where $\mtx{L}$ is the linear operator defined in \eqref{eqt:gCLM-perturbation_linear}. Our goal is to establish finite-codimensional stability of \eqref{eqt:gCLM-linearized-equation} by decomposing $\mtx{L}$ as the sum of a damping term $\mtx{L}_{\rm stb}$ which satisfies certain decay estimate $\|\econst^{t\mtx{L}_{\rm stb}} \|\lesssim \econst^{-\eta_0 t}$ and a compact term $\mtx{K}$ which contributes at most finitely many unstable modes. Write
\begin{equation}\label{eqt:gCLM-decomposition-linear1} 
      \mtx{L}(f)=-\bar{V}f_x+(\bar{U}_x-1)f+\mtx{P}(f)+\mtx{K}_{\rm pr}(f),\quad
 \mtx{P}(f):=\bar\Om \mtx{H}_0(f),\quad \mtx{K}_{\rm pr}(f):=-a\bar\Om_x\mtx{V}(f).
\end{equation}
We will show later that the local term $-\bar{V}f_x+(\bar{U}_x-1)f$ is a damping term that satisfies the desired decay estimates. On the other hand, the nonlocal term is difficult to estimate as we don't have the explicit expression of $\bar{\Omega}$. To establish finite-codimensional stability, one expects to prove the compactness of the nonlocal term. However, \(\mtx{P}\) is not compact in general since \(\mtx{P}(f)\) has the same order as the Hilbert transform \(\mtx{H}(f)\). To overcome this difficulty, we need to extract some compactness by using symmetry. Note that the Hilbert transform is skew-symmetric under the unweighted \(L^2\) inner product. This motivates us to decompose \(\mtx{P}\) into its skew-symmetric and symmetric parts under \(\cX_0\) and establish the compactness of the symmetric part. Let \(\mtx{P}^*\) denote the adjoint of \(\mtx{P}\) with respect to \(\cX_0\), and set
\[ \mtx{P}_{\rm sym}:=\frac{\mtx{P}+\mtx{P}^*}{2},\quad \mtx{P}_{\rm skew}:=\frac{\mtx{P}-\mtx{P}^*}{2}.\]
Following the above strategy, we define
\begin{equation}\label{eqt:gCLM-decomposition-linear}
 \mtx{L}_{\rm stb}(f):=-\bar{V}f_x+(\bar{U}_x-1)f+\mtx{P}_{\rm skew}(f),\quad
 \mtx{K}(f):=\mtx{K}_{\rm pr}(f)+\mtx{P}_{\rm sym}(f),
\end{equation}
and apply standard spectral theory to $\mtx{L}_{\rm stb}$ and $\mtx{L}=\mtx{L}_{\rm stb}+\mtx{K}$ to establish finite-codimensional stability of \eqref{eqt:gCLM-linearized-equation}.
\subsubsection{Compact and bounded perturbations}We first establish the compactness of the nonlocal term $\mtx{K}(f)$ in the following lemma.
\begin{lemma} \label{lem:gCLM_symmetric_compact}
 $\mtx{P}$ is bounded on $\cX_0$. Moreover, $\mtx{P}_{\rm sym}$ and $\mtx{K}_{\rm pr}$ are compact operators on $\cX_0$. 
\end{lemma}
\begin{proof} 
Let us first consider the case $\bar c_l>0$. For any $f\in \cX_0$, we compute that
\begin{equation}\label{eqt:gCLM-P-bounded}
      \begin{aligned}
      \|\mtx{P}(f)\|_{\cX_0}&\leq \left\|\bar{\Omega}\mtx{H}_0(f)\phi^{1/2}\right\|_{L^2[-1,1]}+ \left\|\bar{\Omega}\mtx{H}_0(f)\phi^{1/2}\right\|_{L^2(\R\backslash[-1,1])}\\
      & \lesssim \left\|\frac{\mtx{H}_0(f)}{x}\right\|_{L^2[-1,1]}+ \left\|\mtx{H}(f)|x|^{\alpha/2-1/\bar{c}_l}\right\|_{L^2(\R\backslash[-1,1])}+\left\|\mtx{H}(f)(0)|x|^{\alpha/2-1/\bar{c}_l}\right\|_{L^2(\R\backslash[-1,1])}\\
      & \lesssim \left\|\mtx{H}\left(\frac{f}{x}\right)\right\|_{L^2(\R)}+\left\|\mtx{H}(f) \la x\ra ^{\alpha/2}\right\|_{L^2(\R)}+ |\mtx{H}(f)(0)|\\
      & \lesssim  \left\|\frac{f}{x}\right\|_{L^2(\R)}+\left\|f \la x\ra ^{\alpha/2}\right\|_{L^2(\R)}+ |\mtx{H}(f)(0)|\\
      & \lesssim \|f\|_{\cX_0}.
      \end{aligned}
\end{equation}
We have used the asymptotics of $\bar{\Omega}$ in Lemma~\ref{lem:profile-derivative-estimates} and the asymptotics of the weight $\phi$ in the second inequality, \eqref{eqt:common-hilbert-identities} and the $L^2$ integrability of $|x|^{\alpha/2-1/\bar{c}_l}$ for $|x|\geq 1$ in the third inequality, the $A_2$ estimate \eqref{eqt:common-A2-weight} in the fourth inequality, and \eqref{eqt:gCLM-low-coeff} in the last inequality. This proves that $\mtx{P}$ is bounded on $\cX_0$.

Next, we establish the compactness of $\mtx{P}_{\rm sym}$. Let $L^2_{\rm odd}(\R)$ denote the subspace of $L^2(\R)$ which consists of all odd symmetric functions. We define an isometry
\begin{equation}
      \mtx{Q}: \cX_0\longrightarrow L^2_{\rm odd}(\R), \quad \mtx{Q}(f)=\phi^{1/2}f. 
\end{equation}
To prove that $\mtx{P}_{\rm sym}$ is a compact operator on $\cX_0$, we only need to show that $\mtx{Q} \mtx{P}_{\rm sym}\mtx{Q}^{-1}$ is compact on $L^2_{\rm odd}(\R)$. Note that
\[ \mtx{H}_0(f)(x)=\mtx{H}(f)(x)-\mtx{H}(f)(0)=\frac1\pi\pv\int_\R \frac{x}{y(x-y)}f(y)\idiff y. \]
Therefore the kernel of \(\mtx{Q} \mtx{P} \mtx{Q}^{-1}\) is 
\[K_P(x,y)= \frac{\bar\Om(x)x\phi(x)^{1/2}}{\pi y(x-y)\phi(y)^{1/2}}\]
and the kernel of $\mtx{Q} \mtx{P}_{\rm sym}\mtx{Q}^{-1}$ is given by
\[K_s(x,y) =\frac12\bigl(K_P(x,y)+{K_P(y,x)}\bigr).\]
We are going to show that $K_s(x,y)\in L^2(\R\times \R)$ and therefore $\mtx{P}_{\rm sym}$ is Hilbert--Schmidt. Set
\[ N(x):=x^3\phi(x)\bar\Om(x).\]
Then $N$ is smooth on $\R$. Moreover, by \eqref{eqt:gCLM-profile-two-derivatives} and the definition of $\phi$,
\begin{equation}\label{eqt:gCLM-N-bound}
      |N(x)|+|xN_x(x)|\lesssim \la x\ra^{3+\alpha-1/\bar c_l}.
\end{equation}
We write the symmetric kernel $K_s$ as
\begin{equation}  
      K_s(x,y) =\frac{N(x)y-N(y)x}{2\pi(x-y)x^2y^2\sqrt{\phi(x)\phi(y)}}.
\end{equation}
One can see that the singularity in the original kernel $K_{P}$ is cancelled by the symmetrization operation. Indeed, for any $x,y\in [-1,1]$,
\begin{equation}\label{eqt:gCLM-bound-of-N}
      \left|\frac{N(x)y-N(y)x}{x-y}\right|=\left|\frac{N(x)-N(y)}{x-y}y-N(y)\right|\leq 2\left\|N_x\right\|_{L^{\infty}[-1,1]}+\left\|N\right\|_{L^{\infty}[-1,1]}.
\end{equation}
Since \begin{equation}\label{eqt:gCLM-phi-asymp}
      x^4\phi\asymp \la x \ra^{4+\alpha},
\end{equation} we deduce $K_s$ is bounded on $[-1,1]^2$. We then estimate $K_s$ outside $[-1,1]^2$. As $K_s(x,y)$ is symmetric in $(x,y)$, it suffices to consider the following two cases separately.
\begin{itemize}
      \item $|y|<|x|/2$: By \eqref{eqt:gCLM-N-bound} and \eqref{eqt:gCLM-phi-asymp}, we have 
\begin{align*}
      |K_s(x,y)|& \lesssim \frac{|N(x)y|+|N(y)x|}{|x|x^2\sqrt{\phi(x)}y^2\sqrt{\phi(y)}}\lesssim \frac{|x|^{3+\alpha-1/\bar c_l}\la y \ra+\la y\ra^{3+\alpha-1/\bar c_l}|x|}{|x|^{3+\alpha/2}\la y\ra^{2+\alpha/2}}\\
      & = |x|^{\alpha/2-1/\bar c_l}\la y\ra^{-1-\alpha/2} +|x|^{-2-\alpha/2}\la y\ra^{1+\alpha/2-1/\bar c_l}.
\end{align*}
      \item $|x|/2<|y|<|x|$: If $xy>0$, we use the mean-value theorem to obtain
  \[|K_s(x,y)|\lesssim \frac{|N(x)y-N(y)x|/|x-y|}{x^2\sqrt{\phi(x)}y^2\sqrt{\phi(y)}}\lesssim \frac{|N(x)|+|x|\|N_x\|_{L^{\infty}(\min(x,y),\max\{x,y\})}}{|x|^{2+\alpha/2}\la y\ra^{2+\alpha/2}}\lesssim |x|^{-1-1/\bar c_l}. \]
  If $xy<0$, then 
  \[ |K_s(x,y)|\lesssim \frac{(|N(x)y|+|N(y)x|)/(|x|+|y|)}{x^2\sqrt{\phi(x)}y^2\sqrt{\phi(y)}}\lesssim \frac{|N(x)|+|N(y)|}{|x|^{2+\alpha/2}\la y\ra^{2+\alpha/2}}\lesssim |x|^{-1-1/\bar c_l}. \]
\end{itemize}
Combining the above two cases, we deduce that for any $x,y$ such that $|x|\geq \max\{1,|y|\}$,
\[|K_s(x,y)|\lesssim|x|^{\alpha/2-1/\bar c_l}\la y\ra^{-1-\alpha/2} +|x|^{-2-\alpha/2}\la y\ra^{1+\alpha/2-1/\bar c_l}.\]
Therefore,
\begin{align*}
      \int_{\R^2}|K_s(x,y)|^2\idiff x\idiff y &=\int_{[-1,1]^2}|K_s(x,y)|^2\idiff x\idiff y+2 \int_{|x|>\max\{1,|y|\}}|K_s(x,y)|^2\idiff x\idiff y\\
     & \lesssim 1+ \int_0^{+\infty}\idiff y\int_{\max\{1,|y|\}}^{+\infty} \left(x^{\alpha-2/\bar c_l}\la y\ra^{-2-\alpha}+x^{-4-\alpha}\la y\ra^{2+\alpha-2/\bar c_l}\right)\idiff x\\
     & \lesssim 1+ \int_0^{+\infty} \la y\ra^{-1-2/\bar c_l}\idiff y
     <\infty,
\end{align*}
where we have used the fact that $\alpha-2/\bar{c}_l<-1$ in the second inequality. It follows that $\mtx{P}_{\rm sym}$ is Hilbert--Schmidt, hence compact.

It remains to prove that \(\mtx{K}_{\rm pr}\) is compact on $\cX_0$. To this end, we show that there exists a constant $C$ such that for any $ w\in L^2_{\rm odd}(\R)$
\[\left\|\mtx{Q}\mtx{K}_{\rm pr}\mtx{Q}^{-1}(w)\right\|_{H^1(\R)}\leq C \left\|w\right\|_{L^2(\R)}.\]
Equivalently, we show that for any $f\in \cX_0$,
\begin{equation}\label{eqt:gCLM-K-pr-compact1}
    \left\|\mtx{K}_{\rm pr}(f)\right\|_{\cX_0}\leq C \left\|f\right\|_{\cX_0}, \quad  \left\| \phi^{1/2}\mtx{K}_{\rm pr}(f)  \right\|_{\dot{H}^1(\R)}\leq C \left\|f\right\|_{\cX_0}.
\end{equation}
To prove the first estimate in \eqref{eqt:gCLM-K-pr-compact1}, we use the $L^2$ integrability of $|x|^{\alpha/2-1/\bar{c}_l}$ for $|x|\geq 1$ and \eqref{eqt:gCLM-low-coeff} to obtain
\begin{align*}
      \left\|\mtx{K}_{\rm pr}(f)\right\|_{\cX_0} &\leq |a| \left\|\mtx{V}(f) \bar{\Omega}_x \phi^{1/2}\right\|_{L^2[-1,1]}+|a| \left\|\mtx{V}(f) \bar{\Omega}_x \phi^{1/2}\right\|_{L^2(\R\backslash[-1,1])} \\
      &\lesssim \left\|\frac{\mtx{V}(f)}{x^2}\right\|_{L^{2}[-1,1]}+\left\|\frac{\mtx{V}(f)}{x}|x|^{\alpha/2-1/\bar{c}_l}\right\|_{L^2(\R\backslash[-1,1])}\\
      &\lesssim \left\|\frac{\mtx{V}(f)}{x^2}\right\|_{L^{\infty}(\R)}+\left\|\frac{\mtx{V}(f)}{x}\right\|_{L^{\infty}(\R)}\lesssim \left\|f\right\|_{\cX_0}.
\end{align*}
We further estimate the $\dot{H}^1$ seminorm as 
\[\left\| \phi^{1/2}\mtx{K}_{\rm pr}(f)  \right\|_{\dot{H}^1(\R)} = |a| \left\|\left(\mtx{V}(f) \bar{\Omega}_x \phi^{1/2}\right)_x\right\|_{L^2(\R)} \lesssim \left\|\mtx{H}_0(f) \bar{\Omega}_x \phi^{1/2}\right\|_{L^2(\R)}+\left\|\mtx{V}(f) \left(\bar{\Omega}_x \phi^{1/2}\right)_x\right\|_{L^2(\R)}.\]
We bound the first term using arguments similar to \eqref{eqt:gCLM-P-bounded}:
\begin{align*}
      \left\|\mtx{H}_0(f) \bar{\Omega}_x \phi^{1/2}\right\|_{L^2(\R)} 
      &\lesssim \left\|\frac{\mtx{H}_0(f)}{x^2}\right\|_{L^2[-1,1]}+ \left\|\frac{\mtx{H}_0(f)}{x}|x|^{\alpha/2-1/\bar{c}_l}\right\|_{L^2(\R\backslash[-1,1])}\\
      &\lesssim \left\|\mtx{H}\left(\frac{f}{x^2}\right)\right\|_{L^2(\R)}+ \left\|\mtx{H}\left(\frac{f}{x}\right)\la x\ra^{\alpha/2}\right\|_{L^2(\R)}\\
      & \lesssim \left\|\frac{f}{x^2}\right\|_{L^2(\R)}+ \left\|\left(\frac{f}{x}\right)\la x\ra^{\alpha/2}\right\|_{L^2(\R)}\lesssim \|f\|_{\cX_0},
\end{align*}
where we have used \eqref{eqt:common-hilbert-identities2} in the second inequality and \eqref{eqt:common-A2-weight} in the third inequality. By \eqref{eqt:gCLM-low-coeff} and \eqref{eqt:gCLM-low-coeff1}, the second term can be estimated as
\begin{align*}
      \left\|\mtx{V}(f) \left(\bar{\Omega}_x \phi^{1/2}\right)_x\right\|_{L^2(\R)} &\lesssim \left\|\frac{\mtx{V}(f)}{x^3}\right\|_{L^{2}[-1,1]}+\left\|\frac{\mtx{V}(f)}{x}|x|^{\alpha/2-1/\bar{c}_l-1}\right\|_{L^2(\R\backslash[-1,1])}\\
      &\lesssim \left\|\frac{\mtx{V}(f)}{x^3}\right\|_{L^{2}(\R)}+\left\|\frac{\mtx{V}(f)}{x}\right\|_{L^{\infty}(\R)}\lesssim \left\|f\right\|_{\cX_0}.
\end{align*}
 This proves \eqref{eqt:gCLM-K-pr-compact1}. We also need a uniform tail estimate. In fact, for any $R>1$, we have 
\begin{equation}\label{eqt:gCLM-K-pr-uniform-tail}
      \begin{aligned}
      \left\|\mtx{K}_{\rm pr}(f)\phi^{1/2}\right\|_{L^2(\R\backslash[-R,R])}&=|a|\left\|\frac{\mtx{V}(f)}{x} x\bar{\Omega}_x \phi^{1/2}\right\|_{L^2(\R\backslash[-R,R])}\lesssim \left\|\frac{\mtx{V}(f)}{x}\right\|_{L^{\infty}(\R)}\left\|\la x\ra^{{\alpha/2-1/\bar{c}_l}}\right\|_{L^2(\R\backslash[-R,R])}\\
       &\lesssim  R^{\alpha/2-1/\bar{c}_l+1/2}\left\|f\right\|_{\cX_0}.
      \end{aligned}
\end{equation}
The exponent is negative by \eqref{eqt:gCLM-alpha-range-intro}, so the right-hand side tends to zero as $R\to\infty$, uniformly for $f$ in the unit ball of $\cX_0$.

Now let $\{f_n\}$ be a bounded sequence in $\cX_0$ and set $g_n=\mtx{Q}(\mtx{K}_{\rm pr}(f_n))$. The estimates above show that $\{g_n\}$ is bounded in $H^1(\R)$. By the Rellich theorem and a diagonal argument, we get a subsequence which converges in $L^2([-R,R])$ for every $R>0$. The uniform tail estimate \eqref{eqt:gCLM-K-pr-uniform-tail} then makes this subsequence Cauchy in $L^2(\R)$. Thus $\mtx{Q}\mtx{K}_{\rm pr}\mtx{Q}^{-1}$ is compact on $L^2_{\rm odd}(\R)$, and hence $\mtx{K}_{\rm pr}$ is compact on $\cX_0$. This completes the proof of Lemma \ref{lem:gCLM_symmetric_compact} for $\bar{c}_l>0$.

The case $\bar c_l=0$ can be proved by a similar argument, with the Gaussian tail bound for $\bar{\Om}$ replacing the algebraic tail estimates above. 
\end{proof}

\begin{remark}\label{rem:gCLM-commutator}
To understand the compactness of $\mtx{P}_{\rm sym}$, we consider the following simplified example. On the unweighted space $L^2(\R)$, let $b$ be a suitable function and define $\mtx{T}(f):=b\mtx{H}(f)$. Since $\mtx{H}^*=-\mtx{H}$, its adjoint is $\mtx{T}^*(f)=-\mtx{H}(bf)$. Its symmetric part is therefore
\[\mtx{T}_{\rm sym}(f)=\frac12\bigl(b\mtx{H}(f)-\mtx{H}(bf)\bigr)=\frac12[b,\mtx{H}](f).\]
The commutator has the kernel representation
\[[b,\mtx{H}](f)(x)=\frac1\pi\int_\R\frac{b(x)-b(y)}{x-y}f(y)\idiff y.\]
The factor $b(x)-b(y)$ cancels the diagonal singularity, so the commutator is more regular, which yields compactness under suitable decay assumptions. We also remark that a similar fact, that commutators of operators arising from the Biot--Savart law with multiplication by suitable functions are more regular, was used by Chen and Hou \cite{chen2022stable} in their study of the Hou--Luo blowup scenario for the 3D Euler equations.
\end{remark}

\subsubsection{Linear damping estimates}\label{sec:gCLM-linear-damping}
In this subsection we show that $\mtx{L}_{\rm stb}$ defined in \eqref{eqt:gCLM-decomposition-linear} generates a strongly continuous semigroup which satisfies certain decay estimate. We define the domain of $\mtx{L}_{\rm stb}$ as 
\[ \mathcal D(\mtx{L}_{\rm stb}):=\{f\in\cX_0:\bar{V}f_x\in\cX_0\}.\]
The next lemma illustrates the construction and the decay estimate of the semigroup $\econst^{t\mtx{L}_{\rm stb}}$ in detail. 
\begin{lemma}\label{lem:Lstb-semigroup}
 $\mtx{L}_{\rm stb}$ generates a strongly continuous semigroup $\econst^{t\mtx{L}_{\rm stb}}$ on $\cX_0$. Moreover, the semigroup satisfies the decay estimate
\begin{equation}
 \left\|\econst^{t\mtx{L}_{\rm stb}}\right\|_{\cX_0\to\cX_0} \leq \econst^{-\eta_0t},
\end{equation}
where $\eta_0$ is given by Lemma \ref{lem:weight-parameters}.
\end{lemma}

\begin{proof}
Denote the pure transport operator
\[\mtx{L}_{\rm tr}(f):=-\bar{V}f_x.\] 
We first show that $\mtx{L}_{\rm tr}$ generates a strongly continuous semigroup on $\cX_0$. In the proof we identify odd functions with their restrictions to \((0,\infty)\). Let us introduce the characteristic coordinate
\[ r(x)=\int_1^x\frac{\idiff y}{\bar{V}(y)}.\]
By \eqref{eqt:gCLM-A-comparable}, \(r:(0,\infty)\to\R\) is a \(C^1\)-diffeomorphism. Denote its inverse by \(x=x(r)\). Define 
\[ \mtx{J}(f)(r) :=\sqrt{2\bar{V}(x(r))\phi(x(r))}\,f(x(r)).\]
Then $\mtx{J}$ is an isometry from \(\cX_0\) onto \(L^2(\R)\). A direct calculation yields
\begin{equation}
\mtx{J}\mtx{L}_{\rm tr}\mtx{J}^{-1}
 =-\partial_r+g_0,\quad g_0(r):=\frac12\left[\bar{V}_x(x(r))+\bar{V}(x(r))\frac{\phi_x(x(r))}{\phi(x(r))}\right].
\end{equation}
The profile and weight bounds ensure that \(g_0\in L^\infty(\R)\). It is well known that the operator \(g\mapsto-g_r\), with domain \(H^1(\R)\), is the generator of the translation group on \(L^2(\R)\). By the bounded perturbation theorem, $\mtx{J}\mtx{L}_{\rm tr}\mtx{J}^{-1}$ generates a
strongly continuous semigroup on \(L^2(\R)\) with domain \(H^1(\R)\). Hence, $\mtx{L}_{\rm tr}$ generates a strongly continuous semigroup on $\cX_0$ with domain 
\[\mtx{J}^{-1}\left(H^1(\R)\right)= \{f\in\cX_0:\bar{V}f_x\in\cX_0\}.\]

Next we consider the semigroup generated by $\mtx{L}_{\rm stb}=\mtx{L}_{\rm tr}+(\bar{U}_x-1)+\mtx{P}_{\rm skew}$. Note that \((\bar{U}_x-1)\) and $\mtx{P}_{\rm skew}$ are both bounded on $\cX_0$. The bounded perturbation theorem implies that \(\mtx{L}_{\rm stb}\) with domain $\{f\in\cX_0:\bar{V}f_x\in\cX_0\}$ generates a strongly continuous semigroup $\econst^{t\mtx{L}_{\rm stb}}$ on $\cX_0$. For $f_0\in \mathcal{D}(\mtx{L}_{\rm stb})$, $\econst^{t\mtx{L}_{\rm stb}}(f_0)$ gives the solution to the linear equation 
\[ f_t=\mtx{L}_{\rm stb}(f)=-\bar{V}f_x+(\bar{U}_x-1)f+\mtx{P}_{\rm skew}(f)\]
with the initial data $f_0$. We perform an energy estimate to obtain
\[\frac{1}{2}\frac{\idiff {}}{\idiff t}\left\|f\right\|_{\cX_0}^2=\la -\bar{V}f_x+(\bar{U}_x-1)f,f \ra_{\cX_0}=\int_\R \left((\bar{U}_x-1)+\frac{\bar{V}_x}2+\frac{\bar{V}}{2}\frac{\phi_x}\phi\right) |f|^2\phi\idiff x \leq-\eta_0\left\|f\right\|_{\cX_0}^2,\]
where we used \(\la \mtx{P}_{\rm skew}(f),f\ra_{\cX_0}=0\). The second equality follows from integration by parts, and the last inequality follows from Lemma~\ref{lem:weight-parameters}. Thus,
\[\left\|\econst^{t\mtx{L}_{\rm stb}}(f_0)\right\|_{\cX_0}\leq \econst^{-\eta_0t}\left\|f_0\right\|_{\cX_0} \] 
for any $f_0\in \mathcal{D}(\mtx{L}_{\rm stb})$. The estimate extends to $\cX_0$ by a density argument. Our proof is thus completed.
\end{proof}

\subsubsection{Linear finite-codimensional stability}
To apply spectral theory, we need to complexify the real Hilbert space considered before. More specifically, let \(\cX_0^{\CC}:=\cX_0\oplus \iunit\cX_0\) be the complexification of \(\cX_0\), equipped with its canonical Hermitian inner product. We use the same symbol for a real operator and its canonical complexification; for bounded operators this preserves norms, compactness, and adjoints.

Let \(\mtx{G}\) be a generator of a strongly continuous semigroup on a complex Banach space \(X\). Define
\begin{align*}
 s(\mtx{G})&:=\sup\{\operatorname{Re} z:z\in\sigma(\mtx{G})\},\\
 s_0(\mtx{G})&:=\inf\left\{\theta\in\R:\exists M_\theta\geq1,\
       \left\|\econst^{t\mtx{G}}\right\|\leq M_\theta \econst^{\theta t}\ \text{for all }t\geq0\right\},\\
 \left\|\mtx{T}\right\|_{\rm ess}&:=\inf\left\{\left\|\mtx{T}-\mtx{C}\right\|:\mtx{C}:X\to X\ \text{is compact}\right\},\\
 s_{\rm ess}(\mtx{G})&:=\inf_{t>0}t^{-1}\log\left\|\econst^{t\mtx{G}}\right\|_{\rm ess}.
\end{align*}
These are the spectral bound, growth bound, essential norm, and essential growth bound, respectively. We shall use the following standard result without proof.
\begin{proposition}\cite[Corollary~IV.2.11]{engel2000one}
\label{prop:semigroup-spectral-bound}
If \(\mtx{G}\) generates a strongly continuous semigroup, then
\[ s_0(\mtx{G})=\max\{s_{\rm ess}(\mtx{G}),s(\mtx{G})\}.\]
Moreover, for every \(\theta>s_{\rm ess}(\mtx{G})\), the set \(\sigma(\mtx{G})\cap\{\operatorname{Re} z\geq\theta\}\) is finite and its spectral projection has finite rank.
\end{proposition}
For the operators in \eqref{eqt:gCLM-decomposition-linear}, Lemmas~\ref{lem:gCLM_symmetric_compact} and \ref{lem:Lstb-semigroup} show that \(\mtx{L}_{\rm stb}\) generates an exponentially decaying semigroup on \(\cX_0^{\CC}\) and \(\mtx{K}\) is compact. Hence \(\mtx{L}=\mtx{L}_{\rm stb}+\mtx{K}\) generates a strongly continuous semigroup with \(\mathcal D(\mtx{L})=\mathcal D(\mtx{L}_{\rm stb})\). In particular, since $s_{\rm ess}$ is invariant under compact perturbations (\cite[Proposition~IV.2.12]{engel2000one}), we have
\[ s_{\rm ess}(\mtx{L})=s_{\rm ess}(\mtx{L}_{\rm stb})\leq-\eta_0.\]
Let \(\lambda \in (0,\eta_0)\) be a fixed positive number, and define
\[\Sigma_{\lambda}:=\sigma(\mtx{L})\cap\{z\in\CC:\operatorname{Re} z>-\lambda\}.\]
By Proposition~\ref{prop:semigroup-spectral-bound}, \(\Sigma_{\lambda}\) is a finite set of isolated eigenvalues with finite algebraic multiplicity.
If \(\Sigma_{\lambda}\neq\varnothing\), choose a finite union \(\Gamma\) of positively oriented disjoint circles in \(\CC\setminus\sigma(\mtx{L})\) enclosing exactly the eigenvalues in \(\Sigma_{\lambda}\). The associated Riesz projection is defined by
\[\mtx{\Pi}_{\lambda}:=\frac{1}{2\pi \iunit}\int_\Gamma(z-\mtx{L})^{-1}\idiff z,\quad \mtx{\Pi}_s:=\mtx{I}-\mtx{\Pi}_{\lambda}. \]
If \(\Sigma_{\lambda}=\varnothing\), set \(\mtx{\Pi}_{\lambda}=0\) and \(\mtx{\Pi}_s=\mtx{I}\).
Then $\mtx{\Pi}_{\lambda}$ and $\mtx{\Pi}_s$ are bounded projections on \(\cX_0^{\CC}\) that commute with \(\mtx{L}\) and \(\econst^{t\mtx{L}}\). In particular,
\[ \mtx{\Pi}_j^2=\mtx{\Pi}_j\quad(j=\lambda,s),\quad \mtx{\Pi}_{\lambda}\mtx{\Pi}_s=\mtx{\Pi}_s\mtx{\Pi}_{\lambda}=0,\quad \mtx{\Pi}_{\lambda}+\mtx{\Pi}_s=\mtx{I}.\]
Define
\[ E_j^{\CC}:=\mtx{\Pi}_j\cX_0^{\CC}, \quad E_j:=E_j^{\CC}\cap\cX_0, \quad j=\lambda,s. \]
These subspaces are closed and invariant under \(\mtx{L}\) and \(\econst^{t\mtx{L}}\). The restrictions of \(\mtx{L}\) to
\(E_{\lambda}^{\CC}\) and \(E_s^{\CC}\) have spectra \(\Sigma_{\lambda}\) and \(\sigma(\mtx{L})\setminus\Sigma_{\lambda}\), respectively
\cite[Proposition~IV.1.16]{engel2000one}. Moreover, \(E_{\lambda}^{\CC}\) is finite dimensional and \(\mtx{L}|_{E_{\lambda}^{\CC}}\) is bounded. 

Notice that $\mtx{L}$ commutes with complex conjugation. As a result, for any $z\in\rho(\mtx{L})$ and $g\in\cX_0^{\CC}$,
\[ \overline{(z-\mtx{L})^{-1}g} =(\bar{z}-\mtx{L})^{-1}\bar{g}.\]
It follows that \(\sigma(\mtx{L})\), and hence \(\Sigma_{\lambda}\), is invariant under complex conjugation. When \(\Sigma_{\lambda}\neq\varnothing\), we may therefore choose the contour \(\Gamma\) to be invariant under conjugation as a set. Since conjugation reverses its orientation, direct computation gives
\[ \overline{\mtx{\Pi}_{\lambda}g}=-\frac{1}{2\pi\iunit}\int_{\overline\Gamma}(z-\mtx{L})^{-1}\bar{g}\idiff z=\frac{1}{2\pi\iunit}\int_\Gamma(z-\mtx{L})^{-1}\bar{g}\idiff z=\mtx{\Pi}_{\lambda}\bar{g}.\]
Thus \(\mtx{\Pi}_{\lambda}\) and \(\mtx{\Pi}_s=\mtx{I}-\mtx{\Pi}_{\lambda}\) commute with complex conjugation. For any \(g\in\cX_0\) and $j=\lambda,s$, we have $\overline{\mtx{\Pi}_jg}=\mtx{\Pi}_j\bar{g}=\mtx{\Pi}_jg$, which implies \(\mtx{\Pi}_j\cX_0\subset\cX_0\). As a result,
\[\cX_0=E_{\lambda}\oplus E_s,\quad E_j=\mtx{\Pi}_j\cX_0 \quad j=\lambda,s.\]
We remark that the projections \(\mtx{\Pi}_{\lambda},\mtx{\Pi}_s\) need not be orthogonal, and the decomposition above is a topological direct sum. The next lemma gives the decay estimate of the semigroup $\econst^{t\mtx{L}}$ on $E_s$.

\begin{lemma}
\label{lem:linear-package}
For every \(\eta_s\in(0,\lambda)\), there is a constant
\(C_{\eta_s}>0\) such that
\begin{equation}\label{eqt:gCLM-stable-semigroup-proof}
 \left\|\econst^{t\mtx{L}}\mtx{\Pi}_s\right\|_{\cX_0\to\cX_0}
 \leq C_{\eta_s}\econst^{-\eta_s t},\quad t\geq0.
\end{equation}
\end{lemma}

\begin{proof}
Note that \(E_s^{\CC}\) is closed and invariant under \(\econst^{t\mtx{L}}\). By \cite[Paragraphs~I.5.12 and~II.2.3]{engel2000one}, the restricted semigroup \(\econst^{t\mtx{L}}|_{E_s^{\CC}}\) is strongly continuous with generator
\[ \mtx{L}_s:=\mtx{L}|_{E_s^{\CC}},\quad \mathcal D(\mtx{L}_s)=\mathcal D(\mtx{L})\cap E_s^{\CC}.\]
Since \(E_s^{\CC}=\ker\mtx{\Pi}_{\lambda}\) and \(\mtx{\Pi}_{\lambda}\) is a finite-rank spectral projection, the argument in \cite[proof of Theorem~V.3.1]{engel2000one} gives
\[ s_{\rm ess}(\mtx{L}_s)=s_{\rm ess}(\mtx{L})\leq-\eta_0.\]
The spectral decomposition above gives \(s(\mtx{L}_s)\leq-\lambda\), so Proposition~\ref{prop:semigroup-spectral-bound} yields
\[ s_0(\mtx{L}_s)=\max\{s_{\rm ess}(\mtx{L}_s),s(\mtx{L}_s)\}\leq-\lambda.\]
For every \(\eta_s\in(0,\lambda)\), by the definition of \(s_0\) we have 
\[\left\|\econst^{t\mtx{L}}\right\|_{E_s^{\CC}\to E_s^{\CC}} \leq M_{\eta_s}\econst^{-\eta_st}\]
 for some \(M_{\eta_s}>0\). Therefore, 
\[ \left\|\econst^{t\mtx{L}}\mtx{\Pi}_s\right\|_{\cX_0^{\CC}\to\cX_0^{\CC}} \leq M_{\eta_s}\econst^{-\eta_st}\, \|\mtx{\Pi}_s\|_{\cX_0^{\CC}\to\cX_0^{\CC}}. \]
Restricting the above estimate to \(\cX_0\) proves \eqref{eqt:gCLM-stable-semigroup-proof}.
\end{proof}
\subsection{Nonlinear finite-codimensional stability}
In this subsection, we prove nonlinear finite-codimen\-sional stability of \eqref{eqt:gCLM-perturbation} via a shooting method and establish asymptotically self-similar blowups of the gCLM equation \eqref{eqt:gCLM-model} from $C_c^{\infty}$ initial data for the non-expanding case $\bar{c}_l\geq 0$. 
\subsubsection{Nonlinear energy estimates}
We first perform an energy estimate for the nonlinear equation \eqref{eqt:gCLM-perturbation}. Similar to \cite{chen2021finite,chen2022asymptotically}, to control the nonlinear terms in the energy estimate, we need to obtain some weighted $H^1$ estimates in the next lemma.
\begin{lemma}\label{lem:nonlinear_estimate}
For the nonlinear operator $\mtx{N}$ defined in \eqref{eqt:gCLM-perturbation_nonlinear}, there exists a constant $C$ such that
\begin{equation}\label{eqt:gCLM-nonlinear-low}
\left\|\mtx{N}(f)\right\|_{\cX_0}\leq C\left\|f\right\|_{\cX_0}\left\|f\right\|_{\cX_1},
\end{equation}
and 
\begin{equation}\label{eqt:gCLM-nonlinear-high}
\la \mtx{N}(f),f\ra_{{{\cX}}_1}\leq C \left\|f\right\|_{\cX_1}^3.
\end{equation}
\end{lemma}
\begin{proof}
We use \eqref{eqt:gCLM-high-coeff} to deduce
\begin{align*}
      \left\|\mtx{N}(f)\right\|_{\cX_0}&=\left\| \left(-a\mtx{V}(f)f_x+\mtx{H}_0(f)f\right)\phi^{1/2}\right\|_{L^2(\R)}\\
      &\lesssim\left\|\frac{\mtx{V}(f)}{x}\right\|_{L^{\infty}(\R)}\left\|xf_x\phi^{1/2}\right\|_{L^2(\R)}+\left\|\mtx{H}_0(f)\right\|_{L^{\infty}(\R)}\left\|f\phi^{1/2}\right\|_{L^2(\R)}\\
      &\lesssim \left\|f\right\|_{\cX_0}\left\|f\right\|_{\cX_1},
\end{align*}
which implies \eqref{eqt:gCLM-nonlinear-low}. To prove \eqref{eqt:gCLM-nonlinear-high}, it suffices to establish the $\dot{\cX}_1$ estimate. Integration by parts gives
\begin{align*}
 \la \mtx{N}(f),f\ra_{\dot{\cX}_1}&=\int_\R\bigl(-a\mtx{H}_0(f)f_x-a\mtx{V}(f)f_{xx}+\mtx{H}_0(f)_xf+\mtx{H}_0(f)f_x\bigr)
 f_x\psi\idiff x\\
 &=\int_\R\left(\left(1-\frac a2\right)\mtx{H}_0(f)+\frac{ax\psi_x}{2\psi}\frac{\mtx{V}(f)}x\right)f_x^2\psi\idiff x+\int_\R ff_x\mtx{H}(f_x)\psi\idiff x.
\end{align*}
 The first term above is bounded by $\|f\|_{\cX_1}^3$ because $\|\mtx{H}_0(f)\|_{L^{\infty}(\R)}$ and $\|\mtx{V}(f)/x\|_{L^{\infty}(\R)}$ are both bounded by $\|f\|_{\cX_1}$. For the second term, we use \eqref{eqt:common-hilbert-identities} and \eqref{eqt:common-A2-weight} to obtain
 \begin{equation}\label{eqt:gCLM-weighted-hilbert-derivative}
      \begin{aligned}
   \left\|\mtx{H}(f_x)\psi^{1/2}\right\|_{L^2(\R)}&\lesssim \left\| \frac{\mtx{H}(f_x)}{x} \right\|_{L^2(-1,1)}+ \left\| x\mtx{H}(f_x)\phi^{1/2} \right\|_{L^2(\R\backslash(-1,1))}\\
   &\lesssim \left\|\mtx{H}\left(\frac{f_x}{x}\right) \right\|_{L^2(\R)}+ \left\|\mtx{H}(xf_x)\la x\ra^{\alpha/2}\right\|_{L^2(\R)}\\
   &\lesssim \left\|\frac{f_x}{x} \right\|_{L^2(\R)}+ \left\|xf_x\la x\ra^{\alpha/2}\right\|_{L^2(\R)}\lesssim  \left\|f_x\psi^{1/2}\right\|_{L^2(\R)}.
      \end{aligned}
 \end{equation}
Therefore, 
 \[ \left|\int_\R ff_x\mtx{H}(f_x)\psi\idiff x\right| \leq \|f\|_{L^\infty(\R)}  \left\|f_x\psi^{1/2}\right\|_{L^2(\R)} \left\|\mtx{H}(f_x)\psi^{1/2}\right\|_{L^2(\R)}\lesssim \left\|f\right\|_{\cX_1}^3, \]
where we have used \eqref{eqt:gCLM-high-coeff} to bound $\|f\|_{L^{\infty}(\R)}$. We thus deduce \eqref{eqt:gCLM-nonlinear-high}.
\end{proof}

With the established linear and nonlinear estimates, we can now derive the nonlinear energy estimate for the perturbation equation \eqref{eqt:gCLM-perturbation} in the next lemma.
\begin{lemma}[Energy estimate]\label{lem:high-energy}
      For any $\eta_1\in (0,\eta_0)$, there exists a constant $C=C_{\eta_1}$ such that
\begin{equation}\label{eqt:gCLM-high-differential}
 \frac12\frac{\idiff{}}{\idiff t}\left\|f\right\|_{\cX_1}^2 +\eta_1\left\|f\right\|_{\cX_1}^2 \leq C\left\|f\right\|_{\cX_0}^2+C\left\|f\right\|_{\cX_1}^3.
\end{equation}
Consequently,
\begin{equation}\label{eqt:gCLM-high-integral}
 \left\|f(t)\right\|_{\cX_1} \leq C\econst^{-\eta_1 t}\left\|f(0)\right\|_{\cX_1} +C\int_0^t \econst^{-\eta_1(t-s)}\left(\left\|f(s)\right\|_{\cX_0}+\left\|f(s)\right\|_{\cX_1}^2\right)\idiff s .
\end{equation}
\end{lemma}

\begin{proof}
Write 
\[ \frac12\frac{\idiff{}}{\idiff t}\left\|f\right\|_{\cX_1}^2 =\la \mtx{L}(f),f\ra_{\cX_1}+\la \mtx{N}(f),f\ra_{\cX_1}.\]
By Lemma \ref{lem:nonlinear_estimate} it suffices to establish the linear estimate
\begin{equation}\label{eqt:gCLM-nonlinear-energy-estimate-linear}
\la \mtx{L}(f),f\ra_{\cX_1}\leq -\eta_1\left\|f\right\|_{\cX_1}^2+ C\left\|f\right\|_{\cX_0}^2.
\end{equation}

Recall that $\mtx{L}=\mtx{L}_{\rm stb}+\mtx{K}$, where $\mtx{L}_{\rm stb}$ satisfies
\[ \la \mtx{L}_{\rm stb}(f),f\ra_{\cX_0} \leq-\eta_0\left\|f\right\|_{\cX_0}^2 \]
and $\mtx{K}$ is compact (and thus bounded) on $\cX_0$. We immediately have
\[ \la \mtx{L}(f),f\ra_{\cX_0}\leq C\left\|f\right\|_{\cX_0}^2. \]

We then consider the weighted $H^1$ estimate. Direct calculation gives
\begin{align*}
 \la \mtx{L}(f),f\ra_{{\dot{\cX}_1}} &=\int_\R\Bigl(-\bar{V}f_{xx}-\bar{V}_xf_x+(\bar{U}_x-1)f_x+\bar{U}_{xx}f\\
 &\quad-a\mtx{H}_0(f)\bar\Om_x-a\mtx{V}(f)\bar\Om_{xx}+\mtx{H}_0(f)_x\bar\Om+\mtx{H}_0(f)\bar\Om_x\Bigr)f_x\psi\idiff x\\
 &=\int_\R\left(-\bar{V}f_{xx}+(\bar{U}_x-1-\bar{V}_x)f_x\right)f_x\psi\idiff x+\int_\R\mtx{H}(f_x)\bar\Om f_x\psi\idiff x\\
 &\quad+\int_\R\left(\bar{U}_{xx}f+(1-a)\mtx{H}_0(f)\bar\Om_x-a\mtx{V}(f)\bar\Om_{xx}\right)f_x\psi\idiff x\\
 &=:I_1+I_2+I_3.
\end{align*}
For $I_1$, we use integration by parts to compute
\[ I_1=\int_\R\left((\bar{U}_x-1)+\left(\frac{x\psi_x}{2\psi}\right)\frac{\bar{V}}{x}-\frac{\bar{V}_x}{2}\right)f_x^2\psi\idiff x \leq-\eta_0\left\|f\right\|_{\dot{\cX}_1}^2.\]
We have used Lemma~\ref{lem:weight-parameters} in the above estimate. For $I_2$, we have
 \[I_2=\int_\R x\mtx{H}(f_x)\bar\Om xf_x\phi\idiff x=\int_\R\mtx{H}(xf_x)\bar\Om xf_x\phi\idiff x=\la \mtx{P}(xf_x),xf_x\ra_{\cX_0}=\la \mtx{P}_{\rm sym}(xf_x),xf_x\ra_{\cX_0}.\]
Here we have used \eqref{eqt:common-hilbert-identities} in the second equality, and $\mtx{P}$ is given by \eqref{eqt:gCLM-decomposition-linear1}. By Lemma~\ref{lem:gCLM_symmetric_compact}, $\mtx{P}_{\rm sym}$ is compact on $\cX_0$. Since $f\mapsto xf_x$ is bounded from $\cX_1$ to $\cX_0$, the map $f\mapsto\mtx{P}_{\rm sym}(xf_x)$ is compact from $\cX_1$ to $\cX_0$. Applying Lemma~\ref{lem:compactness-estimate} with $\mathcal{A}=\cX_1$ and $\mathcal{B}=\cX_0$, we obtain that, for any $\varepsilon>0$, there exists a constant $C_{\varepsilon}>0$ such that
\[ \left\|\mtx{P}_{\rm sym}(xf_x)\right\|_{\cX_0}\leq\varepsilon\left\|f\right\|_{\cX_1}+C_\varepsilon\left\|f\right\|_{\cX_0}\]
for all $f\in\cX_1$. We therefore have
\[ |I_2|\leq\left\|\mtx{P}_{\rm sym}(xf_x)\right\|_{\cX_0}\left\|f\right\|_{\dot{\cX}_1} \leq(\varepsilon+\delta C_\varepsilon)\left\|f\right\|_{{\cX}_1}^2 +\frac{C_\varepsilon}{4\delta}\left\|f\right\|_{\cX_0}^2.\]
Choosing $\delta=\varepsilon/C_\varepsilon$, we obtain
\[ |I_2|\leq 2\varepsilon\left\|f\right\|_{\dot{\cX}_1}^2 +C_{\varepsilon}\left\|f\right\|_{\cX_0}^2 \leq 2\varepsilon\left\|f\right\|_{\cX_1}^2+C_{\varepsilon}\left\|f\right\|_{\cX_0}^2. \]
For the cross term $I_3$, by the Cauchy--Schwarz inequality we have
\begin{align*}
 |I_3| &=\left|\int_\R\left(\bar{U}_{xx}f+(1-a)\mtx{H}_0(f)\bar\Om_x -a\mtx{V}(f)\bar\Om_{xx}\right)f_x\psi\idiff x\right|\\
 &\lesssim\left(\left\|x\bar{U}_{xx}\right\|_{L^\infty(\R)}\left\|f\right\|_{\cX_0} +\left\|x\bar\Om_x\mtx{H}_0(f)\right\|_{\cX_0} +\left\|x\bar\Om_{xx}\mtx{V}(f)\right\|_{\cX_0}\right) \left\|f_x\psi^{1/2}\right\|_{L^2(\R)}.
\end{align*}
One can show that the three terms in the bracket are bounded by $\|f\|_{\cX_0}$ by the same arguments as in Lemma~\ref{lem:gCLM_symmetric_compact}. As a result, we have 
\[|I_3|\leq \varepsilon\left\|f_x\psi^{1/2}\right\|_{L^2(\R)}^2+C_{\varepsilon}\left\|f\right\|_{\cX_0}^2\leq \varepsilon\left\|f\right\|_{\cX_1}^2+C_{\varepsilon}\left\|f\right\|_{\cX_0}^2.\] 
Combining the estimates for $I_1$, $I_2$, and $I_3$, we get
\[ \la \mtx{L}(f),f\ra_{\dot{\cX}_1} \leq(-\eta_0+3\varepsilon)\left\|f\right\|_{\cX_1}^2+C\left\|f\right\|_{\cX_0}^2.\]
Taking $\varepsilon<(\eta_0-\eta_1)/3$ yields \eqref{eqt:gCLM-nonlinear-energy-estimate-linear}. \eqref{eqt:gCLM-high-differential} is thus proved.
\eqref{eqt:gCLM-high-integral} follows directly from Grönwall's inequality. This completes the proof.
\end{proof}
\subsubsection{Finite-dimensional shooting}
In this subsection we establish nonlinear finite-codimensional stability of the perturbation equation \eqref{eqt:gCLM-perturbation} using the shooting method. Specifically, we show that for any initial perturbation $\xi\in\cX_1$ which is sufficiently small, one can add a small correction in finitely many fixed $C_c^\infty(\R)$ directions so that the solution of \eqref{eqt:gCLM-perturbation} with the corrected initial data converges to zero in $\cX_1$. 

Recall that $\eta_0>0$ is the damping constant in Lemma~\ref{lem:weight-parameters}, and $\lambda\in(0,\eta_0)$ is the parameter chosen in Section~\ref{sec:linear} to separate the eigenvalues with $\operatorname{Re}z>-\lambda$ from the remaining spectrum. For the shooting argument, we need a lower bound for the quadratic form of the linear operator on $E_\lambda$. The following lemma obtains this bound by introducing an equivalent inner product \(\la\cdot,\cdot\ra_*\).

\begin{lemma}[Equivalent norm on the finite-dimensional subspace]
\label{lem:unstable-adapted-norm}
Let \(\mtx{L}_{\lambda}:=\mtx{L}|_{E_{\lambda}}\).  There is an inner product \(\la\cdot,\cdot\ra_*\) on \(E_{\lambda}\) which is equivalent to $\la\cdot,\cdot\ra_{\cX_0}$ such that
\begin{equation}
 \la \mtx{L}_{\lambda}(g),g\ra_*\geq-\lambda\left\|g\right\|_*^2,
 \quad \text{for any } g\in E_{\lambda}.
\end{equation}
\end{lemma}
\begin{proof}
By the definition of \(\Sigma_{\lambda}\), the complex spectrum of \(\mtx{L}_{\lambda}\) lies in \(\{\operatorname{Re} z>-\lambda\}\).  Hence
\(\mtx{B}_{\lambda}:=\mtx{L}_{\lambda}+\lambda \mtx{I}\) has spectrum in the open right half-plane.  Given any real Euclidean inner product \(\la\cdot,\cdot\ra_0\) on \(E_{\lambda}\), set
\[ \la f,g \ra_*:=\int_0^\infty
 \left\langle \econst^{-t\mtx{B}_{\lambda}}(f),\econst^{-t\mtx{B}_{\lambda}}(g)\right\rangle_0\idiff t.\]
Since \(E_{\lambda}\) is finite dimensional and the spectrum of \(\mtx{B}_{\lambda}\) lies in the open right half-plane, \(\econst^{-t\mtx{B}_{\lambda}}\) decays exponentially. Thus the integral converges and defines a positive-definite inner product. Moreover,
\[2\la \mtx{B}_{\lambda}(g),g\ra_*
=-\int_0^{\infty}\frac{\idiff{}}{\idiff t}
\left\|\econst^{-t\mtx{B}_{\lambda}}g\right\|_{0}^2\idiff t
=\left\|g\right\|_0^2.\]
Therefore
\[ \la \mtx{L}_{\lambda}(g),g\ra_* =\la \mtx{B}_{\lambda}(g)-\lambda g,g\ra_*=-\lambda\left\|g\right\|_*^2+\frac12\left\|g\right\|_0^2 \geq-\lambda\left\|g\right\|_*^2.\]
This completes the proof.
\end{proof}

We also need the following bootstrap lemma to improve the bound on the perturbation. In what follows, we take $\eta:=0.9\lambda$.
\begin{lemma}[Bootstrap]\label{lem:bootstrap}
There exist \(\nu_1,\rho_1,\delta_1>0\) such that for any $\delta\in (0,\delta_1)$, if $f$ is a solution to \eqref{eqt:gCLM-perturbation} whose initial data $f(0)$ satisfy
\[\left\|f(0)\right\|_{\cX_1}\leq \nu_1\delta , \]
and 
\begin{equation}
 \left\|f(t)\right\|_{\cX_1} \leq\delta \econst^{-\eta t},\quad
 \left\|\mtx{\Pi}_{\lambda}(f(t))\right\|_*\leq \rho_1 \delta \econst^{-2\eta t}
\end{equation}
on the time interval $[0,T^*]$, then 
\begin{equation}\label{eqt:gCLM-bootstrap}
      \left\|f(t)\right\|_{\cX_1} \leq \frac{\delta}{2} \econst^{-\eta t}
\end{equation}
for any $t\in [0,T^*]$.
\end{lemma}
\begin{proof}
Write \(f_{\lambda}=\mtx{\Pi}_{\lambda}(f)\) and \(f_s=\mtx{\Pi}_s(f)\). By Duhamel's formula, we have 
\[f_s(t)=\mtx{\Pi}_sf(t)=\mtx{\Pi}_s\econst^{t\mtx{L}}f(0)+\int_0^t \mtx{\Pi}_s\econst^{(t-s)\mtx{L}}\mtx{N}(f(s))\idiff s. \]
Hence for any $t\in [0,T^*]$
\begin{align*}
      \left\|f_s(t)\right\|_{\cX_0} & \leq \left\|\mtx{\Pi}_s\econst^{t\mtx{L}}\right\|_{\cX_0\to\cX_0}\left\|f(0)\right\|_{\cX_0}+\int_0^t \left\|\mtx{\Pi}_s\econst^{(t-s)\mtx{L}}\right\|_{\cX_0\to\cX_0}\left\|\mtx{N}(f(s))\right\|_{\cX_0} \idiff s\\  
      & \leq C \econst^{-\eta t}\left\|f(0)\right\|_{\cX_0}+ C\int_0^t \econst^{-\eta (t-s)} \left\|f(s)\right\|_{\cX_1}^2 \idiff s\\
      & \leq C\nu_1 \delta\econst^{-\eta t}+C\delta^2\int_0^t \econst^{-\eta (t-s)} \econst^{-2\eta s} \idiff s\\
      & \leq C(\nu_1+\delta_1)\delta \econst^{-\eta t}. 
\end{align*}
 We have chosen $\eta_s=\eta$ in Lemma \ref{lem:linear-package} to obtain the decay estimate of the semigroup $\mtx{\Pi}_s\econst^{t\mtx{L}}$. The above inequality implies 
\[ \left\|f(t)\right\|_{\cX_0}\leq \left\|f_s(t)\right\|_{\cX_0}+\left\|f_{\lambda}(t)\right\|_{\cX_0}\leq C(\nu_1+\rho_1+\delta_1)\delta \econst^{-\eta t}. \]
We then use the energy estimate established in Lemma \ref{lem:high-energy} to improve the bound of $\|f\|_{\cX_1}$. Choose $\eta_1\in (\eta,\eta_0)$ in \eqref{eqt:gCLM-high-integral}, then
\begin{align*}
\left\|f(t)\right\|_{\cX_1}&\leq C \econst^{-\eta_1t}\left\|f(0)\right\|_{\cX_1}+C\int_0^t \econst^{-\eta_1(t-s)}\left(\left\|f(s)\right\|_{\cX_0}+\left\|f(s)\right\|_{\cX_1}^2\right)\idiff s\\
&\leq C\econst^{-\eta_1 t}\nu_1 \delta+ C\int_0^t \econst^{-\eta_1(t-s)}\left((\nu_1+\rho_1+\delta)\delta \econst^{-\eta s}+\delta^2\econst^{-2\eta s}\right)\idiff s\\
& \leq C(\nu_1+\rho_1+\delta_1)\delta \econst^{-\eta t}.
\end{align*}
We choose $\nu_1$, $\rho_1$ and $\delta_1$ such that $C(\nu_1+\rho_1+\delta_1)\leq 1/2$. This proves \eqref{eqt:gCLM-bootstrap}.

\end{proof}

To preserve the smoothness and compact support of the initial data for the dynamic rescaling equation \eqref{eqt:gCLM-rescaled}, we need the correction directions to lie in $C_c^\infty(\R)$. To this end, we define
\[ \cX_c:=\left\{g\in C_c^\infty(\R): \text{$g$ is odd and compactly supported away from the origin}\right\}, \]
and let $N:=\dim_{\R}E_{\lambda}$. The next lemma constructs an
$N$-dimensional space of correction directions in $\cX_c$.

\begin{lemma}[Choice of correction directions]
\label{lem:correction-directions}
There exist real-valued functions \(\Phi_1,\ldots,\Phi_N\in\cX_c\) and a linear embedding
 \[\mtx{J}_{\lambda}:E_{\lambda}\longrightarrow \operatorname{span}_{\R}\{\Phi_1,\ldots,\Phi_N\}\]
such that
\begin{equation}\label{eqt:gCLM-right-inverse}
\quad \mtx{\Pi}_{\lambda} \mtx{J}_{\lambda}(g)=g
\end{equation}
for any $g\in E_{\lambda}$. Moreover, \(\mtx{J}_{\lambda}\) is bounded from \(E_{\lambda}\), equipped with
\(\|\cdot\|_*\), to \(\cX_1\).
\end{lemma}
\begin{proof}
If \(E_{\lambda}=\{0\}\), we take \(N=0\) and \(\mtx{J}_{\lambda}=0\). By cutoff and mollification, $\cX_c$ is dense in both $\cX_0$ and $\cX_1$. Since \(\mtx{\Pi}_{\lambda}:\cX_0\to E_{\lambda}\) is a bounded real projection, \(\mtx{\Pi}_{\lambda}(\cX_c)\) is a dense subspace of \(E_{\lambda}\). This implies \(\mtx{\Pi}_{\lambda}(\cX_c)=E_{\lambda}\) because \(E_{\lambda}\) is finite dimensional. We therefore choose \(\Phi_1,\ldots,\Phi_N\in\cX_c\) such that \(\mtx{\Pi}_{\lambda}(\Phi_1),\ldots,\mtx{\Pi}_{\lambda}(\Phi_N)\) form a basis of \(E_{\lambda}\). Define
\[ \mtx{J}_{\lambda}\left(\sum_{j=1}^N a_j\mtx{\Pi}_{\lambda}(\Phi_j)\right) :=\sum_{j=1}^N a_j\Phi_j.\]
Then \(\mtx{J}_{\lambda}\) is real-linear and satisfies \eqref{eqt:gCLM-right-inverse}. The boundedness of $\mtx{J}_{\lambda}$ follows because its domain is finite dimensional.
\end{proof}

With the above preparations, we are ready to prove nonlinear finite-codimensional stability of \eqref{eqt:gCLM-perturbation} in the next theorem.
\begin{theorem}\label{thm:brouwer}
There exist \(\delta_0,\nu,\rho>0\) such that for any $\delta\in (0,\delta_0)$ and any initial data $\xi\in \cX_1$ with $\|\xi\|_{\cX_1}\leq \nu \delta$, there exists \(h\in E_{\lambda}\) with \(\|h\|_*\leq\rho\delta\), such that starting from the initial data
\[ f_h(0):=\xi+\mtx{J}_{\lambda}(h-\mtx{\Pi}_{\lambda}(\xi)),\]
the solution $f_h(t)$ to \eqref{eqt:gCLM-perturbation} satisfies the decay estimates  
\begin{equation}\label{eqt:gCLM-nonlinear-decay-estimate1}
 \left\|f_h(t)\right\|_{\cX_1}\leq\delta \econst^{-\eta t},\end{equation}
 \begin{equation}\label{eqt:gCLM-nonlinear-decay-estimate2}
 \left\|\mtx{\Pi}_{\lambda}(f_h(t))\right\|_*\leq\rho\delta \econst^{-2\eta t}.
\end{equation}
\end{theorem}

\begin{proof}
Let \(\lambda,\eta,\delta_1,\nu_1,\rho_1\) be the parameters fixed above. 
If \(E_{\lambda}=\{0\}\), choose
\(0<\nu\leq\min\{\nu_1,1/2\}\), \(0<\delta_0\leq\delta_1\), and any
\(\rho>0\), and take \(h=0\). Then \(\mtx{\Pi}_{\lambda}=0\), so
\eqref{eqt:gCLM-nonlinear-decay-estimate2} is automatic, while
Lemma~\ref{lem:bootstrap} and the standard continuity argument close
\eqref{eqt:gCLM-nonlinear-decay-estimate1}. Henceforth assume \(E_{\lambda}\neq\{0\}\).
We first justify the choice of $\nu$, $\rho$ and $\delta_0$. Assume that $\|\xi\|_{\cX_1}\leq \nu \delta$ and \(\|h\|_*\leq\rho\delta\). A direct calculation gives 
\[\left\|f_h(0)\right\|_{\cX_1}\leq \left\|\xi\right\|_{\cX_1}+ \left\|\mtx{J}_{\lambda}\right\|\left( \left\|h\right\|_*+\left\|\mtx{\Pi}_{\lambda}(\xi)\right\|_*\right)\leq C \left(\left\|\xi\right\|_{\cX_1}+ \left\|h\right\|_*\right)\leq C(\nu+\rho)\delta.\]
We choose $\nu$ and $\rho$ small enough so that $C(\nu+\rho)\leq \nu_1$ and $\rho<\rho_1$. 

We argue that there exists $\delta_0=\delta_0(\rho)$ such that Theorem~\ref{thm:brouwer} holds. 
Let $\delta<\delta_1$ and $\mathbb{B}_{\lambda}(\rho\delta)=\{h\in E_{\lambda}: \|h\|_*\leq \rho\delta\}$. Let $\xi\in \cX_1$ with $\|\xi\|_{\cX_1}\leq \nu \delta$. Suppose $h\in \mathbb{B}_{\lambda}(\rho\delta)$ such that \eqref{eqt:gCLM-nonlinear-decay-estimate1} or \eqref{eqt:gCLM-nonlinear-decay-estimate2} fails at some time. Note that $\|\mtx{\Pi}_{\lambda} f_h(0)\|_*=\|h\|_*\leq \rho\delta$. Therefore, for every such  $h\in \mathbb{B}_{\lambda}(\rho\delta)$, there exists $t_*\geq 0$ such that equality holds for \eqref{eqt:gCLM-nonlinear-decay-estimate1} or \eqref{eqt:gCLM-nonlinear-decay-estimate2} at $t_*$. Without loss of generality we assume $t_*(h)$ is the first time that equality holds for \eqref{eqt:gCLM-nonlinear-decay-estimate1} or \eqref{eqt:gCLM-nonlinear-decay-estimate2}. The choice of parameters guarantees that $f_h$ lies in the bootstrap regime given in Lemma~\ref{lem:bootstrap} for any $h\in \mathbb{B}_{\lambda}(\rho\delta)$. Hence equality must hold for \eqref{eqt:gCLM-nonlinear-decay-estimate2} at $t_*(h)$. 
Define 
\[q(h,t):=\econst^{4\eta t}\left\|\mtx{\Pi}_{\lambda}(f_h(t))\right\|_*^2.\] 
On the boundary of the bootstrap regime $q(h,t)=\rho^2\delta^2$, we compute 
\begin{align*} 
      \frac{\partial q}{\partial t}& =\econst^{4\eta t}\left(4\eta\left\|\mtx{\Pi}_{\lambda}(f_h(t))\right\|_*^2+2 \langle \mtx{L}(\mtx{\Pi}_{\lambda} f_h)+\mtx{\Pi}_{\lambda} \mtx{N}(f_h),\mtx{\Pi}_{\lambda}(f_h) \rangle_*\right)\\
      &\geq \econst^{4\eta t} \left((4\eta-2\lambda)\left\|\mtx{\Pi}_{\lambda}(f_h(t))\right\|_*^2 -C\left\|\mtx{\Pi}_{\lambda}(f_h(t))\right\|_*\left\|f_h(t)\right\|_{\cX_1}^2\right)\\
      &\geq(4\eta-2\lambda)\rho^2\delta^2 -C \rho\delta^3\\
      & \geq (\eta-C\delta/\rho)\rho^2\delta^2.
\end{align*}
We therefore choose $\delta_0$ such that 
\[\delta_0<\min\left\{\delta_1,\frac{\eta\rho}{2C}\right\}.\]

We then show that for such $\nu, \rho,\delta_0$, for any $\xi\in \cX_1$ with $\|\xi\|_{\cX_1}\leq \nu \delta$, there exists $h\in \mathbb{B}_{\lambda}(\rho\delta)$ such that \eqref{eqt:gCLM-nonlinear-decay-estimate1} and \eqref{eqt:gCLM-nonlinear-decay-estimate2} hold. If this is not true, the above argument shows that for any $h\in \mathbb{B}_{\lambda}(\rho\delta)$, there exists a finite time $t_*(h)$ such that $q(h,t_*(h))=\rho^2\delta^2$ and 
\[\left.\frac{\partial q}{\partial t}\right|_{(h,t_*(h))}>0.\]
By continuous dependence of the local flow on its initial data, \(q\) is
continuous in \((h,t)\), and \(q_t\) is continuous near the exit set. The
strict transversality condition above and the standard continuity lemma for
transversal first-exit times imply that \(h\mapsto t_*(h)\) is continuous.
We may therefore define the continuous retraction
\[R(h):=\econst^{2\eta t_*(h)}\mtx{\Pi}_{\lambda}(f_h(t_*(h)))\]
from \(\mathbb{B}_{\lambda}(\rho\delta)\) to its boundary.
For \(h\in\partial\mathbb{B}_{\lambda}(\rho\delta)\), we have $t_*(h)=0$ and thus \(R(h)=h\). By Brouwer's no-retraction theorem, we obtain a contradiction. The proof is thus completed.
\end{proof}

As a direct corollary of the above shooting argument, the next theorem, which is the main result of this section, establishes asymptotically self-similar blowups of the gCLM equation \eqref{eqt:gCLM-model} from $C_c^{\infty}$ initial data for the non-expanding case.
\begin{theorem}\label{thm:gCLM-nonexpanding-blowup}
      For any $a$ such that \eqref{eqt:gCLM-rescaled} admits a steady-state solution \((\bar\Om,\bar U,\bar c_l,-1)\) with $\bar{c}_l\geq 0$ given by Theorem \ref{thm:existence_of_self_similar_profile_gclm}, 
      there exists an odd initial datum $\Omega_0\in C_c^{\infty}(\R)$ such that the solution $\Omega(\cdot,t)$ to \eqref{eqt:gCLM-rescaled} with the normalization \eqref{eqt:gCLM-normalization} and initial data $\Omega_0$ satisfies 
      \[\left\|\Omega(\cdot,t)-\bar{\Omega} \right\|_{\cX_1}\lesssim \econst^{-\eta t}.\]
      As a result, starting from the same initial data $\Omega_0$, the solution $\omega$ of the original gCLM equation \eqref{eqt:gCLM-model} develops an asymptotically self-similar singularity in finite time.
\end{theorem}
\begin{proof}
Let $\delta_0,\nu,\rho$ be given by Theorem~\ref{thm:brouwer}. We can choose a sufficiently small perturbation $\xi\in\cX_1$ such that
\begin{equation}\label{eqt:gCLM-smooth-tail}
 \|\xi\|_{\cX_1}\leq\frac{\nu\delta_0}{2},\quad
 \bar\Omega+\xi\in C_c^\infty(\R).
\end{equation}
Indeed, take an even cutoff function $\chi\in C_c^\infty(\R)$ with $0\leq\chi\leq1$, $\chi=1$ on $[-1,1]$, and $\chi=0$ outside $[-2,2]$. Set $\xi_R:=(\chi(x/R)-1)\bar\Omega$. Then $\bar\Omega+\xi_R=\chi(x/R)\bar\Omega$ is smooth and supported in $[-2R,2R]$. Direct calculation shows that
\[ \|\xi_R\|_{\cX_1}^2 \lesssim\int_{|x|\geq R}\left(\left|\bar\Omega(x)\right|^2+x^2\left|\bar\Omega_x(x)\right|^2\right)\phi(x)\idiff x\longrightarrow0\quad\text{as }R\to+\infty.\]
Thus one can choose $R$ sufficiently large so that \eqref{eqt:gCLM-smooth-tail} holds. We therefore apply Theorem~\ref{thm:brouwer} with $\delta=\delta_0/2$ to obtain $h\in E_\lambda$ and $f_h$ with the stated decay. Define
\[\Omega_0:=\bar\Omega+\xi+\mtx{J}_\lambda\bigl(h-\mtx{\Pi}_\lambda(\xi)\bigr).\]
 Then we have $\Omega_0\in C_c^\infty(\R)$. Moreover, by Theorem~\ref{thm:brouwer}, the solution $\Omega(\cdot,t)$ to \eqref{eqt:gCLM-rescaled} with initial data $\Omega_0$ satisfies 
 \[\left\|\Omega(\cdot,t)-\bar{\Omega} \right\|_{\cX_1}=\|f_h(\cdot,t)\|_{\cX_1} \leq \delta\econst^{-\eta t},\]
 which is the desired decay estimate.

By \eqref{eqt:gCLM-normalization} and Lemma~\ref{lem:coefficient-bounds}, we have 
\[ \left\|\Omega(\cdot,t)-\bar\Omega\right\|_{L^\infty(\R)} +\left|c_l(t)-\bar c_l\right|+|c_\om(t)+1|
 \lesssim\delta\econst^{-\eta t}.\] 
This implies the asymptotically self-similar blowup of the physical variable $\om$ by \cite[Proposition~2.1]{huang2026singular}. The proof is thus completed.
\end{proof}

\section{Finite-codimensional stability of gCLM: expanding case}
\label{sec:compact}

In this section, we establish finite-codimensional stability of the dynamic rescaling equation \eqref{eqt:gCLM-rescaled} for each $a$ such that \eqref{eqt:gCLM-rescaled} admits a steady-state solution \((\bar\Om,\bar U,\bar c_l,-1)\) of \hyperref[item:profile-type-3]{type~(\ref*{item:profile-type-3})} given by Theorem \ref{thm:existence_of_self_similar_profile_gclm}. That is the case where $\bar c_l<0$ and the self-similar profile $\bar\Om$ is compactly supported on $[-L,L]$. In this case, the rescaled flow field 
\[ \bar{V}(x):=a\bar U(x)+\bar c_lx\]
changes sign at the support boundary $x=L$. Intuitively, this drives perturbations outside the support towards the boundary, where compression can amplify their spatial derivatives and complicate the weighted Sobolev estimates. We overcome this difficulty by choosing singular weights with different behavior inside and outside the support to capture the out-going property of $\bar V$. Combined with the boundary damping in \eqref{eqt:gCLM-npc-BL}, these weights control the perturbations transported there and allow us to establish finite-codimensional stability.

Denote the perturbation around the steady state by $f:=\Om-\bar\Om$ and define
\[ \mtx{H}_0(f):=\mtx{H}(f)-\mtx{H}(f)(0),\quad
 \mtx{V}(f)(x):=\int_0^x\mtx{H}_0(f)(y)\idiff y, \]
and 
\begin{equation}
 \mtx{m}(f):=\mtx{H}(f)(0)+\frac{\mtx{V}(f)(L)}L =\frac1L\int_0^L\mtx{H}(f)(y)\idiff y.
\end{equation}
We impose the following normalization conditions in the dynamic rescaling equation \eqref{eqt:gCLM-rescaled}:
\begin{equation}\label{eqt:gCLM-npc-normalization}
 c_l=\bar c_l-a\mtx{m}(f),\quad
 c_\om=-1-a\mtx{m}(f)-(1-a)\mtx{H}(f)(0).
\end{equation}
It is not hard to see that under the above normalization conditions,
\[\Omega_{x,t}(0,t)\equiv 0.\]
Furthermore, if we assume the initial perturbation $f(x,0)$ is supported in $[-L,L]$, then $f(\pm L,t)\equiv 0$, and therefore $f$ is supported in $[-L,L]$ for all time. However, we will not assume that the perturbation is supported in $[-L,L]$ in the rest of this section.

In the rest of this section, we study the perturbation equation for $f$:
\begin{equation}
 f_t=\mtx{L}(f)+\mtx{N}(f),
 \label{eqt:gCLM-npc-perturbation}
\end{equation}
where
\begin{align}
 \mtx{L}(f)&=-\bar{V}f_x+(\bar{U}_x-1)f
 +\bar\Om\left(\mtx{H}_0(f)-a\frac{\mtx{V}(f)(L)}L\right)-a\bar\Om_x
 \left(\mtx{V}(f)-\frac{x}{L}\mtx{V}(f)(L)\right),\label{eqt:gCLM-npc-perturbation_linear}\\
 \mtx{N}(f)&=-a\left(\mtx{V}(f)-\frac{x}{L}\mtx{V}(f)(L)\right)f_x+\left(\mtx{H}_0(f)-a\frac{\mtx{V}(f)(L)}L\right)f.\label{eqt:gCLM-npc-perturbation_nonlinear}
\end{align}

\subsection{Choice of singular weights}
Suppose that the initial perturbation satisfies $f_x(0,0)=0$. By the normalization conditions above we know that
\begin{equation}\label{eqt:gCLM-npc-slope0}
 f_x(0,t)=0
\end{equation}
for all $t>0$.
To capture the out-going effect established in Theorem \ref{thm:out-going}, we introduce singular weighted $L^2$ and $H^1$ spaces which are adapted to \eqref{eqt:gCLM-npc-slope0} as follows. Define
\[\phi_0(x):=\frac{\la x\ra^4}{x^4}.\]
We choose the singular $L^2$ weight as $\phi_0$ times a corrector:
\begin{equation}
 \phi(x):=\econst^{\kappa(x)}\phi_0(x),
\end{equation}
and define the singular $H^1$ weight as
\[\psi(x):=x^2\phi(x).\]
Here $\kappa(x)$ is a bounded smooth function to be determined later. One can check that $\phi$ and $\psi$ satisfy the asymptotics
\begin{align*}
 \phi(x)&\asymp |x|^{-4},\quad \text{as }x\to0,\quad \phi(x)\asymp 1,\quad \text{as }x\to+\infty,\\
 \psi(x)&\asymp |x|^{-2},\quad \text{as }x\to0,\quad \psi(x)\asymp x^2,\quad \text{as }x\to+\infty.
\end{align*}
We will choose $\kappa$ to be nonincreasing on $(0,L)$ and nondecreasing on $(L,+\infty)$, while remaining constant near $0$, $L$, and infinity. Its variation elsewhere will be sufficiently fast to capture the out-going property of $\bar V$. Together with the above asymptotics of $\phi$ and the boundary damping, this will give the damping estimates for the local part of the linearized operator $\mtx{L}$. The choice of $\kappa$ is determined in the next lemma.

\begin{lemma}\label{lem:npc-weight}
There are a bounded smooth even function $\kappa$, constant near infinity, and a constant $\eta_0>0$ such that, with
\begin{equation}\label{eqt:gCLM-npc-weight}
      \phi=\econst^{\kappa}\phi_0, \quad \psi=x^2\phi,
\end{equation}
the following inequalities hold 
\begin{align}
 d_0(x)&:=(\bar{U}_x-1)+\frac{\bar{V}_x}{2}+\frac{\bar{V}}{2}\frac{\phi_x}{\phi}\leq-\eta_0,\label{eqt:gCLM-npc-d0}\\
 d_1(x)&:=(\bar{U}_x-1)-\frac{\bar{V}_x}{2}+\frac{\bar{V}}{2}\frac{\psi_x}{\psi}\leq-\eta_0.\label{eqt:gCLM-npc-d1}
\end{align}
\end{lemma}

\begin{proof}
Substituting \eqref{eqt:gCLM-npc-weight} into the formulas for $d_0$ and $d_1$ yields
\begin{align*}
 d_0(x)&=(\bar{U}_x-1)+\frac{\bar{V}_x}{2}+\frac {\bar{V}}{2x}\frac{x(\phi_0)_x}{\phi_0}+\frac{\bar{V}}{2}\kappa_x=:d_0^{(0)}(x)+\frac{\bar{V}}{2}\kappa_x,\\
 d_1(x)&=(\bar{U}_x-1)-\frac{\bar{V}_x}{2}+ \frac{\bar{V}}{x}+ \frac {\bar{V}}{2x}\frac{x(\phi_0)_x}{\phi_0}+\frac{\bar{V}}{2}\kappa_x=:d_1^{(0)}(x)+\frac{\bar{V}}{2}\kappa_x.
\end{align*}
By Theorem~\ref{thm:out-going} and Lemma~\ref{lem:profile-derivative-estimates}, we have
\[d_0^{(0)}(0)=d_1^{(0)}(0)=-\frac{\bar{V}_x(0)}2<0.\]
\[d_0^{(0)}(+\infty)=d_1^{(0)}(+\infty)= -1+\frac{\bar c_l}{2}<0.\]
\[d_0^{(0)}(L)<d_1^{(0)}(L)=(\bar{U}_x(L)-1)-\frac{\bar{V}_x(L)}{2}\leq \frac{(\bar{U}_x(L)-1)}{2}<0.\]
Therefore, there exist $\delta>0$ and $\eta_0$ such that 
\[\max\left\{d_0^{(0)},d_1^{(0)}\right\}\leq-\eta_0\] 
holds on $\mathcal U:= (0,\delta)\cup(L-\delta,L+\delta)\cup(1/\delta,+\infty)$.
We then take $\mathcal U'=(0,\delta/2)\cup(L-\delta/2,L+\delta/2)\cup(2/\delta,+\infty)$ and choose a smooth cut-off function 
$\zeta\in C_c^\infty(0,+\infty)$ such that
\[0\leq\zeta\leq1,\quad \zeta=0\ \hbox{on }\mathcal U',\quad \zeta=1\ \hbox{on }(0,+\infty)\setminus\mathcal U.\]
The compact set $\mathcal{K}=(0,+\infty)\setminus\mathcal U$ contains no zero of $\bar{V}$. Hence, we have 
\[ \bar{V}_*:=\min_{x\in\mathcal{K}}\left|\bar{V}\right|>0,\quad
 M_*:=\max_{x\in \mathcal K,\,j=0,1}\left\{d_j^{(0)}\right\}<\infty. \]
Choose a sufficiently large $M>0$ such that $M_*-M\bar{V}_*/2\leq-\eta_0$. We define
\begin{equation}
 \kappa(x)=-M\int_0^x \operatorname{sgn}(\bar{V}(s))\zeta(s)\idiff s, \quad \text {for } x\geq 0
\end{equation}
and extend the definition to $(-\infty,0)$ evenly. Then one can see that $\kappa$ is a bounded smooth function such that \eqref{eqt:gCLM-npc-d0} and \eqref{eqt:gCLM-npc-d1} hold. The proof is thus completed.
\end{proof}

In the rest of this section we choose the weight functions $\phi$ and $\psi$ according to Lemma \ref{lem:npc-weight}. With these weight functions, we will establish nonlinear stability in a weighted $H^1$-norm of the perturbation. Let us define the weighted $L^2$ space
\[ \cX_0:=\left\{f\in L^2_{\rm loc}(\R):f\text{ is odd},\ \|f\|_{\cX_0}<\infty\right\},\]
with
\[ \left\|f\right\|_{\cX_0}^2:=\int_\R f(x)^2\phi(x)\idiff x, \]
and the weighted $H^1$ space
\[ \cX_1:=\left\{f\in\cX_0\cap H^1_{\rm loc}(\R):\|f\|_{\cX_1}<\infty\right\},\]
with
\[ \left\|f\right\|_{\dot{\cX}_1}^2:=\int_\R f_x(x)^2\psi(x)\idiff x,\quad
 \left\|f\right\|_{\cX_1}^2=\left\|f\right\|_{\dot{\cX}_1}^2+\left\|f\right\|_{\cX_0}^2.\]
For the stability analysis below, we first record some coefficient estimates in these weighted spaces. The proof of the following lemma is the same as that of Lemma~\ref{lem:coefficient-bounds} with $\alpha=0$. We omit the details for brevity.
\begin{lemma}
\label{lem:npc-coeff}
There is a constant $C$ such that 
\begin{align}
 |\mtx{H}(f)(0)|+ \left\|\frac{\mtx{V}(f)}x\right\|_{L^\infty(\R)}+\left\|\frac{\mtx{V}(f)}{x^2}\right\|_{L^\infty(\R)}
 &\leq C\left\|f\right\|_{\cX_0},\label{eqt:gCLM-npc-low-coeff}\\
 \left\|\frac{\mtx{H}_0(f)}{x^2}\right\|_{L^2(\R)}
 +\left\|\frac{\mtx{V}(f)}{x^3}\right\|_{L^2(\R)}
 &\leq C\left\|f\right\|_{\cX_0},\label{eqt:gCLM-npc-low-L2}\\
 \left\|f\right\|_{L^\infty(\R)}+\left\|\mtx{H}(f)\right\|_{L^\infty(\R)}
 &\leq C\left\|f\right\|_{\cX_1}.
\end{align}
\end{lemma}

\subsection{Linear finite-codimensional stability}
In this subsection we establish finite-codimensional stability of the linearized equation
\begin{equation}\label{eqt:gCLM-npc-linearized}
 f_t=\mtx{L}(f),
\end{equation}
where $\mtx{L}$ is defined in \eqref{eqt:gCLM-npc-perturbation_linear}. Similar to the previous section, we will decompose the linearized operator $\mtx{L}$ into a damping part $\mtx{L}_{\rm stb}$ and a compact part $\mtx{K}$, and then apply the semigroup theory to establish finite-codimensional stability. Denote
\[\mtx{P}(f):=\bar\Om \mtx{H}_0(f),\quad \mtx{K}_{\rm pr}(f):=-a\bar\Om_x\mtx{V}(f),\quad \mtx{K}_e(f):=-a\frac{\mtx{V}(f)(L)}L\left(\bar\Om-x\bar\Om_x\right). \]
Similar to the non-expanding case, $\mtx{P}(f)$ is a nonlocal and noncompact term. To handle this term, we further decompose it into a symmetric part and a skew-symmetric part with respect to $\cX_0$. Denote
\[ \mtx{P}_{\rm sym}:=\frac{\mtx{P}+\mtx{P}^*}{2},\quad \mtx{P}_{\rm skew}:=\frac{\mtx{P}-\mtx{P}^*}{2},\]
where $\mtx{P}^*$ is the adjoint of $\mtx{P}$ with respect to $\cX_0$. We define
 \begin{align*}
 \mtx{L}_{\rm stb}(f)&:=-\bar{V}f_x+(\bar{U}_x-1)f+\mtx{P}_{\rm skew}(f),\\
 \mtx{K}(f)&:=\mtx{P}_{\rm sym}(f)+\mtx{K}_{\rm pr}(f)+\mtx{K}_e(f).
 \end{align*}
Then $\mtx{L}=\mtx{L}_{\rm stb}+\mtx{K}$. The next lemma establishes the compactness of the nonlocal term $\mtx{K}$.

\begin{lemma}
\label{lem:npc-compact}
The operator $\mtx{P}$ is bounded on $\cX_0$, and $\mtx{P}_{\rm sym}$ is Hilbert--Schmidt. Moreover, $\mtx{K}_{\rm pr}$ and $\mtx{K}_e$ are compact on $\cX_0$.
\end{lemma}

\begin{proof}
The operators $\mtx{P}$ and $\mtx{P}_{\rm sym}$ have the same form as in the non-expanding case. The boundedness of $\mtx{P}$ and the Hilbert--Schmidt property of $\mtx{P}_{\rm sym}$ therefore follow by arguments similar to those in Lemma~\ref{lem:gCLM_symmetric_compact}. In fact, these arguments are simpler here because $\bar\Om$ is compactly supported. We omit the details for brevity.

By \eqref{eqt:gCLM-npc-low-coeff} and \eqref{eqt:gCLM-npc-low-L2},
\[\left\|\frac{\mtx{V}(f)}{x^2}\right\|_{L^2(-L,L)} \lesssim \left\|\frac{\mtx{V}(f)}{x^3}\right\|_{L^2(\R)}\lesssim\left\|f\right\|_{\cX_0},\]
and
\[ \left\|\left(\frac{\mtx{V}(f)}{x^2}\right)_x\right\|_{L^2(-L,L)}
 =\left\|\frac{\mtx{H}_0(f)}{x^2}-2\frac{\mtx{V}(f)}{x^3}\right\|_{L^2(-L,L)}\lesssim \left\|\frac{\mtx{H}_0(f)}{x^2}\right\|_{L^2(\R)}+ \left\|\frac{\mtx{V}(f)}{x^3}\right\|_{L^2(\R)} \lesssim \left\|f\right\|_{\cX_0}. \]
Consequently, we obtain 
\[
 \left\|\frac{\mtx{V}(f)}{x^2}\right\|_{H^1(-L,L)}
 \lesssim\left\|f\right\|_{\cX_0}.
\]
By the Rellich theorem, the map 
\[\mtx{K}_0(f):=\frac{\mtx{V}(f)}{x^2}\]
is compact from $\cX_0$ to $L^2_{\rm odd}(-L,L)$. Since $x^2\phi^{1/2}\bar\Om_x\in L^\infty(\R)$ and $\bar\Om_x$ is supported on $[-L,L]$, the map
\[\mtx{B}_0(g):=-ax^2\bar\Om_x g\]
is bounded from $L^2_{\rm odd}(-L,L)$ to $\cX_0$, where the product is extended by zero outside $(-L,L)$. Thus $\mtx{K}_{\rm pr}=\mtx{B}_0\mtx{K}_0$ is compact on $\cX_0$.

Note that $\bar\Om-x\bar\Om_x\in\cX_0$ and  
\[f \mapsto\frac{\mtx{V}(f)(L)}{L} \]
is a bounded linear functional on $\cX_0$ by \eqref{eqt:gCLM-npc-low-coeff}. Therefore $\mtx{K}_e$ is finite rank and thus compact. This completes the proof.
\end{proof}

The next lemma shows that the damping part $\mtx{L}_{\rm stb}$ generates a strongly continuous semigroup on $\cX_0$ with exponential decay.
\begin{lemma}
\label{lem:npc-Lstb-semigroup}
The operator $\mtx{L}_{\rm stb}$, with domain
\begin{equation}\label{eqt:gCLM-npc-domain}
 \mathcal D(\mtx{L}_{\rm stb})=
 \{f\in\cX_0:\bar{V}f_x\in\cX_0\},
\end{equation}
generates a strongly continuous semigroup on $\cX_0$ satisfying
\begin{equation}
 \left\|\econst^{t\mtx{L}_{\rm stb}}\right\|_{\cX_0\to\cX_0}
 \leq\econst^{-\eta_0t},
 \label{eqt:gCLM-npc-Lstb-decay}
\end{equation}
where $\eta_0$ is given by Lemma \ref{lem:npc-weight}.
\end{lemma}

\begin{proof}
We first study the generator of the pure transport semigroup. Set
\[ \mtx{L}_{\rm tr}(f):=-\bar{V}f_x. \]
 Identify odd functions with their restrictions to $(0,\infty)$. We set $I_1=(0,L)$ and $I_2=(L,\infty)$ and introduce the characteristic coordinates
\[ r_1(x):=\int_{L/2}^x\frac{\idiff y}{\bar{V}(y)},\quad 0<x<L, \quad r_2(x):=\int_{2L}^x\frac{\idiff y}{\bar{V}(y)},\quad x>L. \]
Thus each $r_j$ maps $I_j$ onto $\R$ and $(r_j)_x=\bar{V}^{-1}$. Write the inverse of each map as $x=x(r_j)$ and define
\[ \mtx{J}_j(f)(r_j):=\sqrt{2\left|\bar{V}(x(r_j))\right|\phi(x(r_j))}\,f(x(r_j)), \quad \mtx{J}:=\mtx{J}_1\oplus \mtx{J}_2.\]
Then $\mtx{J}$ is an isometry from $\cX_0$ onto $L^2(\R)\oplus L^2(\R)$, and a direct calculation gives
\[ \mtx{J}\mtx{L}_{\rm tr}\mtx{J}^{-1}=\left(-\partial_{r_1}+g_1(r_1),-\partial_{r_2}+g_2(r_2) \right),\quad
 g_j(r_j):=\frac12\left(\bar{V}_x+\bar{V}\frac{\phi_x}{\phi}\right)(x(r_j)),\quad j=1,2.\]
Note that $(-\partial_{r_1},-\partial_{r_2})$ with domain $H^1(\R)\oplus H^1(\R)$ generates the standard translation semigroup on $L^2(\R)\oplus L^2(\R)$. Since $g_j\in L^\infty(\R)$, bounded perturbation theorem yields $\mtx{J}\mtx{L}_{\rm tr}\mtx{J}^{-1}$ with the same domain generates a strongly continuous semigroup on $L^2(\R)\oplus L^2(\R)$. Therefore, $\mtx{L}_{\rm tr}$ generates a strongly continuous semigroup on $\cX_0$. The domain of $\mtx{L}_{\rm tr}$ is given by 
\[\mathcal D(\mtx{L}_{\rm tr})= \mtx{J}^{-1}\bigl(H^1(\R)\oplus H^1(\R)\bigr)=\left\{f\in\cX_0:\phi^{1/2}\bar{V}f_x\in L^2(0,L), \phi^{1/2}\bar{V}f_x\in L^2(L,\infty)\right\}.\]
Since $\bar{V}(L)=0$, the domain above is exactly \eqref{eqt:gCLM-npc-domain}.

The rest of the proof is identical to that of Lemma \ref{lem:Lstb-semigroup}. Since $(\bar{U}_x-1)+\mtx{P}_{\rm skew}$ is bounded on $\cX_0$, the bounded perturbation theorem gives that $\mtx{L}_{\rm stb}$ generates a strongly continuous semigroup with $\mathcal D(\mtx{L}_{\rm stb})=\mathcal D(\mtx{L}_{\rm tr})$. The energy argument
in Lemma~\ref{lem:Lstb-semigroup} applies verbatim and gives
\[\langle \mtx{L}_{\rm stb}(f),f\rangle_{\cX_0}=\int_\R d_0|f|^2\phi\idiff x\leq-\eta_0\left\|f\right\|_{\cX_0}^2,\]
and hence \eqref{eqt:gCLM-npc-Lstb-decay}. This completes the proof.
\end{proof}

Notice that $\mtx{L}$ is a bounded perturbation of $\mtx{L}_{\rm stb}$ by the compact operator $\mtx{K}$. By the bounded perturbation theorem, $\mtx{L}$ with domain 
\[\mathcal D(\mtx{L})=\mathcal D(\mtx{L}_{\rm stb})\]
 also generates a strongly continuous semigroup $\econst^{t\mtx{L}}$ on $\cX_0$. As in Subsection~\ref{sec:linear}, we fix $\lambda\in(0,\eta_0)$ and define the spectral projections $\mtx{\Pi}_\lambda$, $\mtx{\Pi}_s$ and the corresponding spaces $E_\lambda$, $E_s$ for $\mtx{L}$. In particular, $E_\lambda$ is finite dimensional, and
\[ \cX_0=E_\lambda\oplus E_s.\]
An argument similar to that in the proof of Lemma~\ref{lem:linear-package} gives the following decay estimate of the semigroup $\econst^{t\mtx{L}}$ restricted to $E_s$, which establishes finite-codimensional stability of \eqref{eqt:gCLM-npc-linearized}.
\begin{lemma}
\label{lem:npc-linear}
For every $\eta_s\in(0,\lambda)$, there is a constant $C_{\eta_s}>0$ such that
\begin{equation}
 \left\|\econst^{t\mtx{L}}\mtx{\Pi}_s\right\|_{\cX_0\to\cX_0}
 \leq C_{\eta_s}\econst^{-\eta_st},\quad t\geq0.
\end{equation}
\end{lemma}

\subsection{Nonlinear finite-codimensional stability via shooting}
In this subsection, we prove nonlinear finite-codimensional stability of \eqref{eqt:gCLM-npc-perturbation} via a shooting method and establish the asymptotically self-similar blowups of the gCLM equation \eqref{eqt:gCLM-model} from $C_c^{\infty}$ initial data for the expanding case $\bar{c}_l< 0$. Similar to the non-expanding case, let us first perform nonlinear energy estimates as follows.
\begin{lemma}
\label{lem:npc-high}
For the nonlinear operator $\mtx{N}$ defined in \eqref{eqt:gCLM-npc-perturbation_nonlinear}, there is a constant $C$ such that
\begin{equation} \label{eqt:gCLM-npc-nonlinear}
 \left\|\mtx{N}(f)\right\|_{\cX_0}\leq C\left\|f\right\|_{\cX_0}\left\|f\right\|_{\cX_1}, \quad \la \mtx{N}(f),f\ra_{\cX_1}\leq C\left\|f\right\|_{\cX_1}^3.
\end{equation}
Moreover, for every $0<\eta_1<\eta_0$ there is a constant $C=C_{\eta_1}$ such that
\begin{equation}\label{eqt:gCLM-npc-high-energy}
 \frac12\frac{\idiff{}}{\idiff t}\left\|f\right\|_{\cX_1}^2+\eta_1\left\|f\right\|_{\cX_1}^2\leq C_{\eta_1}\left\|f\right\|_{\cX_0}^2+C_{\eta_1}\left\|f\right\|_{\cX_1}^3.
\end{equation}
Consequently,
\begin{equation}\label{eqt:gCLM-npc-high-integral}
 \left\|f(t)\right\|_{\cX_1}
 \leq \econst^{-\eta_1t}\left\|f(0)\right\|_{\cX_1}
 +C_{\eta_1}\int_0^t\econst^{-\eta_1(t-s)}
 \left(\left\|f(s)\right\|_{\cX_0}
 +\left\|f(s)\right\|_{\cX_1}^2\right)\idiff s.
\end{equation}
\end{lemma}

\begin{proof}
By \eqref{eqt:gCLM-npc-low-coeff}, 
\[ \frac{|\mtx{V}(f)(L)|}{L}+\left\|\frac{\mtx{V}(f)}x-\frac{\mtx{V}(f)(L)}L\right\|_{L^\infty(\R)} \lesssim\left\|f\right\|_{\cX_0}. \]
We thus obtain \eqref{eqt:gCLM-npc-nonlinear} by applying the proof of Lemma~\ref{lem:nonlinear_estimate} with the corrected coefficients in \eqref{eqt:gCLM-npc-perturbation}.

We next consider the linear estimate. The linear weighted $L^2$ estimate is identical to that in the proof of Lemma~\ref{lem:high-energy}. For the weighted $\dot{H}^1$ estimate, a direct calculation gives
\[ \la \mtx{L}(f),f\ra_{\dot{\cX}_1}=I_1+I_2+I_3,\]
where
\begin{align*}
 I_1&:=\int_\R\left((\bar{U}_x-1-\bar{V}_x)f_x^2-\bar{V}f_{xx}f_x\right)\psi\idiff x,\\
 I_2&:=\int_\R\mtx{H}(f_x)\bar\Om f_x\psi\idiff x=\la \mtx{P}_{\rm sym}(xf_x),xf_x\ra_{\cX_0},\\
 I_3&:=\int_\R\left(\bar{U}_{xx}f+(1-a)\bar\Om_x\mtx{H}_0(f) -a\bar\Om_{xx} \left(\mtx{V}(f)-\frac{x}{L}\mtx{V}(f)(L)\right)\right)f_x\psi\idiff x.
\end{align*}
Using the same estimates as in Lemma~\ref{lem:high-energy}, we obtain
\[I_1\leq-\eta_0\left\|f\right\|_{\dot{\cX}_1}^2,\]
while Lemma~\ref{lem:compactness-estimate}, applied to $f\mapsto\mtx{P}_{\rm sym}(xf_x)$, gives
 \[ |I_2|\leq\varepsilon\left\|f\right\|_{{\cX}_1}^2+C_\varepsilon\left\|f\right\|_{\cX_0}^2\]
for any $\varepsilon>0$. In fact, the only part of the proof that differs when $\bar c_l<0$ is the estimate of $I_3$, since $\bar{U}_{xx}$ and $\bar\Om_{xx}$ may be unbounded in this case. By Cauchy--Schwarz,
\begin{equation}\label{eqt:gCLM-npc-bound-of-I3}
      |I_3|\leq \left(\left\|x\bar{U}_{xx}f\right\|_{\cX_0}+\left\| (1-a)x\bar\Om_x\mtx{H}_0(f)\right\|_{\cX_0}+\left\|ax\bar\Om_{xx} \left(\mtx{V}(f)-\frac{x}{L}\mtx{V}(f)(L)\right)\right\|_{\cX_0}\right)\left\|f\right\|_{\dot{\cX}_1}.
\end{equation}

We have $\bar{U}_{xx}=\mtx{H}(\bar\Om_x)\in L^2(\R)$. By the oddness and interior smoothness of $\bar\Om$, $x\bar{U}_{xx}=O(x^2)$ at the origin, while Lemma~\ref{lem:profile-derivative-estimates} gives $x\bar{U}_{xx}=O(x^{-2})$ at infinity. Therefore, for every sufficiently large $R$,
\begin{equation}\label{eqt:gCLM-npc-uniform-tail}
      \left\|x\bar{U}_{xx}f\phi^{1/2}\right\|_{L^2((0,1/R)\cup(R,+\infty))}\lesssim \frac{1}{R^2}\left\|f\right\|_{\cX_0}.
\end{equation}  
On any finite interval $(1/R,R)$ with $R>1$, since $H^{1}(1/R,R)$ is compactly embedded into $L^{\infty}(1/R,R)$ and $x\bar{U}_{xx}\in L^2(1/R,R)$, the map $f \mapsto x\bar{U}_{xx} f$ is compact from $H^1(1/R,R)$ to $L^2(1/R,R)$. 
Therefore, given any sequence $\{f_n\}\subset \cX_0$ such that $\|f_n\|_{\cX_0}+\|x(f_n)_x\|_{\cX_0}$ is uniformly bounded, by a standard diagonal argument, we obtain a subsequence of $x\bar{U}_{xx} f_n$ which converges in $L^2(1/R,R)$ for any $R>1$. \eqref{eqt:gCLM-npc-uniform-tail} then makes this subsequence Cauchy in $\cX_0$. Thus the map $f \mapsto x\bar{U}_{xx} f$ is compact from $\cX_1$ to $\cX_0$. Lemma~\ref{lem:compactness-estimate} then yields
\[\left\|x\bar{U}_{xx}f\right\|_{\cX_0} \leq\varepsilon\left\|xf_x\right\|_{\cX_0}+C_\varepsilon\left\|f\right\|_{\cX_0}\]
for any $\varepsilon>0$.

Note that $(1-a)x^2\phi^{1/2}\bar\Om_x\in L^{\infty}(-L,L)$. Hence the second term in \eqref{eqt:gCLM-npc-bound-of-I3} can be controlled as follows:
\begin{align*}
      \left\|(1-a)x\bar\Om_x\mtx{H}_0(f)\right\|_{\cX_0}&=\left\|(1-a)x\bar\Om_x\mtx{H}_0(f)\phi^{1/2}\right\|_{L^2(-L,L)}\\
      &\leq \left\|(1-a)x^2\phi^{1/2}\bar\Om_x\right\|_{L^{\infty}(-L,L)}\left\|\frac{\mtx{H}_0(f)}{x}\right\|_{L^2(\R)}\\
      & \lesssim \left\| \mtx{H}\left(\frac{f}{x}\right)\right\|_{L^2(\R)}\lesssim \left\|f\right\|_{\cX_0}.
\end{align*}

It remains to estimate the third term in \eqref{eqt:gCLM-npc-bound-of-I3}. Note that 
\[\left.\left(\mtx{V}(f)-\frac{\mtx{V}(f)(L)x}{L}\right)\right|_{x=L}=0. \]
Therefore, for $x\in[L/2,L]$, 
\begin{align*}
      \left|\mtx{V}(f)(x)-\frac{x}{L}\mtx{V}(f)(L)\right|&=\left|\int_x^L\left( \mtx{H}_0(f)-\frac{\mtx{V}(f)(L)}{L}\right)\idiff y \right|\leq (L-x)^{1/2} \left\|\mtx{H}_0(f)-\frac{\mtx{V}(f)(L)}L\right\|_{L^2(L/2,L)}\\
      &\lesssim (L-x)^{1/2} \left(\left\|\mtx{H}(f)\right\|_{L^2(\R)}+|\mtx{H}(f)(0)|+\left|\frac{\mtx{V}(f)(L)}{L}\right|\right)\lesssim (L-x)^{1/2}\left\|f\right\|_{\cX_0}.
\end{align*}
We have used the $L^2$ boundedness of the Hilbert transform and \eqref{eqt:gCLM-npc-low-coeff}. 
Therefore, 
\begin{align*}
      &\left\|x\bar\Om_{xx} \left(\mtx{V}(f)-\frac{x}{L}\mtx{V}(f)(L)\right)\right\|_{\cX_0}\\
      &=\sqrt{2}\left\|\phi^{1/2}x\bar\Om_{xx} \left(\mtx{V}(f)-\frac{x}{L}\mtx{V}(f)(L)\right)\right\|_{L^2(0,L)}\\
      & \lesssim  \left\|\phi^{1/2}x\bar\Om_{xx} \left(\mtx{V}(f)-\frac{x}{L}\mtx{V}(f)(L)\right)\right\|_{L^2(0,L/2)}+ \left\|\phi^{1/2}x\bar\Om_{xx} \left(\mtx{V}(f)-\frac{x}{L}\mtx{V}(f)(L)\right)\right\|_{L^2(L/2,L)}\\
      & \lesssim \left\|\phi^{1/2}x^2\bar\Om_{xx}\right\|_{L^{2}(0,L/2)}\left\|\frac{\mtx{V}(f)}{x}-\frac{\mtx{V}(f)(L)}{L} \right\|_{L^{\infty}(\R)}+\left\|\phi^{1/2}x\bar\Om_{xx}(L-x)^{1/2}\right\|_{L^2(L/2,L)}\left\|f\right\|_{\cX_0}\\
      & \lesssim \left\|\frac{\mtx{V}(f)}{x}\right\|_{L^{\infty}(\R)}+\left\|\bar\Om_{xx}(L-x)^{1/2}\right\|_{L^2(L/2,L)}\left\|f\right\|_{\cX_0}\\
      & \lesssim \left\|f\right\|_{\cX_0}.
\end{align*}
We have used \eqref{eqt:gCLM-npc-low-coeff} and Lemma~\ref{lem:profile-derivative-estimates} to deduce the last inequality. The first equality holds because $\mtx{V}(f)-x\mtx{V}(f)(L)/L$ vanishes at $x=\pm L$, so the possible singular terms at the endpoints make no contribution.
Combining the above estimates with Young's inequality yields
\[|I_3|\leq\varepsilon\left\|f\right\|_{\dot{\cX}_1}^2+C_\varepsilon\left\|f\right\|_{\cX_0}^2\leq \varepsilon\left\|f\right\|_{\cX_1}^2+C_\varepsilon\left\|f\right\|_{\cX_0}^2 .\]
As a result, we obtain 
\[ \la \mtx{L}(f),f \ra_{\dot{\cX}_1}\leq (-\eta_0+2\varepsilon)\left\|f\right\|_{\cX_1}^2+C_{\varepsilon}\left\|f\right\|_{\cX_0}^2.\]
We choose $\varepsilon<(\eta_0-\eta_1)/2$. Then the above estimates, together with the weighted $L^2$ estimate and nonlinear estimates \eqref{eqt:gCLM-npc-nonlinear} imply \eqref{eqt:gCLM-npc-high-energy}.
Finally, \eqref{eqt:gCLM-npc-high-integral} follows from \eqref{eqt:gCLM-npc-high-energy} and Grönwall's inequality. The proof is thus completed.
\end{proof}

As in Section~\ref{sec:noncompact}, let $\cX_c$ denote the space of odd $C_c^\infty(\R)$ functions supported away from the origin. We choose correction directions in this space to preserve the smoothness and compact support of $\Omega_0$. By an argument similar to that in the proof of Lemma~\ref{lem:unstable-adapted-norm}, we construct an inner product $\la\cdot,\cdot\ra_*$ on the present space $E_\lambda$, whose induced norm $\|\cdot\|_*$ is equivalent to $\|\cdot\|_{\cX_0}$ and satisfies
\[\la\mtx{L}(h),h\ra_*\geq-\lambda\|h\|_*^2,\quad h\in E_\lambda.\]
Similarly, adapting the proof of Lemma~\ref{lem:correction-directions}, we construct a real-linear correction map $\mtx{J}_\lambda:E_\lambda\to\cX_c$ satisfying
\[\mtx{\Pi}_\lambda\mtx{J}_\lambda(h)=h,\quad\|\mtx{J}_\lambda(h)\|_{\cX_1}\leq C\|h\|_*,\quad h\in E_\lambda.\]
Take $\eta=0.9\lambda$. Using arguments similar to those in the proofs of Lemma~\ref{lem:bootstrap} and Theorem~\ref{thm:brouwer}, with Lemma~\ref{lem:npc-linear} replacing Lemma~\ref{lem:linear-package} and Lemma~\ref{lem:npc-high} replacing Lemmas~\ref{lem:nonlinear_estimate} and~\ref{lem:high-energy}, we obtain finite-codimensional stability of \eqref{eqt:gCLM-npc-perturbation} in the following theorem.
\begin{theorem}\label{thm:npc-main}
There exist \(\delta_0,\nu,\rho>0\) such that for any $\delta\in (0,\delta_0)$ and any initial data $\xi\in \cX_1$ with $\|\xi\|_{\cX_1}\leq \nu \delta$, there exists \(h\in E_{\lambda}\) with \(\|h\|_*\leq\rho\delta\), such that starting from the initial data
\[ f_h(0):=\xi+\mtx{J}_{\lambda}(h-\mtx{\Pi}_{\lambda}(\xi))\]
the solution $f_h(t)$ to \eqref{eqt:gCLM-npc-perturbation} satisfies the decay estimates  
\begin{equation}
 \left\|f_h(t)\right\|_{\cX_1}\leq\delta \econst^{-\eta t},\end{equation}
 \begin{equation}
 \left\|\mtx{\Pi}_{\lambda}(f_h(t))\right\|_*\leq\rho\delta \econst^{-2\eta t},\quad t\geq0.
\end{equation}
\end{theorem}

Finally, we use an argument similar to that in the proof of Theorem~\ref{thm:gCLM-nonexpanding-blowup}, with Theorem~\ref{thm:npc-main} replacing Theorem~\ref{thm:brouwer}, to obtain the main result of this section, which establishes asymptotically self-similar finite-time blowups of the gCLM model \eqref{eqt:gCLM-model} from smooth compactly supported initial data in the expanding case.
\begin{theorem}\label{thm:gCLM-expanding-blowup}
For any $a$ such that \eqref{eqt:gCLM-rescaled} admits a steady-state solution $(\bar\Omega,\bar U,\bar c_l,-1)$ with $\bar c_l<0$ given by Theorem~\ref{thm:existence_of_self_similar_profile_gclm}, there exists an odd initial datum $\Omega_0\in C_c^\infty(\R)$ such that the solution $\Omega(\cdot,t)$ to \eqref{eqt:gCLM-rescaled} with the normalization \eqref{eqt:gCLM-npc-normalization} and initial data $\Omega_0$ satisfies
\[\left\|\Omega(\cdot,t)-\bar\Omega\right\|_{\cX_1}\lesssim\econst^{-\eta t}.\]
As a result, starting from the same initial data $\Omega_0$, the solution $\omega$ of the original gCLM equation \eqref{eqt:gCLM-model} develops an asymptotically self-similar singularity in finite time.
\end{theorem}

\section{Dynamic rescaling formulation and exact self-similar blowup of the HL model}
\label{sec:hl-rescaling}

In this section, we introduce the dynamic rescaling formulation of the 1D Hou--Luo (HL) model \eqref{eqt:1Dhouluo}. We then recall the existence and properties of the self-similar blowup profiles constructed in \cite{huang2025exact}, which will be useful in the subsequent stability analysis.

As in the gCLM case, by \cite[Lemma~2.1]{chen2026novel}, solutions $(\omega,\theta)$ of \eqref{eqt:1Dhouluo} with the same initial data but different choices of zeros of $u$ differ only by a time-dependent spatial translation. Therefore, we may impose $u(0,t)=0$ without loss of generality. 

Following \cite{chen2022asymptotically,huang2025exact}, we will work with the variable $\upsilon:=\theta_x$ instead of $\theta$ throughout the discussion. Under this change of variable, \eqref{eqt:1Dhouluo} is equivalent to
\begin{equation}\label{eqt:HL-model}
    \begin{aligned}
    &\om_t+u\om_x=\upsilon,\\
    &\upsilon_t+u\upsilon_x=-u_x\upsilon,\\
    &u_x=\mtx{H}(\om),\quad u(0,t)=0.
    \end{aligned}
\end{equation}

\subsection{Dynamic rescaling formulation}
Let us consider the change of variables
\begin{equation}
\begin{aligned}
 \omega(x,t)&=C_\om(\tau)^{-1}\Omega(C_l(\tau)x,\tau(t)),\quad 
 u(x,t)=C_\om(\tau)^{-1}C_l(\tau)^{-1}U(C_l(\tau)x,\tau(t)),\\
 \upsilon(x,t)&=C_l(\tau)C_\theta(\tau)^{-1}\Upsilon(C_l(\tau)x,\tau(t))
\end{aligned}
\end{equation}
with
\[
\begin{aligned}
 C_\om(\tau)&=\exp\left(\int_0^\tau c_\om(s)\idiff s\right),\quad
 C_\theta(\tau)=\exp\left(\int_0^\tau c_\theta(s)\idiff s\right),\\
 C_l(\tau)& =\exp\left(\int_0^\tau c_l(s)\idiff s\right),\quad 
 t(\tau)=\int_0^\tau C_\om(s)\idiff s,\quad c_\theta=c_l+2c_\om.
\end{aligned}
\]
Then \eqref{eqt:HL-model} is equivalent to the following dynamic rescaling equation:
\begin{equation}\label{eqt:HL-rescaled}
\begin{aligned}
 &\Omega_\tau+(U+c_lX)\Omega_X=c_\om\Omega+\Upsilon,\\
 &\Upsilon_\tau+(U+c_lX)\Upsilon_X=(2c_\om-U_X)\Upsilon,\\
 &U_X=\mtx{H}(\Omega),\quad U(0,\tau)=0.
\end{aligned}
\end{equation}
The steady state equation of \eqref{eqt:HL-rescaled} is given by
\begin{equation}\label{eqt:HL-steady}
\begin{aligned}
 &(U+c_lX)\Omega_X=c_\om\Omega+\Upsilon,\\
 &(U+c_lX)\Upsilon_X=(2c_\om-U_X)\Upsilon,\\
 &U_X=\mtx{H}(\Omega),\quad U(0)=0.
\end{aligned}
\end{equation}
As in the gCLM case, a nontrivial steady state $(\bar\Omega,\bar\Upsilon,\bar c_l,\bar c_\om)$ of \eqref{eqt:HL-steady} with $\bar c_\om<0$ corresponds to an exact self-similar finite-time blowup of the HL model \eqref{eqt:HL-model}, with spatial shrinking rate $\gamma=-\bar c_l/\bar c_\om$.
Note that if $(\Omega(X),\Upsilon(X),c_l,c_\om)$ is a steady state solution of \eqref{eqt:HL-rescaled}, then
\[\bigl(\alpha\Omega(\beta X),\alpha^2\Upsilon(\beta X),\alpha c_l,\alpha c_\om\bigr)\]
is also a steady state solution for $\alpha>0$ and $\beta>0$, which indicates that this equation admits two degrees of freedom. Hence, in order to uniquely identify the profile up to rescaling, we need to impose two normalization conditions to fix these two degrees of freedom.

Moreover, if under some normalization conditions the solution tuple $(\Omega,\Upsilon,c_l,c_\om)$ converges sufficiently fast to a nontrivial steady state $(\bar\Omega,\bar\Upsilon,\bar c_l,\bar c_\om)$ of \eqref{eqt:HL-rescaled} with $\bar c_\om<0$ as $\tau\to\infty$, then the solution would blow up in an asymptotically self-similar manner. We refer the readers to \cite[Lemma~2.2]{chen2026novel} for a rigorous proof of this statement. Therefore, the study of asymptotically self-similar blowup of \eqref{eqt:HL-model} is reduced to the stability analysis of the dynamic rescaling equation \eqref{eqt:HL-rescaled} around some steady state.

\subsection{Existence of the steady state solutions}
Following a strategy similar to that used for the gCLM profiles in \cite{huang2024self}, Huang, Qin, Wang, and Wei \cite{huang2025exact} established the existence of smooth exact self-similar profiles for the HL model by a purely analytic fixed-point construction. We briefly restate this result in the following theorem; see \cite[Theorems~1.1, 4.2, and~4.4]{huang2025exact} for details.

\begin{theorem}[Existence of exact self-similar profiles]\label{thm:hl-existence}
Equation \eqref{eqt:HL-steady} admits a solution $(\bar\Omega,\bar\Upsilon,\bar c_l,\bar c_\om)$ with $\bar c_\om=-1$ and $2<\bar c_l<4.53$. The profiles $\bar\Omega$ and $\bar\Upsilon$ are smooth odd functions on $\R$, strictly positive on $(0,+\infty)$, and satisfy
\[\bar\Omega_X(0)>0,\quad \bar\Upsilon_X(0)>0.\]
Moreover, $\bar\Omega/X$ and $\bar\Upsilon/X$ are decreasing in $X$ and convex in $X^2$ on $[0,+\infty)$, and
\[\lim_{X\to+\infty}X^{1/\bar c_l}\bar\Omega(X)=C_\Omega,\quad\lim_{X\to+\infty}X^{2/\bar c_l}\bar\Upsilon(X)=C_\Upsilon\]
for some $C_\Omega,C_\Upsilon>0$. Both profiles belong to $L^\infty(\R)\cap\dot H^m(\R)$ for every integer $m\geq1$.
\end{theorem}

As in the gCLM case, the out-going property plays an important role in the stability analysis below. This property follows from the fixed-point representation in \cite[Section~2.2 and Corollary~3.6]{huang2025exact}. We state it in the following theorem and defer the proof to Appendix~\ref{app:hl-profiles}.

\begin{theorem}[Strictly out-going property]\label{thm:hl-out-going}
Let $(\bar\Omega,\bar\Upsilon,\bar c_l,\bar c_\om)$ be the self-similar profile given by Theorem~\ref{thm:hl-existence} with $\bar{c}_\om=-1$. Then
\[\bar U(X)+\bar c_lX\geq\left(\frac{\bar c_l}{2}-1\right)X>0,\quad X>0.\]
At the origin, the following identities hold:
\begin{equation}\label{eqt:HL-origin-identities}
 \bar c_l=2\frac{\bar\Upsilon_X(0)}{\bar\Omega_X(0)},\quad
\bar U_X(0)=-1-\frac{\bar c_l}{2},\quad
 \bar U_X(0)+\bar c_l=\frac{\bar c_l}{2}-1>0.
\end{equation}
\end{theorem}

\subsection{Further properties of the steady state solutions}
We also need further properties of the profiles, which we summarize in the following lemma. These properties follow from the steady-state equations \eqref{eqt:HL-steady} and Theorems~\ref{thm:hl-existence} and~\ref{thm:hl-out-going}. We again defer the proof to Appendix~\ref{app:hl-profiles}.

\begin{lemma}\label{lem:hl-profile-derivatives}
Let $(\bar\Omega,\bar\Upsilon,\bar c_l,\bar c_\om)$ be the self-similar profile given by Theorem~\ref{thm:hl-existence} with $\bar{c}_\om=-1$.  Then
\[\bar U_X(X)\to0,\quad \frac{\bar U(X)}{X}\to0\quad\text{as }X\to+\infty,\]
and $\|X\bar U_{XX}\|_{L^\infty(\R)}<\infty$. Moreover, there exists a constant $C>0$ such that
\begin{equation}\label{eqt:HL-profile-derivatives}
\begin{aligned}
 \left|\bar\Omega(X)\right|+X\left|\bar\Omega_X(X)\right|+X^2\left|\bar\Omega_{XX}(X)\right|&\leq CX^{-1/\bar c_l},\\
 \left|\bar\Upsilon(X)\right|+X\left|\bar\Upsilon_X(X)\right|+X^2\left|\bar\Upsilon_{XX}(X)\right|&\leq CX^{-2/\bar c_l},\quad X\geq1.
\end{aligned}
\end{equation}
\end{lemma}

\section{Finite-codimensional stability analysis of the HL model}
\label{sec:hl-stability}
In this section, we establish finite-codimensional stability of the dynamic rescaling equation \eqref{eqt:HL-rescaled} around the smooth self-similar profiles given by Theorem~\ref{thm:hl-existence}. 
As in Section~\ref{sec:noncompact}, we denote the rescaled space-time variable $(X,\tau)$ by $(x,t)$. In these variables, \eqref{eqt:HL-rescaled} becomes
\begin{equation}\label{eqt:HL-rescaled-xt}
\begin{aligned}
 &\Omega_t+(U+c_lx)\Omega_x=c_\om\Omega+\Upsilon,\\
 &\Upsilon_t+(U+c_lx)\Upsilon_x=(2c_\om-U_x)\Upsilon,\\
 &U_x=\mtx{H}(\Omega),\quad U(0,t)=0.
\end{aligned}
\end{equation}
As in the non-expanding gCLM case, the background rescaled flow 
\[\bar{V}(x):=\bar U(x)+\bar c_lx\] 
has the same sign as the spatial variable $x$ on the whole line by Theorem~\ref{thm:hl-out-going}. This out-going property allows us to choose singular weights to obtain local damping estimates. However, the noncompact nonlocal term appears as a cross term coupling the two perturbation components in the energy estimate. We therefore cannot directly apply the gCLM decomposition into symmetric and skew-symmetric parts to extract a compact symmetric term. Instead, we decompose this noncompact nonlocal term into two parts: a compact part, which we set aside in the damping estimate, and a remaining part, which contributes a cross term that can be absorbed into the local damping using the out-going property.

Denote the perturbation as
\[ (f,g):=(\Omega-\bar\Omega,\Upsilon-\bar\Upsilon),\]
and write $F = (f,g)$. Following \cite{chen2022asymptotically}, we impose the normalization conditions
\begin{equation}\label{eqt:HL-normalization}
 c_l=2\frac{\Upsilon_x(0,t)}{\Omega_x(0,t)},\quad c_\om=\frac{c_l}{2}+U_x(0,t).
\end{equation}
Differentiating \eqref{eqt:HL-rescaled-xt} at the origin gives
\begin{align*}
    \partial_t\Omega_x(0,t) &=(c_\om-c_l-U_x(0,t))\Omega_x(0,t)+\Upsilon_x(0,t),\\
    \partial_t\Upsilon_x(0,t) &=(2c_\om-c_l-2U_x(0,t))\Upsilon_x(0,t).
\end{align*}
The normalization conditions \eqref{eqt:HL-normalization} thus guarantee that $\Omega_x(0,t)$ and $\Upsilon_x(0,t)$ remain constant in time. If the initial slopes $\Omega_x(0,0)$ and $\Upsilon_x(0,0)$ agree with those of the steady state, \eqref{eqt:HL-normalization} reduces to
\begin{equation}\label{eqt:HL-normalization-reduced}
 c_l=\bar c_l,\quad c_\om=-1+\mtx{H}(f)(0),\quad f_x(0,t)=g_x(0,t)=0.
\end{equation}

Let
\[ \mtx{H}_0(f):=\mtx{H}(f)-\mtx{H}(f)(0),\quad \mtx{V}(f)(x):=\int_0^x\mtx{H}_0(f)(y)\idiff y.\]
Then the equations for the perturbation can be written as
\begin{equation}\label{eqt:HL-perturbation}
 F_t=\mtx{L}(F)+\mtx{N}(F),
\end{equation}
where
\begin{align}
 \mtx{L}(f,g)&:=\begin{pmatrix}
 -\bar{V} f_x-f+g-\bar\Omega_x\mtx{V}(f)+\mtx{H}(f)(0)(\bar\Omega-x\bar\Omega_x)\\[5pt]
 -\bar{V} g_x-(2+\bar U_x)g-\bar\Upsilon\mtx{H}_0(f)
 -\bar\Upsilon_x\mtx{V}(f)+\mtx{H}(f)(0)(\bar\Upsilon-x\bar\Upsilon_x)
 \end{pmatrix}^T,\label{eqt:HL-linear-operator}\\
 \mtx{N}(f,g)&:=\begin{pmatrix}
 -\left(\mtx{V}(f)+x\mtx{H}(f)(0)\right)f_x+\mtx{H}(f)(0)f\\[5pt]
 -\left(\mtx{V}(f)+x\mtx{H}(f)(0)\right)g_x+(2\mtx{H}(f)(0)-\mtx{H}(f))g
 \end{pmatrix}^T.\label{eqt:HL-nonlinear-operator}
\end{align}
Here the superscript $T$ denotes transpose, so both $\mtx{L}(f,g)$ and $\mtx{N}(f,g)$ are row vectors, which is consistent with $F=(f,g)$. We will study \eqref{eqt:HL-perturbation} in the rest of this section.

\subsection{Choice of singular weights}
Suppose that the initial perturbation satisfies $f_x(0,0)=g_x(0,0)=0$. By \eqref{eqt:HL-normalization} we know that
\begin{equation}\label{eqt:HL-slope0}
 f_x(0,t)=g_x(0,t)=0
\end{equation}
for all $t>0$. To capture the out-going effect of the rescaled flow field $\bar{V}$, we introduce singular weighted $L^2$ and $H^1$ spaces which are adapted to \eqref{eqt:HL-slope0} as follows. Take
\[\alpha\in\left(-1,2/\bar c_l-1\right)\]
and define
\[\phi_0(x):=\frac{\la x\ra^{\alpha+4}}{x^4}.\]
We choose the singular $L^2$ weight as $\phi_0$ times a corrector:
\begin{equation}
 \phi(x):=\econst^{\kappa(x)}\phi_0(x),
\end{equation}
and define the singular $H^1$ weight as
\[\psi(x):=x^2\phi(x).\]
Here $\kappa(x)$ is a bounded smooth function to be determined later. One can check that $\phi$ and $\psi$ satisfy the asymptotics
\begin{align*}
 \phi(x)&\asymp |x|^{-4},\quad \text{as }x\to0,\quad \phi(x)\asymp x^\alpha,\quad \,\,\,\,\text{as } x\to+\infty,\\
 \psi(x)&\asymp |x|^{-2},\quad \text{as }x\to0,\quad \psi(x)\asymp x^{\alpha+2},\quad \text{as }x\to+\infty.
\end{align*}
The choice of $\kappa$ is determined in the next lemma. Compared with the gCLM case, the energy estimates for the HL model involve additional cross terms. To control these terms, we introduce additional parameters in the damping estimate below. We will explain their roles in the next subsection.
\begin{lemma}\label{lem:hl-weight}
There are a bounded smooth even function $\kappa$, constant near the origin and infinity, and constants $\varepsilon_0>0$, $M\geq1$, and $\eta_0>0$ such that, with
\begin{equation}\label{eqt:HL-weight}
 \phi=\econst^{\kappa}\phi_0,\quad \psi=x^2\phi,
\end{equation}
the following inequalities hold
\begin{equation}\label{eqt:HL-weight-damping}
 \max\{d_{f,j}(x),d_{g,j}(x)\}+\varepsilon_0+\frac{1}{2\sqrt M}+\frac{M}{4\varepsilon_0}\bar\Upsilon(x)^2\leq-\eta_0,
 \quad x\neq0,\quad j=0,1,
\end{equation}
where, for $j=0,1$,
\begin{align*}
 d_{f,j}(x)&:=-1+\frac{\bar{V}_x}{2}+\frac{\bar{V}}{2}\frac{\phi_x}{\phi}+j\left(\frac{\bar{V}}{x}-\bar{V}_x\right),\\
 d_{g,j}(x)&:=-2-\bar U_x+\frac{\bar{V}_x}{2}+\frac{\bar{V}}{2}\frac{\phi_x}{\phi}+j\left(\frac{\bar{V}}{x}-\bar{V}_x\right).
\end{align*}
\end{lemma}
\begin{proof}
Substituting \eqref{eqt:HL-weight} into the formulas for $d_{f,j}$ and $d_{g,j}$ yields
\begin{align*}
 d_{f,j}(x)&=-1+\frac{\bar{V}_x}{2}+\frac{\bar{V}}{2x}\frac{x(\phi_0)_x}{\phi_0}+j\left(\frac{\bar{V}}x-\bar{V}_x\right)+\frac{\bar{V}}2\kappa_x\mathrel{=:}d_{f,j}^{(0)}(x)+\frac{\bar{V}}2\kappa_x,\\
 d_{g,j}(x)&=-2-\bar U_x+\frac{\bar{V}_x}{2}+\frac{\bar{V}}{2x}\frac{x(\phi_0)_x}{\phi_0}+j\left(\frac{\bar{V}}x-\bar{V}_x\right)+\frac{\bar{V}}2\kappa_x\mathrel{=:}d_{g,j}^{(0)}(x)+\frac{\bar{V}}2\kappa_x.
\end{align*}
By Theorem~\ref{thm:hl-out-going} and Lemma~\ref{lem:hl-profile-derivatives}, we have, for $j=0,1$,
\begin{align*}
 d_{f,j}^{(0)}(0)&=-\frac{3\bar c_l-2}{4}<0,\quad d_{g,j}^{(0)}(0)=-\frac{\bar c_l-2}{4}<0,\\
 d_{f,j}^{(0)}(+\infty)&=-1+\frac{(1+\alpha)\bar c_l}{2}<0,\quad d_{g,j}^{(0)}(+\infty)=-2+\frac{(1+\alpha)\bar c_l}{2}<0.
\end{align*}
Choose $\varepsilon_0$ small and then $M$ large so that these limits remain strictly negative after adding $\varepsilon_0+1/(2\sqrt M)$. Since $\bar\Upsilon$ vanishes at both ends, there exist $\delta>0$ and $\eta_0>0$ such that
\[\max_{j=0,1}\left\{d_{f,j}^{(0)},d_{g,j}^{(0)}\right\} +\varepsilon_0+\frac1{2\sqrt M}+\frac{M}{4\varepsilon_0}\bar\Upsilon^2 \leq-\eta_0\]
holds on $\mathcal U:=(0,\delta)\cup(1/\delta,+\infty)$. We then take $\mathcal U'=(0,\delta/2)\cup(2/\delta,+\infty)$ and choose a smooth cut-off function $\zeta\in C_c^\infty(0,+\infty)$ such that
\[0\leq\zeta\leq1,\quad \zeta=0\ \hbox{on }\mathcal U',\quad \zeta=1\ \hbox{on }(0,+\infty)\setminus\mathcal U.\]
The compact set $\mathcal K=(0,+\infty)\setminus\mathcal U$ contains no zero of $\bar{V}$. Hence, we have $\bar{V}_*:=\min_{x\in\mathcal K}\bar{V}>0$ and
\[ M_*:=\max_{x\in\mathcal K}\left\{\max_{j=0,1}\left\{d_{f,j}^{(0)}(x),d_{g,j}^{(0)}(x)\right\}+\varepsilon_0+\frac1{2\sqrt M}+\frac{M}{4\varepsilon_0}\bar\Upsilon(x)^2\right\}<\infty.\]
Choose a sufficiently large $A>0$ such that $M_*-A\bar{V}_*/2\leq-\eta_0$. We define $\kappa(x)=-A\int_0^x\zeta(s)\idiff s$ for $x\geq0$ and extend the definition to $(-\infty,0)$ evenly. Then one can see that $\kappa$ is a bounded smooth function such that \eqref{eqt:HL-weight-damping} holds. The proof is thus completed.
\end{proof}

In the rest of this section, we will fix our choice of the weight functions $\phi$ and $\psi$ as in \eqref{eqt:HL-weight} with the constants $\varepsilon_0$, $\eta_0$, and $M$ determined as in Lemma \ref{lem:hl-weight}. With these weight functions, we will establish nonlinear stability in a weighted $H^1$-norm of the perturbation. Let us define the weighted $L^2$ space
\[ \cX_0:=\left\{F=(f,g)\in L^2_{\rm loc}(\R)^2:f,g\text{ are odd},\ \|F\|_{\cX_0}<\infty\right\},\]
with
\[ \left\|F\right\|_{\cX_0}^2:=\int_\R\bigl(f(x)^2+Mg(x)^2\bigr)\phi(x)\idiff x, \]
and the weighted $H^1$ space
\[ \cX_1:=\left\{F\in\cX_0\cap H^1_{\rm loc}(\R)^2:\|F\|_{\cX_1}<\infty\right\},\]
with
\[ \left\|F\right\|_{\dot{\cX}_1}^2:=\int_\R\bigl(f_x(x)^2+Mg_x(x)^2\bigr)\psi(x)\idiff x,\quad
 \left\|F\right\|_{\cX_1}^2:=\left\|F\right\|_{\dot{\cX}_1}^2+\left\|F\right\|_{\cX_0}^2.\]
Here the constant $M$, as in Lemma \ref{lem:hl-weight}, is introduced to balance the two components in the weighted energy.

For the stability analysis below, we first record some coefficient estimates in these weighted spaces. The proof of the following lemma is the same as that of Lemma~\ref{lem:coefficient-bounds} applied to each component, and is therefore omitted.
\begin{lemma}\label{lem:hl-coefficients}
For $F=(f,g)\in\cX_1$,
\begin{align}
 |\mtx{H}(f)(0)|+\left\|\frac{\mtx{V}(f)}x\right\|_{L^\infty} +\left\|\frac{\mtx{V}(f)}{x^2}\right\|_{L^\infty}&\leq C\|F\|_{\cX_0},\label{eqt:HL-low-coefficients}\\
 \|f\|_{L^\infty}+\|g\|_{L^\infty}+\|\mtx{H}(f)\|_{L^\infty}+\left\|\frac{\mtx{V}(f)}x\right\|_{L^\infty} &\leq C\|F\|_{\cX_1}.
\end{align}
\end{lemma}

\subsection{Linear stability estimates via semigroups}
In this subsection, we establish finite-codimensional stability of the linearized equation 
\begin{equation}\label{eqt:HL-linearized}
 F_t=\mtx{L}(F),
\end{equation}
where $\mtx{L}$ is defined in \eqref{eqt:HL-linear-operator}. Note that in the linear part $\mtx{L}(f,g)$, terms involving $\mtx{V}(f)$ and $\mtx{H}(f)(0)$ have the same structure as the compact terms in Lemma~\ref{lem:gCLM_symmetric_compact}. Specifically, write
\[ \mtx{K}_{\rm pr}(f,g):=\bigl(-\bar\Omega_x\mtx{V}(f)+\mtx{H}(f)(0)(\bar\Omega-x\bar\Omega_x),\, -\bar\Upsilon_x\mtx{V}(f)+\mtx{H}(f)(0)(\bar\Upsilon-x\bar\Upsilon_x)\bigr).\]
This part will be proved to be compact later. Hence, one seeks to obtain damping estimates for the remaining terms of $\mtx{L}$. To illustrate the main difference from the gCLM case, let us perform an energy estimate as follows:
\begin{align*}
 \la(\mtx{L}-\mtx{K}_{\rm pr})F,F\ra_{\cX_0}&=\int_\R\bigl[-\bar{V}(ff_x+Mgg_x)-f^2-M(2+\bar U_x)g^2+fg-M\bar\Upsilon g\mtx{H}_0(f)\bigr]\phi\idiff x\\
 &=\int_\R\left[\frac{(\bar{V}\phi)_x}{2\phi}\left(f^2+Mg^2\right)-f^2-M(2+\bar U_x)g^2+fg-M\bar\Upsilon g\mtx{H}_0(f)\right]\phi\idiff x\\
 &=\int_\R\bigl(d_{f,0}f^2+Md_{g,0}g^2\bigr)\phi\idiff x+\int_\R fg\phi\idiff x-M\int_\R\bar\Upsilon g\mtx{H}_0(f)\phi\idiff x,
\end{align*}
where we have used integration by parts in the second equality. The first integral gives the local damping term. Compared with the gCLM case, we need to deal with two additional cross terms coupling $f$ and $g$. To this end, we seek to handle these cross terms using the damping from the out-going property. For the first cross term, we choose $M$ sufficiently large to make its contribution small relative to the weighted energy:
\[ \left|\int_\R fg\phi\idiff x\right| \leq\frac{1}{2\sqrt M}\int_\R\left(f^2+Mg^2\right)\phi\idiff x.\]
This yields a large coefficient $M$ in the second cross term, which however is harmless near the origin or in the very far field due to the vanishing of $\bar\Upsilon$ in those regions. Eventually, as in the proof of Lemma~\ref{lem:hl-weight}, a suitable choice of $\phi$ significantly strengthens the local damping in the intermediate region (away from the origin and infinity) and therefore takes care of the largeness of $M$.
However, as the second cross term is nonlocal in $f$, a direct weighted estimate would introduce a constant depending on $\phi$, which in turn may compromise the choice of $\phi$.
In order to obtain an estimate that is independent of $\phi$, we move the factor $\sqrt{\phi}$ in the energy pairing inside the Hilbert transform. By Lemma~\ref{lem:core-identities}, we can write
\begin{equation}\label{eqt:HL-hilbert-conjugation}
 \bar\Upsilon\mtx{H}_0(f)=\frac{\bar\Upsilon}{\sqrt\phi}\mtx{H}\left(\sqrt\phi f\right)+\mtx{R}(f),
\end{equation}
where
\[\mtx{R}(f):=\frac{\bar\Upsilon}{\sqrt\phi}\left[x\sqrt\phi,\mtx{H}\right]\left(\frac{f}{x}\right).\]
As explained in Remark~\ref{rem:gCLM-commutator}, the commutator cancels the diagonal singularity and is more regular. This motivates us to establish the compactness of $\mtx{R}$ and estimate the principal term without introducing a weight-dependent constant:
\[\begin{aligned}
 M\left|\int_\R\frac{\bar\Upsilon\mtx{H}\left(\sqrt\phi f\right)}{\sqrt\phi}g\phi\idiff x\right|
 &\leq M\left\|\mtx{H}\left(\sqrt\phi f\right)\right\|_2\left\|\bar\Upsilon\sqrt\phi\,g\right\|_2
 =M\left\|\sqrt\phi f\right\|_2\left\|\bar\Upsilon\sqrt\phi\,g\right\|_2\\
 &\leq\varepsilon_0\int_\R f^2\phi\idiff x
 +\frac{M^2}{4\varepsilon_0}\int_\R\bar\Upsilon^2g^2\phi\idiff x,
\end{aligned}\]
where the equality follows from the $L^2$ isometry of $\mtx{H}$. As shown in Lemma~\ref{lem:hl-weight}, the out-going property allows us to choose the weight so that the local damping absorbs the $g^2$ term above. 

Motivated by the above discussion, let us decompose the linear operator $\mtx{L}$ as 
\[\mtx{L} = \mtx{L}_{\rm stb} + \mtx{K},\]
where 
\[ \mtx{K}(F):=\mtx{K}_{\rm pr}(F) -(0,\mtx{R}(f)), \quad \mtx{L}_{\rm stb}:=\mtx{L}-\mtx{K}\]
for any $F=(f,g)\in \cX_0$.
Combining the above energy estimate and Lemma~\ref{lem:hl-weight}, we obtain
\begin{equation}\label{eqt:HL-energy-remainder}
 \la\mtx{L}_{\rm stb}F,F\ra_{\cX_0}\leq\int_\R\left[\left(d_{f,0}+\varepsilon_0+\frac{1}{2\sqrt M}\right)f^2+M\left(d_{g,0}+\frac{1}{2\sqrt M}+\frac{M}{4\varepsilon_0}\bar\Upsilon^2\right)g^2\right]\phi\idiff x\leq-\eta_0\|F\|_{\cX_0}^2.
\end{equation}
Thus $\mtx{L}_{\rm stb}$ is a damping term. To establish the linear stability of \eqref{eqt:HL-linearized} via semigroups, we need to show the compactness of the remaining nonlocal term $\mtx{K}$. 
\begin{lemma}\label{lem:hl-compact}
$\mtx{K}$ is compact on $\cX_0$.
\end{lemma}
\begin{proof}
Using an argument similar to that in Lemma~\ref{lem:gCLM_symmetric_compact}, one can show that the map
\[(f,g) \mapsto \bigl(-\bar\Omega_x\mtx{V}(f),\, -\bar\Upsilon_x\mtx{V}(f)\bigr)\]
is compact on $\cX_0$. The remaining term in $\mtx{K}_{\rm pr}$ is a rank-one operator
\[ (f,g) \mapsto \mtx{H}(f)(0)\bigl(\bar\Omega-x\bar\Omega_x,\,\bar\Upsilon-x\bar\Upsilon_x\bigr).\]
It is not hard to see that $(\bar\Omega-x\bar\Omega_x,\,\bar\Upsilon-x\bar\Upsilon_x)\in \cX_0$. Furthermore, Lemma~\ref{lem:hl-coefficients} shows that the functional $f\mapsto\mtx{H}(f)(0)$ is bounded. This establishes the compactness of $\mtx{K}_{\rm pr}$ on $\cX_0$.

To show that the remainder term $(f,g)\mapsto -(0,\mtx{R}(f))$ is compact on $\cX_0$, it suffices to prove that 
\begin{equation}\label{eqt:HL_Rcompact}
h\mapsto\sqrt\phi\,\mtx{R}\left(\phi^{-1/2}h\right)
\end{equation}
is Hilbert--Schmidt on $L^2(\R)$. To this end, we write the kernel expansion of \eqref{eqt:HL_Rcompact} as
\[\sqrt{\phi(x)}\,\mtx{R}\left(\phi^{-1/2}h\right)(x)=\int_\R K_R(x,y)h(y)\idiff y,\]
where 
 \[ K_R(x,y):=\frac{\bar\Upsilon(x)}{\pi x}\frac{1}{y^2\sqrt{\phi(y)}}\frac{x^2\sqrt{\phi(x)}y-y^2\sqrt{\phi(y)}x}{x-y}.\]
By an argument similar to that in Lemma~\ref{lem:gCLM_symmetric_compact}, $K_R$ is bounded on $[-1,1]^2$ and satisfies the following estimates outside this square:
\[|K_R(x,y)|\lesssim
 \begin{cases}
 |x|^{\alpha/2-2/\bar c_l}\la y\ra^{-1-\alpha/2},& |x|\geq2|y|,\\
\la x\ra^{-2/\bar c_l}|y|^{-1},& |y|\geq2|x|,\\
 |x|^{-1-2/\bar c_l},& \frac12|x|<|y|<2|x|.
 \end{cases}\]
The third estimate can be combined with either of the first two: the first bound then holds on $|y|\leq|x|$, and the second on $|x|<|y|$. Squaring and integrating over these two regions gives
\begin{align*}
 \int_{\R^2}|K_R(x,y)|^2\idiff x\idiff y
 &\lesssim 1+\int_{|x|>\max\{1,|y|\}}|x|^{\alpha-4/\bar c_l}\la y\ra^{-2-\alpha}\idiff x\idiff y+\int_{|y|>\max\{1,|x|\}}\la x\ra^{-4/\bar c_l}|y|^{-2}\idiff x\idiff y\\
 &\lesssim 1+\int_{\R}\la x\ra ^{-1-4/\bar c_l}\idiff x<\infty.
\end{align*}
Here we used $-1<\alpha<2/\bar c_l-1$. Thus $h\mapsto\sqrt\phi\,\mtx{R}(\phi^{-1/2}h)$ is Hilbert--Schmidt on $L^2(\R)$, proving compactness of the remainder term. This completes the proof.
\end{proof}

With the above damping estimate and compactness result, we can establish finite-codimensional stability of \eqref{eqt:HL-linearized} via semigroups.
Adapting the proof of Lemma~\ref{lem:Lstb-semigroup} componentwise, we deduce that $\mtx{L}_{\rm stb}$, with domain
\[ \mathcal D(\mtx{L}_{\rm stb}):=\{F\in\cX_0:\bar{V} F_x\in\cX_0\},\]
generates a strongly continuous semigroup $\econst^{t\mtx{L}_{\rm stb}}$ on $\cX_0$. By the damping estimate \eqref{eqt:HL-energy-remainder}, we obtain the exponential decay of this semigroup:
\[ \left\|\econst^{t\mtx{L}_{\rm stb}}\right\|_{\cX_0\to\cX_0}\leq\econst^{-\eta_0t},\quad t\geq0,\]
where $\eta_0$ is given by Lemma \ref{lem:hl-weight}. 

Notice that $\mtx{L}$ is a bounded perturbation of $\mtx{L}_{\rm stb}$ by the compact operator $\mtx{K}$. By the bounded perturbation theorem, $\mtx{L}$ with domain 
\[\mathcal D(\mtx{L})=\mathcal D(\mtx{L}_{\rm stb})\]
 also generates a strongly continuous semigroup $\econst^{t\mtx{L}}$ on $\cX_0$. As in Subsection~\ref{sec:linear}, we fix $\lambda\in(0,\eta_0)$ and define the spectral projections $\mtx{\Pi}_\lambda$, $\mtx{\Pi}_s$ and the corresponding spaces $E_\lambda$, $E_s$ for $\mtx{L}$. In particular, $E_\lambda$ is finite dimensional, and
\[ \cX_0=E_\lambda\oplus E_s.\]
An argument similar to that in the proof of Lemma~\ref{lem:linear-package} gives the following decay estimate of the semigroup $\econst^{t\mtx{L}}$ restricted to $E_s$, which establishes finite-codimensional stability of \eqref{eqt:HL-linearized}.

\begin{lemma}\label{lem:hl-linear-package}
For every $\eta_s\in(0,\lambda)$, there is a constant $C_{\eta_s}>0$ such that
\begin{equation}
 \left\|\econst^{t\mtx{L}}\mtx{\Pi}_s\right\|_{\cX_0\to\cX_0}\leq C_{\eta_s}\econst^{-\eta_st},\quad t\geq0.
\end{equation}
\end{lemma}

\subsection{Nonlinear finite-codimensional stability}
In this subsection, we prove nonlinear finite-codimensional stability of \eqref{eqt:HL-perturbation} and establish asymptotically self-similar blowup of the HL model \eqref{eqt:1Dhouluo} from $C_c^{\infty}$ initial data. The proof is based on a shooting argument which is similar to that in Section~\ref{sec:noncompact}. To start, we perform nonlinear energy estimates for \eqref{eqt:HL-perturbation} in the following lemma.

\begin{lemma}\label{lem:hl-high-energy}
For the nonlinear operator $\mtx{N}$ defined in \eqref{eqt:HL-nonlinear-operator}, there exists a constant $C>0$ such that, for any $F\in\cX_1$,
\begin{equation}\label{eqt:HL-nonlinear-estimates}
 \|\mtx{N}(F)\|_{\cX_0}\leq C\|F\|_{\cX_0}\|F\|_{\cX_1},\quad\la\mtx{N}(F),F\ra_{\cX_1}\leq C\|F\|_{\cX_1}^3.
\end{equation}
Moreover, for every $\eta_1\in(0,\eta_0)$, a smooth solution of \eqref{eqt:HL-perturbation} satisfies
\begin{equation}\label{eqt:HL-high-energy}
 \frac12\frac{\idiff{}}{\idiff t}\|F\|_{\cX_1}^2 +\eta_1\|F\|_{\cX_1}^2 \leq C_{\eta_1}\|F\|_{\cX_0}^2+C_{\eta_1}\|F\|_{\cX_1}^3.
\end{equation}
Consequently,
\begin{equation}\label{eqt:HL-high-integral}
 \|F(t)\|_{\cX_1}\leq C\econst^{-\eta_1t}\|F(0)\|_{\cX_1}+C\int_0^t\econst^{-\eta_1(t-s)}\left(\|F(s)\|_{\cX_0}+\|F(s)\|_{\cX_1}^2\right)\idiff s.
\end{equation}
\end{lemma}
\begin{proof}
We first prove the nonlinear estimates. Lemma~\ref{lem:hl-coefficients} gives
\begin{align*}
 \|\mtx{N}(F)\|_{\cX_0}
 &\lesssim\left\|\frac{\mtx{V}(f)}{x}+\mtx{H}(f)(0)\right\|_{L^\infty(\R)}\|F\|_{\dot{\cX}_1}+\left(\|\mtx{H}(f)\|_{L^\infty}+|\mtx{H}(f)(0)|\right)\|F\|_{\cX_0}\\
 &\lesssim\|F\|_{\cX_0}\|F\|_{\cX_1}.
\end{align*}
Thus $|\la\mtx{N}(F),F\ra_{\cX_0}|\lesssim\|F\|_{\cX_0}^2\|F\|_{\cX_1}$. To estimate the $\dot{\cX}_1$ terms, integration by parts in \eqref{eqt:HL-nonlinear-operator} gives
\begin{align*}
 \la\mtx{N}(F),F\ra_{\dot{\cX}_1}
 &=\int_\R\left(\mtx{H}(f)(0)-\frac12\mtx{H}(f)+\frac{\mtx{V}(f)+x\mtx{H}(f)(0)}{2x}\frac{x\psi_x}{\psi}\right)f_x^2\psi\idiff x\\
 &\quad+M\int_\R\left(2\mtx{H}(f)(0)-\frac32\mtx{H}(f)+\frac{\mtx{V}(f)+x\mtx{H}(f)(0)}{2x}\frac{x\psi_x}{\psi}\right)g_x^2\psi\idiff x\\
 &\quad -M\int_\R gg_x\mtx{H}(f_x)\psi\idiff x.
\end{align*}
By Lemma~\ref{lem:hl-coefficients} and the boundedness of $x\psi_x/\psi$, the coefficients in the first two integrals are bounded in $L^\infty(\R)$ by $C\|F\|_{\cX_1}$. The argument for \eqref{eqt:gCLM-weighted-hilbert-derivative} in the proof of Lemma~\ref{lem:nonlinear_estimate} gives
\[\left\|\mtx{H}(f_x)\sqrt\psi\right\|_{L^2(\R)}\lesssim\left\|\frac{f_x}{x}\right\|_{L^2(\R)}+\left\|xf_x\la x\ra^{\alpha/2}\right\|_{L^2(\R)}\lesssim\left\|f_x\sqrt\psi\right\|_{L^2(\R)}.\]
Consequently,
\[M\left|\int_\R gg_x\mtx{H}(f_x)\psi \idiff x\right|\leq M\|g\|_{L^\infty(\R)}\left\|g_x\sqrt\psi\right\|_{L^2(\R)}\left\|\mtx{H}(f_x)\sqrt\psi\right\|_{L^2(\R)}\lesssim\|F\|_{\cX_1}^3,\]
which, together with the $\cX_0$ estimate, proves the cubic estimate in \eqref{eqt:HL-nonlinear-estimates}.

For the linear $L^2$ estimate, \eqref{eqt:HL-energy-remainder} and Lemma~\ref{lem:hl-compact} imply
\[\la\mtx{L}(F),F\ra_{\cX_0}=\la\mtx{L}_{\rm stb}(F),F\ra_{\cX_0}+\la\mtx{K}(F),F\ra_{\cX_0}\leq C\|F\|_{\cX_0}^2.\]
For the linear weighted $H^1$ estimate, we decompose the energy pairing as in the proof of Lemma~\ref{lem:high-energy}:
\begin{align*}
 \la\mtx{L}(F),F\ra_{\dot{\cX}_1}
 &=I_1+I_2+I_3,\\
 I_1&:=\int_\R\bigl(d_{f,1}f_x^2+Md_{g,1}g_x^2+f_xg_x\bigr)\psi\idiff x
       -M\int_\R\frac{x\bar\Upsilon}{\sqrt\phi}\mtx{H}\left(x\sqrt\phi f_x\right)g_x\phi\idiff x,\\
 I_2&:=-M\int_\R x\mtx{R}(xf_x)g_x\phi\idiff x,\\
 I_3&:=-\int_\R\mtx{C}_1(F)f_x\psi\idiff x-M\int_\R\mtx{C}_2(F)g_x\psi\idiff x,
\end{align*}
where
\begin{align*}
 \mtx{C}_1(F)&:=\bar\Omega_x\mtx{H}_0(f)+\bar\Omega_{xx}\bigl(\mtx{V}(f)+x\mtx{H}(f)(0)\bigr),\\
 \mtx{C}_2(F)&:=\bar U_{xx}g+2\bar\Upsilon_x\mtx{H}_0(f)+\bar\Upsilon_{xx}\bigl(\mtx{V}(f)+x\mtx{H}(f)(0)\bigr).
\end{align*}
As in Lemma~\ref{lem:high-energy}, $I_1$ is the damping term. Lemma~\ref{lem:hl-weight} with $j=1$ gives
\[ I_1\leq-\eta_0\|F\|_{\dot{\cX}_1}^2.\]
The term $I_2$ is the compact remainder. Lemmas~\ref{lem:hl-compact} and~\ref{lem:compactness-estimate} give, for every $\varepsilon>0$,
\[ |I_2|\leq\varepsilon\|F\|_{\cX_1}^2+C_\varepsilon\|F\|_{\cX_0}^2.\]
The term $I_3$ consists of cross terms. By the Cauchy--Schwarz inequality,
\[|I_3|\leq\left(\left\|\mtx{C}_1(F)\sqrt\psi\right\|_{L^2(\R)}+\sqrt M\left\|\mtx{C}_2(F)\sqrt\psi\right\|_{L^2(\R)}\right)\|F\|_{\dot{\cX}_1}.\]
As in Section~\ref{sec:noncompact}, \eqref{eqt:common-hilbert-identities}, \eqref{eqt:common-A2-weight}, and Lemmas~\ref{lem:hl-profile-derivatives} and~\ref{lem:hl-coefficients} yield
\[\left\|\mtx{C}_1(F)\sqrt\psi\right\|_{L^2(\R)}+\sqrt M\left\|\mtx{C}_2(F)\sqrt\psi\right\|_{L^2(\R)}\leq C\|F\|_{\cX_0}.\]
Young's inequality therefore gives
\[ |I_3|\leq\varepsilon\|F\|_{\dot{\cX}_1}^2+C_\varepsilon\|F\|_{\cX_0}^2.\]
Choosing $0<2\varepsilon<\eta_0-\eta_1$ and combining these estimates yields
\[\la\mtx{L}(F),F\ra_{\dot{\cX}_1}\leq-\eta_1\|F\|_{\dot{\cX}_1}^2+C_{\eta_1}\|F\|_{\cX_0}^2.\]
Combining the $\cX_0$ and $\dot{\cX}_1$ estimates and enlarging $C_{\eta_1}$ proves \eqref{eqt:HL-high-energy}. Finally, \eqref{eqt:HL-high-integral} follows from \eqref{eqt:HL-high-energy} and Grönwall's inequality. The proof is thus completed.
\end{proof}

As in Section~\ref{sec:noncompact}, let $\cX_c$ denote the space of pairs of odd $C_c^\infty(\R)$ functions supported away from the origin. We choose correction directions in this space to preserve the smoothness and compact support of $\Omega_0$ and $\Upsilon_0$. By an argument similar to that in the proof of Lemma~\ref{lem:unstable-adapted-norm}, we construct an inner product $\la\cdot,\cdot\ra_*$ on the present space $E_\lambda$, whose induced norm $\|\cdot\|_*$ is equivalent to $\|\cdot\|_{\cX_0}$ and satisfies
\[\la\mtx{L}(h),h\ra_*\geq-\lambda\|h\|_*^2,\quad h\in E_\lambda.\]
Similarly, adapting the proof of Lemma~\ref{lem:correction-directions}, we construct a real-linear correction map $\mtx{J}_\lambda:E_\lambda\to\cX_c$ satisfying
\[\mtx{\Pi}_\lambda\mtx{J}_\lambda(h)=h,\quad\|\mtx{J}_\lambda(h)\|_{\cX_1}\leq C\|h\|_*,\quad h\in E_\lambda.\]
Take $\eta=0.9\lambda$. Using arguments similar to those in the proofs of Lemma~\ref{lem:bootstrap} and Theorem~\ref{thm:brouwer} with Lemma~\ref{lem:hl-linear-package} replacing Lemma~\ref{lem:linear-package} and Lemma~\ref{lem:hl-high-energy} replacing Lemmas~\ref{lem:nonlinear_estimate} and~\ref{lem:high-energy}, we obtain finite-codimensional stability of \eqref{eqt:HL-perturbation} in the following theorem.

\begin{theorem}\label{thm:hl-shooting}
There exist \(\delta_0,\nu,\rho>0\) such that for any $\delta\in (0,\delta_0)$ and any initial data $\xi\in \cX_1$ with $\|\xi\|_{\cX_1}\leq \nu \delta$, there exists \(h\in E_{\lambda}\) with \(\|h\|_*\leq\rho\delta\), such that starting from the initial data
\[ F_h(0):=\xi+\mtx{J}_{\lambda}(h-\mtx{\Pi}_{\lambda}(\xi))\]
the solution $F_h(t)$ to \eqref{eqt:HL-perturbation} satisfies the decay estimates  
\begin{equation}
 \left\|F_h(t)\right\|_{\cX_1}\leq\delta \econst^{-\eta t},\end{equation}
 \begin{equation}
 \left\|\mtx{\Pi}_{\lambda}(F_h(t))\right\|_*\leq\rho\delta \econst^{-2\eta t},\quad t\geq0.
\end{equation}
\end{theorem}

Finally, we use an argument similar to that in the proof of Theorem~\ref{thm:gCLM-nonexpanding-blowup}, with Theorem~\ref{thm:hl-shooting} replacing Theorem~\ref{thm:brouwer} and \cite[Lemma~2.2]{chen2026novel} replacing \cite[Proposition~2.1]{huang2026singular}, to obtain the main result of this section, which establishes asymptotically self-similar finite-time blowup of the HL model \eqref{eqt:1Dhouluo} from smooth compactly supported initial data.

\begin{theorem}\label{thm:hl-blowup}
Let $(\bar\Omega,\bar\Upsilon,\bar c_l,-1)$ be a steady-state solution of \eqref{eqt:HL-rescaled-xt} given by Theorem~\ref{thm:hl-existence}. There exist odd initial data $\Omega_0,\Upsilon_0\in C_c^\infty(\R)$ such that the solution $(\Omega(\cdot,t),\Upsilon(\cdot,t))$ to \eqref{eqt:HL-rescaled-xt} with the normalization \eqref{eqt:HL-normalization} and initial data $(\Omega_0,\Upsilon_0)$ satisfies
\[\left\|\left(\Omega(\cdot,t)-\bar\Omega,\Upsilon(\cdot,t)-\bar\Upsilon\right)\right\|_{\cX_1}\lesssim\econst^{-\eta t}.\]
As a result, starting from the same initial data $(\Omega_0,\Upsilon_0)$, the solution $(\omega,\upsilon)$ of the original HL model \eqref{eqt:HL-model} develops an asymptotically self-similar singularity in finite time.
\end{theorem}

\appendix 
\section{Auxiliary lemmas}
We gather some auxiliary lemmas in this appendix.
\label{app:hilbert}
\begin{lemma}\label{lem:core-identities}
For any suitable function $f$,
\begin{equation}\label{eqt:common-hilbert-identities}
 \frac{\mtx{H}(f)(x)-\mtx{H}(f)(0)}x=\mtx{H}\left(\frac{f}{x}\right)(x),
\quad \mtx{H}(xf_x)=x\mtx{H}(f_x).
\end{equation}
Moreover, suppose $f$ is odd and $f_x(0)=0$, then
\begin{equation}\label{eqt:common-hilbert-identities2}
 \frac{\mtx{H}(f)(x)-\mtx{H}(f)(0)}{x^2}=\mtx{H}\left(\frac{f}{x^2}\right)(x).
\end{equation}
\end{lemma}
\begin{proof}
For smooth $f$ with sufficient decay and $x\neq0$, we have
\[\frac{\mtx{H}(f)(x)-\mtx{H}(f)(0)}x=\frac1{\pi x}\pv\int_\R f(y)\left(\frac1{x-y}+\frac1y\right)\idiff y=\frac1\pi\pv\int_\R\frac{f(y)}{y(x-y)}\idiff y=\mtx{H}\left(\frac f x\right)(x).\]
Applying this identity to $xf_x$ proves the second identity in \eqref{eqt:common-hilbert-identities}, since 
\[\mtx{H}(xf_x)(0)=-\frac{1}{\pi}\int_\R f_x(y)\idiff y=0.\] 
If $f$ is odd and $f_x(0)=0$, then $f/x$ is even and vanishes at the origin, so $\mtx{H}(f/x)(0)=0$. Thus,
\[\frac{\mtx{H}(f)(x)-\mtx{H}(f)(0)}{x^2}=\frac{\mtx{H}\left(f/x\right)(x)-\mtx{H}(f/x)(0)}{x}=\mtx{H}\left(\frac{f}{x^2}\right)(x), \]
which proves \eqref{eqt:common-hilbert-identities2}.
\end{proof}
\begin{lemma}
For any $\alpha\in(-1,1)$, there exists a constant $C=C_{\alpha}$ such that
\begin{equation}\label{eqt:common-A2-weight}
      \int_{\R} |\mtx{H}(f)|^2\la x\ra^{\alpha} \idiff x\leq C \int_{\R} |f|^2\la x\ra^{\alpha} \idiff x.
\end{equation}
\end{lemma}
\begin{proof}
      The above inequality follows from the fact that $\la x\ra^{\alpha}$ is an $A_2$ weight.
\end{proof}

\begin{lemma}[Compactness estimate]\label{lem:compactness-estimate}
Let $\mathcal{A}$ and $\mathcal{B}$ be Hilbert spaces such that $\mathcal{A}$ is continuously embedded in $\mathcal{B}$. If $\mtx{K}:\mathcal{A}\to\mathcal{B}$ is a compact linear operator, then, for every $\varepsilon>0$, there exists $C_\varepsilon>0$ such that
\[
 \|\mtx{K}x\|_{\mathcal{B}}\leq\varepsilon\|x\|_{\mathcal{A}}+C_\varepsilon\|x\|_{\mathcal{B}},\quad x\in\mathcal{A}.
\]
\end{lemma}
\begin{proof}
If the estimate fails for some $\varepsilon>0$, there exists a sequence $\{x_n\}\subset\mathcal{A}$ such that
\[
 \|x_n\|_{\mathcal{A}}=1,\quad \|\mtx{K}x_n\|_{\mathcal{B}}>\varepsilon+n\|x_n\|_{\mathcal{B}}.
\]
Since $\mtx{K}$ is bounded, $x_n\to0$ in $\mathcal{B}$. Passing to a subsequence, $x_n\rightharpoonup x$ in $\mathcal{A}$, and hence in $\mathcal{B}$, so $x=0$. Compactness then implies $\mtx{K}x_n\to0$ in $\mathcal{B}$, a contradiction.
\end{proof}

\section{Proofs of the profile properties}
In this appendix, we complete the proofs of the profile properties stated in Sections~\ref{sec:rescaling} and~\ref{sec:hl-rescaling}.
\label{app:profile-estimates}
\subsection{The gCLM model}\label{app:gclm-profiles}
Throughout this subsection, we write $\bar V=a\bar U+\bar c_lX$ as in Sections \ref{sec:noncompact} and \ref{sec:compact}. The steady equation of the dynamic rescaling equation is given by 
\begin{equation}\label{eqt:gCLM-steady-equation}
 \bar{V}\bar\Om_X=(\bar{U}_X-1)\bar\Om.
\end{equation}

\begin{proof}[Proof of Theorem~\ref{thm:out-going}]
The fixed-point representation in
\cite[Section~2.2 and Proposition~3.1]{huang2024self}, together with
\cite[Lemma~3.4 and the proof of Corollary~3.5]{huang2024self}, gives
the asserted signs of $\bar{V}$, with $\bar{V}_X(0)>0$ and, in the compactly
supported case, $\bar{V}_X(L)<0$. For $a=0$, these assertions follow from
$\bar{V}(X)=X$.
Dividing \eqref{eqt:gCLM-steady-equation} by $X$ and letting $X\to0$ gives
$\bar{V}_X(0)=(\bar{U}_X(0)-1)>0$, since $\bar\Om_X(0)\neq0$.

For $\bar c_l<0$, $\bar U_X\in H^1(\R)$, so $\bar U_X-1$ is continuous. Dividing \eqref{eqt:gCLM-steady-equation} by $(L-X)^p$ and letting $X\to L^{-}$, we obtain
\[(\bar{U}_X(L)-1)=p\bar{V}_X(L)\leq \bar{V}_X(L)<0,\]
which proves \eqref{eqt:gCLM-npc-BL}.
\end{proof}

\begin{lemma}\label{lem:boundedness-and-decay}
Suppose that $f\in L^r(\R)\cap \dot{H}^1(\R)$ for some $1<r<\infty$. Then 
\[ \|f\|_{L^\infty(\R)} \leq\|f\|_{L^r(\R)}+\|f_X\|_{L^2(\R)}, \quad \lim_{|X|\to\infty}f(X)=0. \]
\end{lemma}
\begin{proof}
For any $X\in\R$ and $Y\in[X,X+1]$, we have 
\[|f(X)| \leq|f(Y)|+ \int_X^{X+1}|f_X(Y)|\idiff Y. \]
Integrating over $[X,X+1]$ and applying H\"older's inequality, we obtain
\[ |f(X)| \leq \|f\|_{L^1([X,X+1])}+ \|f_X\|_{L^1([X,X+1])} 
 \leq \|f\|_{L^r([X,X+1])}+ \|f_X\|_{L^2([X,X+1])}.\]
This proves the stated bound. Note that both terms on the right tend to zero as $|X|\to\infty$, which implies $f$ vanishes at infinity.
\end{proof}

\begin{proof}[Proof of Lemma~\ref{lem:profile-derivative-estimates}]
Choose $r>\max\{2,\bar c_l\}$. Theorem~\ref{thm:existence_of_self_similar_profile_gclm} and the $L^r$-boundedness of $\mtx{H}$ give $\bar\Om,\bar U_X\in L^r(\R)$, while $\bar U_{XX}=\mtx{H}(\bar\Om_X)\in L^2(\R)$. By Lemma~\ref{lem:boundedness-and-decay}, we have $\bar U_X\in L^\infty(\R)$ and $\bar U_X(X)\to0$ as $|X|\to\infty$, which implies $\bar U(X)/X\to0\, (X\to +\infty)$, and $(\bar{U}_X-1),\bar{V}_X\in L^\infty(\R)$.

\emph{The case $\bar c_l>0$:}
In this case we have $\bar{V}(X)\asymp\bar c_lX$ as $X\to +\infty$. By Theorem~\ref{thm:existence_of_self_similar_profile_gclm} we know that $|\bar{\Om}(X)|\lesssim |X|^{-1/\bar c_l}$. Therefore,
\[\left|\bar\Om_X(X)\right|=\left|\frac{(\bar{U}_X(X)-1)\bar\Om(X)}{\bar{V}(X)}\right|\lesssim X^{-1-1/\bar c_l},\quad X\geq1.\]
Differentiating \eqref{eqt:gCLM-steady-equation} yields
\begin{equation}\label{eqt:gCLM-app-profile-differentiated}
 \bar{V}\bar\Om_{XX}=\bar U_{XX}\bar\Om+((\bar{U}_X-1)-\bar{V}_X)\bar\Om_X,
\end{equation}
and thus 
\[ X\bar\Om_{XX} =\frac{X\bar\Om}{\bar{V}}\bar U_{XX} +\frac{X((\bar{U}_X-1)-\bar{V}_X)}{\bar{V}}\bar\Om_X\in L^2(\R).\]
We therefore have $X\bar{\Om}_X\in L^r(\R)\cap \dot{H}^1(\R)$. By Lemma \ref{lem:core-identities} and the boundedness of the Hilbert transform, we deduce
\[ X\bar U_{XX}=X\mtx{H}(\bar{\Om}_X)=\mtx{H}(X\bar\Om_X)\in L^r(\R)\cap \dot{H}^1(\R).\]
Hence $X\bar U_{XX}\in L^{\infty}(\R)$. Substitution into \eqref{eqt:gCLM-app-profile-differentiated} yields $|\bar\Om_{XX}(X)|\lesssim X^{-2-1/\bar c_l}$, which completes the proof for $\bar{c}_l>0$.

\emph{The case $\bar c_l=0$:}
By Theorem~\ref{thm:existence_of_self_similar_profile_gclm} we have $|\bar\Om(X)|\lesssim\econst^{-\delta X^2}$ for some $\delta>0$. It is not hard to show that 
\[\bar U(X)=\frac1\pi\int_0^\infty\bar\Om(Y)\log\left|\frac{X-Y}{X+Y}\right|\idiff Y=-\frac{1}{\pi X}\int_\R Y\bar\Om(Y)\idiff Y+O\left(X^{-3}\right)\asymp X^{-1},\quad X\to +\infty.\]
Here the leading coefficient is positive since $Y\bar\Om(Y)<0$ for $Y\neq0$. Theorem~\ref{thm:out-going} gives $\bar{V}=a\bar U>0$, so $a>0$ and $\bar{V}\asymp X^{-1}$ as $X\to +\infty$. We can use a similar argument as in the case $\bar c_l>0$ to deduce that 
\[\left|\bar\Om_X(X)\right|=\left|\frac{(\bar{U}_X(X)-1)\bar\Om(X)}{\bar{V}(X)}\right|\lesssim X\econst^{-\delta X^2},\quad X\geq 1,\] and similarly 
\[X\bar{U}_{XX}\in L^{\infty}(\R),\quad \left|\bar\Om_{XX}(X)\right| \lesssim X^2\econst^{-\delta X^2}, \quad X\geq 1.\]

\emph{The case $\bar c_l<0$:}
By Theorem~\ref{thm:out-going} we have $\bar{V}(X)\asymp L-X$ near $X=L$. Hence
\[ \bar\Om_X(X)=\frac{(\bar{U}_X(X)-1)\bar\Om(X)}{\bar{V}(X)} =O\left((L-X)^{p-1}\right)\quad\text{as }X\to L^{-}, \]
which implies $\bar\Om_X\in L^\infty(\R)$ since $\bar\Om_X(X)$ is interiorly smooth and $p\geq 1$.
For the case $a<1$, we have $p>1$. By \eqref{eqt:gCLM-app-profile-differentiated},
\[\left|\bar\Om_{XX}(X)\right|\lesssim(L-X)^{p-1}\left|\bar U_{XX}(X)\right|+(L-X)^{p-2},\]
and thus
 \[\int_{L/2}^L(L-X)\left|\bar\Om_{XX}(X)\right|^2\idiff X
 \lesssim\left\|\bar U_{XX}\right\|_2^2 +\int_{L/2}^L(L-X)^{2p-3}\idiff X<\infty.\]
If $a=1$, the identity at the origin in Theorem~\ref{thm:out-going} gives $\bar c_l=-1$, hence $\bar{V}_X=(\bar{U}_X-1)$. Thus \eqref{eqt:gCLM-app-profile-differentiated} reduces to $\bar\Om_{XX}=(\bar\Om/\bar{V})\bar U_{XX}\in L^2(L/2,L)$. This proves \eqref{eqt:gCLM-npc-edge-H2}. Finally, since $\bar{\Om}$ is compactly supported, for any $X>L$, we have 
\[\bar{U}_X(X)=\frac1\pi\int_0^{L}\frac{2Y\bar{\Om}(Y)}{X^2-Y^2} \idiff Y= O\left(X^{-2}\right), \quad \bar{U}_{XX}(X)=\frac1\pi\int_0^{L}\frac{-4XY\bar{\Om}(Y)}{\left(X^2-Y^2\right)^2} \idiff Y= O\left(X^{-3}\right) \]
as $X\to \infty$. This completes the proof of Lemma~\ref{lem:profile-derivative-estimates}.
\end{proof}

\subsection{The HL model}\label{app:hl-profiles}
Throughout this subsection, we write $\bar V=\bar U+\bar c_lX$, as in Section~\ref{sec:hl-stability}. The steady equations \eqref{eqt:HL-steady}, with $\bar c_\om=-1$, are
\begin{equation}\label{eqt:HL-app-steady-equations}
 \bar V\bar\Omega_X=-\bar\Omega+\bar\Upsilon,\quad \bar V\bar\Upsilon_X=-(2+\bar U_X)\bar\Upsilon.
\end{equation}

\begin{proof}[Proof of Theorem~\ref{thm:hl-out-going}]
By \cite[Corollary~3.6]{huang2025exact}, we have $\bar{V}(X)/(X\bar{V}_X(0))\geq1$ for $X>0$. Differentiating the two equations in \eqref{eqt:HL-app-steady-equations} at the origin gives
\[ (\bar c_l+\bar U_X(0))\bar\Omega_X(0) =-\bar\Omega_X(0)+\bar\Upsilon_X(0),\quad \bar c_l+\bar U_X(0)=-2-\bar U_X(0).\]
These yield \eqref{eqt:HL-origin-identities} and the asserted lower bound. 
\end{proof}

\begin{proof}[Proof of Lemma~\ref{lem:hl-profile-derivatives}]
The argument is similar to the proof of Lemma~\ref{lem:profile-derivative-estimates} in Appendix~\ref{app:gclm-profiles}. Choose $q>\max\{2,\bar c_l\}$. The tail estimates in Theorem~\ref{thm:hl-existence} and the $L^q$-boundedness of $\mtx{H}$ give $\bar U_X\in L^q(\R)$, while $\bar U_{XX}=\mtx{H}(\bar\Omega_X)\in L^2(\R)$. Lemma~\ref{lem:boundedness-and-decay} gives $\bar U_X\in L^\infty(\R)$ and $\bar U_X(X)\to0$, hence $\bar U(X)/X\to0$. In particular, $\bar{V}\asymp X$ and $\bar{V}_X$ is bounded.

The steady equations now give the first derivative estimates in \eqref{eqt:HL-profile-derivatives}. Differentiating their first equation yields
\[
 \bar{V}\bar\Omega_{XX}=-(1+\bar{V}_X)\bar\Omega_X+\bar\Upsilon_X.
\]
Since $X/\bar{V}$ is bounded, this shows $X\bar\Omega_{XX}\in L^2(\R)$ and gives the pointwise bound for $\bar\Omega_{XX}$. Thus $X\bar\Omega_X\in L^q(\R)\cap\dot H^1(\R)$. By Lemma~\ref{lem:core-identities}, we have
\[ X\bar U_{XX}=\mtx{H}(X\bar\Omega_X)\in L^q(\R)\cap\dot H^1(\R).\]
Another application of Lemma~\ref{lem:boundedness-and-decay} yields $\|X\bar U_{XX}\|_\infty<\infty$. Finally, differentiating the second equation, we have
\[ \bar{V}\bar\Upsilon_{XX}=-X\bar U_{XX}\frac{\bar\Upsilon}{X}+(-2-\bar U_X-\bar{V}_X)\bar\Upsilon_X,\]
which gives the pointwise bound for $\bar\Upsilon_{XX}$. This completes the proof of Lemma~\ref{lem:hl-profile-derivatives}.
\end{proof}

\vspace{2mm}

\subsection*{Acknowledgement} The authors are supported by the National Key R\&D Program of China under the grant 2021YFA1001500 and the National Natural Science Foundation of China under the grant NSFC No. 12288101.

\bibliographystyle{myalpha}
\bibliography{reference}

\end{document}

%% file: macro-3.tex
\usepackage[usenames,dvipsnames]{xcolor}
\usepackage{fancyhdr}
\usepackage{amsmath,amsfonts,amsbsy,amsgen,amscd,amssymb,amsthm}
\usepackage{url}

\usepackage{mathtools}

\usepackage[font=small,margin=0.25in,labelfont={bf},labelsep={colon}]{caption}

\usepackage{subcaption}

\usepackage{tikz}
\usepackage{microtype}
\usepackage{enumitem}

\definecolor{dark-gray}{gray}{0.3}
\definecolor{dkgray}{rgb}{.4,.4,.4}
\definecolor{dkblue}{rgb}{0,0,.5}
\definecolor{medblue}{rgb}{0,0,.75}
\definecolor{rust}{rgb}{0.5,0.1,0.1}

\usepackage[colorlinks=true]{hyperref}

\hypersetup{urlcolor=rust}
\hypersetup{citecolor=dkblue}
\hypersetup{linkcolor=dkblue}

\usepackage{setspace}

\usepackage{graphicx}
\usepackage{booktabs,longtable,tabu} 
\usepackage{multirow} 

\usepackage{float}

\usepackage[full]{textcomp}

\usepackage[scaled=.98,sups]{XCharter}
\usepackage[scaled=1.04,varqu,varl]{inconsolata}
\usepackage[type1]{cabin}
\usepackage[charter,vvarbb,scaled=1.07]{newtxmath}
\usepackage[cal=boondoxo]{mathalfa}
\usepackage[T1]{fontenc}

\usepackage{bm} 

\graphicspath{{figures/}}

\newtheorem{theorem}{Theorem}[section]
\newtheorem{lemma}[theorem]{Lemma}

\newtheorem{proposition}[theorem]{Proposition}

\theoremstyle{definition}

\newtheorem{remark}[theorem]{Remark}

\numberwithin{equation}{section} 
\numberwithin{figure}{section}
\numberwithin{table}{section}

\floatstyle{plaintop}
\newfloat{recipe}{thp}{lor}
\floatname{recipe}{Recipe}
\numberwithin{recipe}{section}

\providecommand{\mathbold}[1]{\bm{#1}}  
                                
\renewcommand{\phi}{\varphi}

\newcommand{\econst}{\mathrm{e}}
\newcommand{\iunit}{\mathrm{i}}

\providecommand{\mathbbm}{\mathbb} 

\newcommand{\R}{\mathbbm{R}}

\newcommand{\CC}{\mathbbm{C}}

\newcommand{\diff}[1]{\mathrm{d}{#1}}
\newcommand{\idiff}[1]{\, \diff{#1}}

\newcommand{\mtx}[1]{\mathbold{#1}}

\newcommand{\triplenorm}[1]{{\left\vert\kern-0.25ex\left\vert\kern-0.25ex\left\vert #1
    \right\vert\kern-0.25ex\right\vert\kern-0.25ex\right\vert}}



%% file: reference.bib
@article{luo2014potentially,
  title={{Potentially singular solutions of the 3D axisymmetric Euler equations}},
  author={Luo, Guo and Hou, Thomas Y},
  journal={Proceedings of the National Academy of Sciences},
  volume={111},
  number={36},
  pages={12968--12973},
  year={2014},
  publisher={National Academy of Sciences}
}

@article{luo2014toward,
  title={{Toward the finite-time blowup of the 3D axisymmetric Euler equations: a numerical investigation}},
  author={Luo, Guo and Hou, Thomas Y},
  journal={Multiscale Modeling \& Simulation},
  volume={12},
  number={4},
  pages={1722--1776},
  year={2014},
  publisher={SIAM}
}

@article{choi2015finite,
  title={{Finite time blow up for a 1D model of 2D Boussinesq system}},
  author={Choi, Kyudong and Kiselev, Alexander and Yao, Yao},
  journal={Communications in Mathematical Physics},
  volume={334},
  number={3},
  pages={1667--1679},
  year={2015},
  publisher={Springer}
}

@article{choi2017finite,
  title={{On the finite-time blowup of a one-dimensional model for the three-dimensional axisymmetric Euler equations}},
  author={Choi, Kyudong and Hou, Thomas Y and Kiselev, Alexander and Luo, Guo and {\v{S}}ver{\'a}k, Vladimir and Yao, Yao},
  journal={Communications on Pure and Applied Mathematics},
  volume={70},
  number={11},
  pages={2218--2243},
  year={2017},
  publisher={Wiley Online Library}
}

@article{chen2021finite,
  title={On the finite time blowup of the {D}e {G}regorio model for the 3{D} {E}uler equations},
  author={Chen, Jiajie and Hou, Thomas Y and Huang, De},
  journal={Communications on Pure and Applied Mathematics},
  volume={74},
  number={6},
  pages={1282--1350},
  year={2021},
  publisher={Wiley Online Library}
}

@article{elgindi2021finite,
  title={{Finite-time singularity formation for ${C}^{1,\alpha}$ solutions to the incompressible {E}uler equations on $\mathbb{R}^{3}$}},
  author={Elgindi, Tarek M},
  journal={Annals of Mathematics},
  volume={194},
  number={3},
  pages={647--727},
  year={2021},
  publisher={Department of Mathematics, Princeton University Princeton, New Jersey, USA}
}

@article{elgindi2022stability,
  title={On the stability of self-similar blow-up for ${C}^{1,\alpha}$ solutions to the incompressible {E}uler equations on $\mathbb{R}^{3}$},
  author={Elgindi, Tarek M and Ghoul, Tej-Eddine and Masmoudi, Nader},
  journal={Cambridge Journal of Mathematics},
  volume={9},
  number={4},
  pages={1035--1075},
  year={2021},
  publisher={International Press of Boston, Inc. Somerville, MA 02143, USA}
}

@article{chen2022stable,
  title={{Stable nearly self-similar blowup of the 2D Boussinesq and 3D Euler equations with smooth data I: Analysis}},
  author={Chen, Jiajie and Hou, Thomas Y},
  journal={arXiv preprint arXiv:2210.07191},
  year={2022}
}

@article{chen2025singularity,
  title={Singularity formation in 3{D} {Euler} equations with smooth initial data and boundary},
  author={Chen, Jiajie and Hou, Thomas Y},
  journal={Proceedings of the National Academy of Sciences},
  volume={122},
  number={27},
  pages={e2500940122},
  year={2025},
  publisher={National Academy of Sciences}
}

@article{chen2022asymptotically,
  title={{Asymptotically self-similar blowup of the Hou--Luo model for the 3D Euler equations}},
  author={Chen, Jiajie and Hou, Thomas Y and Huang, De},
  journal={Annals of PDE},
  volume={8},
  number={2},
  pages={24},
  year={2022},
  publisher={Springer}
}

@article{chen2026eulerII,
  author={Chen, Jiajie},
  title={{Asymptotically self-similar blowup for 3D incompressible {Euler}
          with {$C^{1,1/3-}$} velocity II: 3D profiles, blowup, and limiting
          behavior}},
  journal={arXiv preprint arXiv:2605.15130},
  year={2026},
  eprint={2605.15130},
  archivePrefix={arXiv},
  primaryClass={math.AP}
}

@article{shao2026selfsimilar,
  author={Shao, Feng and Wei, Dongyi and Zhang, Ping and Zhang, Zhifei},
  title={{Self-similar blow-up solutions of $d$-dimensional incompressible {Euler} equations with $C^{1,(1-2/d)-}$ velocity}},
  journal={arXiv preprint arXiv:2605.19716},
  year={2026},
  eprint={2605.19716},
  archivePrefix={arXiv},
  primaryClass={math.AP}
}

@article{huang2024self,
  title={{Self-similar finite-time blowups with smooth profiles of the generalized {Constantin--Lax--Majda} model}},
  author={Huang, De and Qin, Xiang and Wang, Xiuyuan and Wei, Dongyi},
  journal={Archive for Rational Mechanics and Analysis},
  volume={248},
  number={2},
  pages={22},
  year={2024},
  publisher={Springer}
}

@article{huang2025exact,
  title={{Exact self-similar finite-time blowup of the Hou--Luo model with smooth profiles}},
  author={Huang, De and Qin, Xiang and Wang, Xiuyuan and Wei, Dongyi},
  journal={Communications in Mathematical Physics},
  volume={406},
  number={10},
  pages={243},
  year={2025},
  publisher={Springer}
}

@article{cordoba2025finite,
  title={Finite Time Singularities to the 3{D} Incompressible {E}uler Equations for Solutions in \({C}^{\infty}(\mathbb{R}^3 \setminus \{0\})\cap {C}^{1,\alpha}\cap {L}^2\)},
  author={C{\'o}rdoba, Diego and Mart{\'i}nez-Zoroa, Luis and Zheng, Fan},
  journal={Annals of PDE},
  volume={11},
  number={2},
  pages={19},
  year={2025},
  publisher={Springer}
}

@book{majda2002vorticity,
  title={{Vorticity and Incompressible Flow}},
  author={Majda, Andrew J and Bertozzi, Andrea L},
  series={Cambridge Texts in Applied Mathematics},
  publisher={Cambridge University Press},
  year={2002}
}

@article{elgindi2021stable,
  author={Elgindi, Tarek M. and Ghoul, Tej-Eddine and Masmoudi, Nader},
  title={{Stable self-similar blowup for a family of nonlocal transport equations}},
  journal={Analysis \& PDE},
  volume={14},
  number={3},
  pages={891--908},
  year={2021},
  publisher={Mathematical Sciences Publishers}
}

@book{engel2000one,
  author={Engel, Klaus-Jochen and Nagel, Rainer},
  title={{One-Parameter Semigroups for Linear Evolution Equations}},
  series={Graduate Texts in Mathematics},
  volume={194},
  publisher={Springer-Verlag},
  address={New York},
  year={2000}
}

@article{chen2020singularity,
  title={Singularity formation and global well-posedness for the generalized {C}onstantin--{L}ax--{M}ajda equation with dissipation},
  author={Chen, Jiajie},
  journal={Nonlinearity},
  volume={33},
  number={5},
  pages={2502--2532},
  year={2020},
  doi={10.1088/1361-6544/ab74b0},
  publisher={IOP Publishing}
}

@article{constantin1985simple,
  title={A simple one-dimensional model for the three-dimensional vorticity equation},
  author={Constantin, Peter and Lax, Peter D and Majda, Andrew},
  journal={Communications on Pure and Applied Mathematics},
  volume={38},
  number={6},
  pages={715--724},
  year={1985},
  publisher={Wiley Online Library}
}

@article{elgindi2020effects,
  title={On the effects of advection and vortex stretching},
  author={Elgindi, Tarek M and Jeong, In-Jee},
  journal={Archive for Rational Mechanics and Analysis},
  volume={235},
  number={3},
  pages={1763--1817},
  year={2020},
  publisher={Springer}
}

@article{lushnikov2021collapse,
  title={Collapse Versus Blow-Up and Global Existence in the Generalized {C}onstantin--{L}ax--{M}ajda Equation},
  author={Lushnikov, Pavel M and Silantyev, Denis A and Siegel, Michael},
  journal={Journal of Nonlinear Science},
  volume={31},
  number={5},
  pages={82},
  year={2021},
  doi={10.1007/s00332-021-09737-x},
  publisher={Springer}
}

@article{de1990one,
  title={On a one-dimensional model for the three-dimensional vorticity equation},
  author={De Gregorio, Salvatore},
  journal={Journal of Statistical Physics},
  volume={59},
  number={5},
  pages={1251--1263},
  year={1990},
  publisher={Springer}
}

@article{okamoto2008generalization,
  title={On a generalization of the {C}onstantin--{L}ax--{M}ajda equation},
  author={Okamoto, Hisashi and Sakajo, Takashi and Wunsch, Marcus},
  journal={Nonlinearity},
  volume={21},
  number={10},
  pages={2447--2461},
  year={2008},
  doi={10.1088/0951-7715/21/10/013},
  publisher={IOP Publishing}
}

@article{chen2021slightly,
  title={On the Slightly Perturbed {D}e {G}regorio Model on $ {S}^1$ },
  author={Chen, Jiajie},
  journal={Archive for Rational Mechanics and Analysis},
  volume={241},
  number={3},
  pages={1843--1869},
  year={2021},
  publisher={Springer}
}

@article{chen2021regularity,
  title={On the regularity of the {D}e {G}regorio model for the 3{D} {E}uler equations},
  author={Chen, Jiajie},
  journal={Journal of the European Mathematical Society},
  volume={27},
  number={4},
  pages={1619--1677},
  year={2025},
  doi={10.4171/JEMS/1399},
  eprint={2107.04777},
  archivePrefix={arXiv},
  primaryClass={math.AP}
}

@article{huang2023self,
  title={On self-similar finite-time blowups of the {D}e {G}regorio model on the real line},
  author={Huang, De and Tong, Jiajun and Wei, Dongyi},
  journal={Communications in Mathematical Physics},
  volume={402},
  number={3},
  pages={2791--2829},
  year={2023},
  publisher={Springer}
}

@article{castro2009self,
  title={Self-similar solutions for a transport equation with nonlocal flux},
  author={Castro, Angel and C{\'o}rdoba, Diego},
  journal={Chinese Annals of Mathematics, Series B},
  volume={30},
  number={5},
  pages={505--512},
  year={2009},
  publisher={Springer}
}

@article{shkoller2026threshold,
  author  = {Shkoller, Steve},
  title   = {{Incompressible Euler blowup at the $C^{1,\frac{1}{3}}$ threshold}},
  journal = {arXiv preprint arXiv:2603.10945},
  year    = {2026},
  doi     = {10.48550/arXiv.2603.10945},
  url     = {https://arxiv.org/abs/2603.10945}
}

@article{shao2026global,
  author  = {Shao, Feng and Wei, Dongyi and Zhang, Zhifei},
  title   = {Global Regularity of Axisymmetric {E}uler Equations Without Swirl in Higher Dimensions},
  journal = {Acta Mathematica Sinica, English Series},
  volume  = {42},
  number  = {3},
  pages   = {663--679},
  year    = {2026},
  doi     = {10.1007/s10114-026-5110-0},
  url     = {https://doi.org/10.1007/s10114-026-5110-0}
}

@article{cordoba2023forced,
  author  = {C{\'o}rdoba, Diego and Mart{\'i}nez-Zoroa, Luis},
  title   = {{Blow-up for the incompressible 3D-Euler equations with uniform $C^{1,\frac{1}{2}-\epsilon}\cap L^2$ force}},
  journal = {arXiv preprint arXiv:2309.08495},
  year    = {2023},
  doi     = {10.48550/arXiv.2309.08495},
  url     = {https://arxiv.org/abs/2309.08495}
}

@misc{openai2026euler,
  author       = {{OpenAI}},
  title        = {Finite time blowup for the {E}uler equation},
  howpublished = {Preprint},
  year         = {2026},
  note         = {\url{https://cdn.openai.com/pdf/315b36cd-ec98-4023-8342-93345194ece1/euler.pdf}},
  url          = {https://cdn.openai.com/pdf/315b36cd-ec98-4023-8342-93345194ece1/euler.pdf}
}

@article{huang2026singular,
  title={Self-similar finite-time blowups with singular profiles of the generalized {Constantin--Lax--Majda} model: theoretical and numerical investigations},
  author={Huang, De and Tong, Jiajun and Wang, Xiuyuan},
  journal={Nonlinearity},
  volume={39},
  number={9},
  pages={095005},
  year={2026}
}

@article{chen2026novel,
  title={Novel self-similar finite-time blowups with singular profiles of the {1D Hou--Luo} model and the {2D Boussinesq} equations: A numerical investigation},
  author={Chen, Bojin and Huang, De and Li, Xiangyuan},
  journal={arXiv preprint arXiv:2604.01868},
  year={2026},
  doi={10.48550/arXiv.2604.01868},
  url={https://arxiv.org/abs/2604.01868}
}

@article{chen2021boundary,
  title={{Finite time blowup of 2D Boussinesq and 3D Euler equations with $C^{1,\alpha}$ velocity and boundary}},
  author={Chen, Jiajie and Hou, Thomas Y.},
  journal={Communications in Mathematical Physics},
  volume={383},
  number={3},
  pages={1559--1667},
  year={2021},
  doi={10.1007/s00220-021-04067-1},
  eprint={1910.00173},
  archivePrefix={arXiv},
  primaryClass={math.AP}
}

@article{hou2024blowup,
  title={{Blowup analysis for a quasi-exact 1D model of 3D Euler and Navier--Stokes}},
  author={Hou, Thomas Y. and Wang, Yixuan},
  journal={Nonlinearity},
  volume={37},
  number={3},
  pages={035001},
  year={2024},
  doi={10.1088/1361-6544/ad1c2f},
  eprint={2306.04146},
  archivePrefix={arXiv},
  primaryClass={math.AP}
}
